\documentclass[preprint,12pt]{elsarticle}

\usepackage{amssymb}
\usepackage{amsthm}
\usepackage{amsmath}
\usepackage{multirow}
\usepackage{subfigure}
\usepackage{float}

\journal{Journal of Computational Physics}

\begin{document}

\begin{frontmatter}



\title{Virtual Interpolation Point Method for Viscous Flows in Complex Geometries}

\author[a,b]{Seong-Kwan Park\corref{park}}
\ead{pakskwan@yonsei.kr}
\author[b,c]{Gahyung Jo}
\author[b]{Hi Jun Choe}
\address[a]{Department of Turbulent Boundary Layer, PARK Seong-Kwan Institute, Seoul, 136-858, Republic of Korea}
\address[b]{Department of Mathematics, Yonsei University, Seoul, 120-749, Republic of Korea}
\address[c]{National Fusion Research Institute, Daejeon 169-148, Republic of Korea}
\cortext[park]{Corresponding author}

\begin{abstract}
A new approach for simulating flows over complex geometries is
developed by introducing an accurate virtual interpolation point
scheme as well as a virtual local stencil approach. The present
method is based on the concept of point collocation on a virtual
staggered structure together with a fractional step method. The use
of a virtual staggered structure arrangement, which stores all the
variables at the same physical location and employs only one set of
nodes using virtual interpolation points, reduces the geometrical
complexity. The virtual staggered structure consists of the virtual
interpolation points and the virtual local stencil. Also,
computational enhancement of the virtual interpolation point method
is considerable since the present method directly discretizes the
strong forms of the incompressible Navier-Stokes equations without
numerical integration. It makes a key difference from others. In the
virtual interpolation point method, the choice of an accurate
interpolation scheme satisfying the spatial approximation in the
complex domain is important because there is the virtual staggered
structure for computation of the velocities and pressure since there
is no explicit staggered structure for stability. In our proposed
method, the high order derivative approximations for constructing
node-wise difference equations are easily obtained. Several
different flow problems (decaying vortices, lid-driven cavity,
triangular cavity, flow over a circular cylinder and a bumpy
cylinder) are simulated using the virtual interpolation point method
and the results agree very well with previous numerical and
experimental results. They verify the accuracy of the present
method.
\end{abstract}


\end{frontmatter}


\section{Introduction}
\label{} The ability to handle complex geometries has been one of
the main issue in computational fluid schemes because most
engineering problems have complex geometries. So far, two kinds of
grid arrangements to simulating complex flow have been known:
staggered grids and non-staggered grids. For the staggered grids,
vector components and scalar variables are stored at different
locations, while for the non-staggered grids, vector variables and
scalar variables are stored at the same locations, being half a
control-volume width apart in each coordinate.
Staggered grid methods are popular because of their ability to
prevent checkerboard pressure in the flow solution as discussed
in~\cite{Patankar}. The main disadvantages of such an arrangement
are the geometrical complexity due to the boundary conditions, and
the difficulty of implementation to non-orthogonal curvilinear
grids~\cite{Peric}. In the non-staggered grid methods, the main
disadvantages are the primitive variables and mass conservation in
order to solve the pressure field, either interpolation cell-face
velocities or interpolate the pressure gradients in a special way,
usually with an upwind-bias to avoid the checkerboard pressure
fields in the flow solution.

Since most engineering problems posed on complex geometries with
rough boundary pervade many fields of research(see, Figure
~\ref{fig:bumpy cylinder}), it is rather difficult for the ordinary
staggered grid method or non-staggered grid method to compute the
solution with numerical integration. For the purpose of numerical
simplicity and efficiency, we introduce a virtual interpolation
point(VIP) method using a moving least-squares(MLS) approximation
without numerical integration.

Indeed, we distribute the global nodes regularly and introduce the
virtual interpolation points on a local stencil(see,
Figure~\ref{fig:1}). In the proposed method, the high order
derivative approximations for constructing node-wise difference
equations are easily obtained. Such capability allows us to get a
local stencil estimate of the flux derivatives and thus preserve all
complicated discontinuous behaviors of solutions. We emphasize that
the mesh generation is unnecessary in our scheme and that we are
convinced that our method is more effective in higher dimensions
such as three-dimensional problem with a complex geometry.

Many researchers have been studying the meshfree method~\cite{D.W.Kim, H.J.Choe2, H.J.Choe3}. The meshfree method bases on the MLS
approximation. The meshfree is attractive because it requires no
connectivity among nodes in constructing approximation. Until now,
however, they are prone to produce a false pressure
field-checkerboard pressure. The meshfree method for hyperbolic
equations has not yet been possible in literatures due to the lack
of an innate dissipation mechanism essential to suppress numerical
oscillations by convective terms in hyperbolic equations.

Recently, an upwind meshfree method using virtual local stencil
approach was presented by Park et al~\cite{S.Park, H.J.Choe, S.Park2}, who simulated the compressible flow for the high voltage gas
blast circuit breaker with the moving boundary.

\subsection{The present contribution}
The objective of the present study is to develop the VIP method that
introduce both the virtual interpolation point and the virtual local
stencil to represent properly on complex geometries. The present
method is based on the MLS approach on a virtual staggered structure
together with a fraction-step method. The virtual staggered
structure consists of the virtual interpolation points and the
virtual local stencil. It makes a key difference from others. In
this implementation, the set of nodes for computation can be
distributed arbitrarily in principle and hence the proposed method
is applied to the flow problems on complex geometry.

The virtual local stencil(as in Park et al~\cite{S.Park}) and the
virtual interpolation points are applied only on the virtual
staggered structure. A second-order accurate interpolation scheme
for evaluating the virtual interpolation point is proposed in this
study, which is numerically stable irrespective of the relative
position between the virtual local stencil and the virtual staggered
structure. It will be also shown that introduction of the virtual
interpolation point is necessary to obtain physical solutions and
enhance accuracy.

\begin{figure}[H]
\centering \mbox{
    \subfigure[$\Psi^{[0,0]}$]{
        \includegraphics[height=3.5 cm]
        {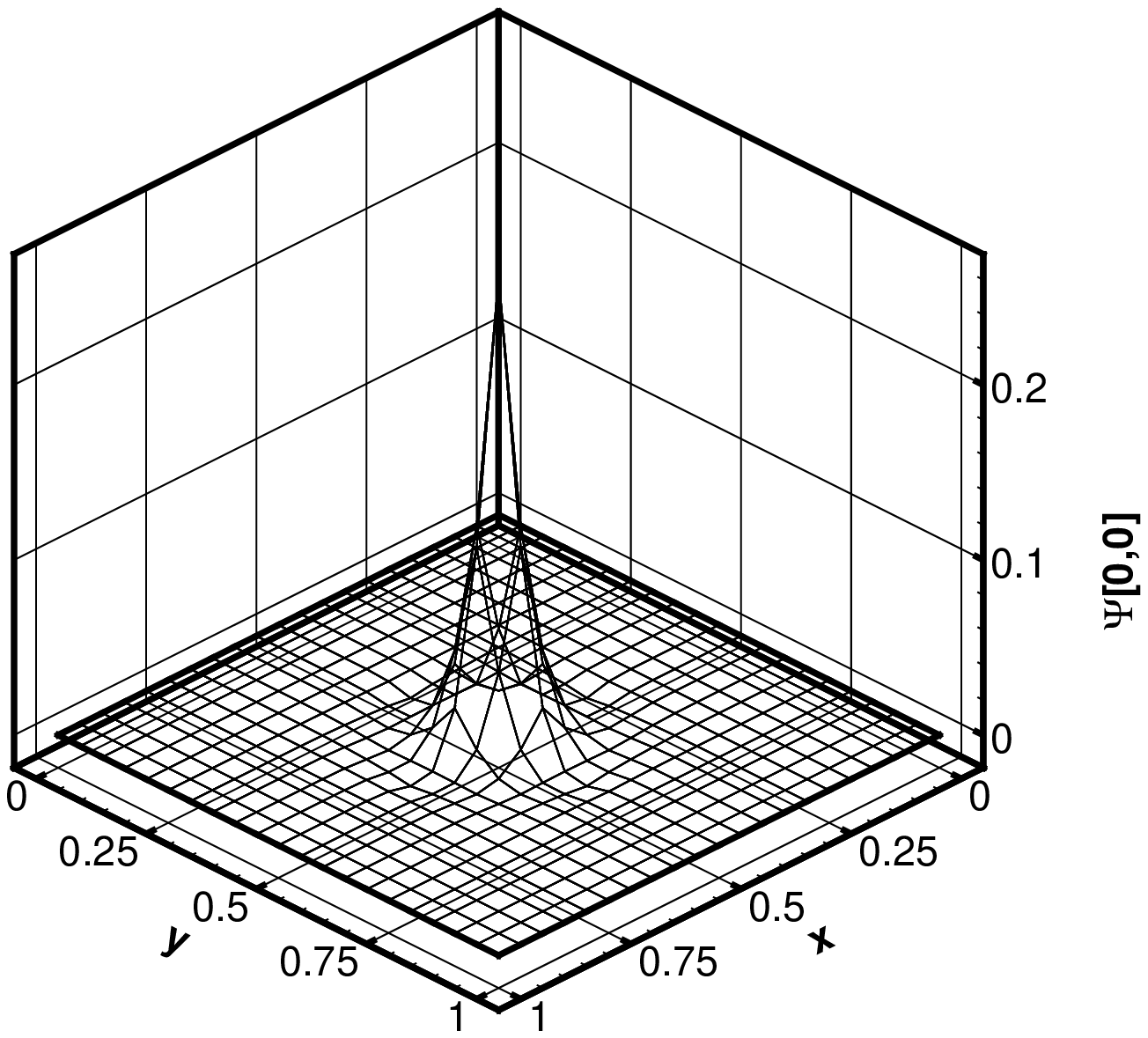}}
    \subfigure[$\Psi^{[1,0]}$]{
        \includegraphics[height=3.5 cm]
        {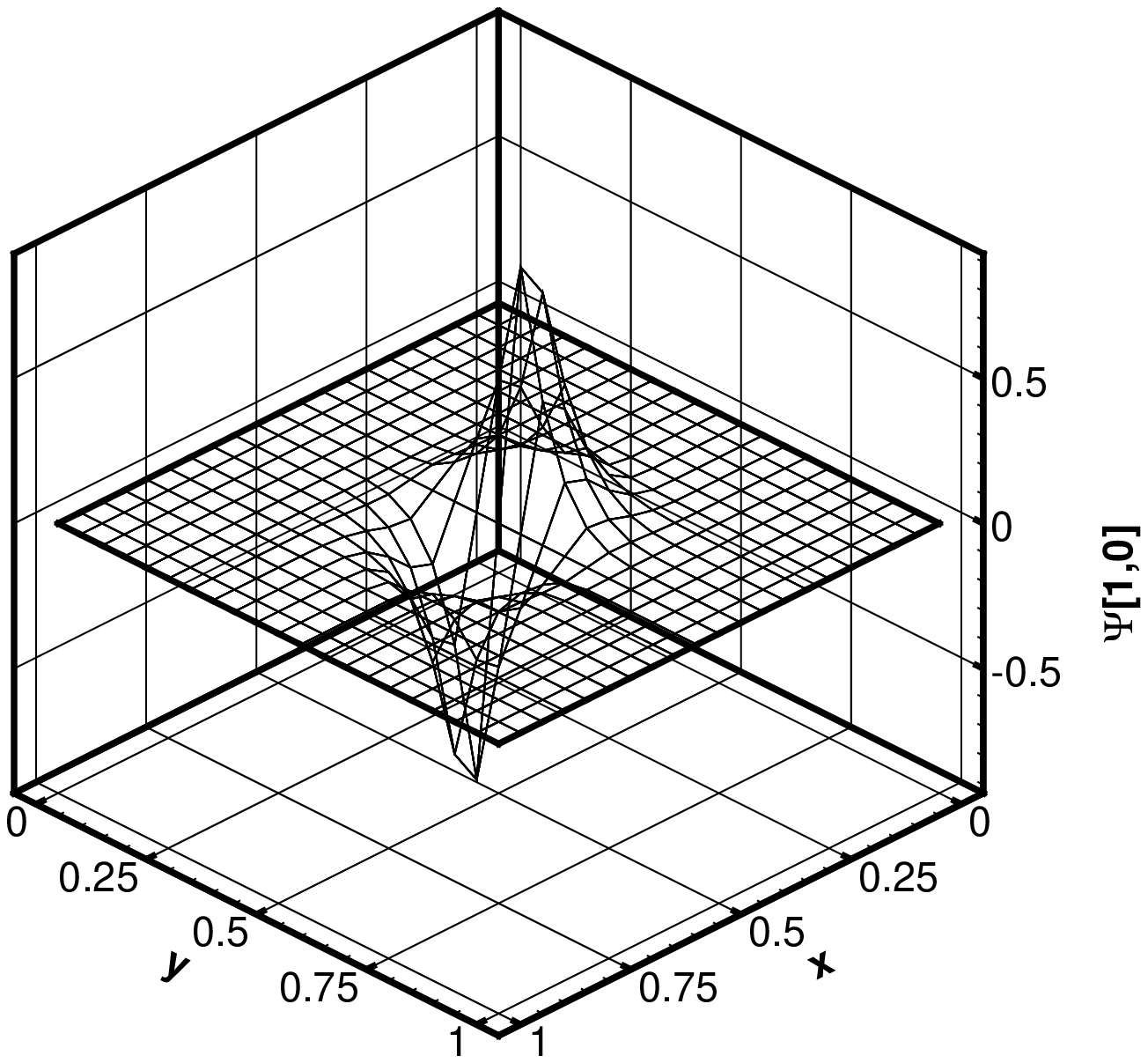}}
    \subfigure[$\Psi^{[0,1]}$]{
        \includegraphics[height=3.5 cm]
        {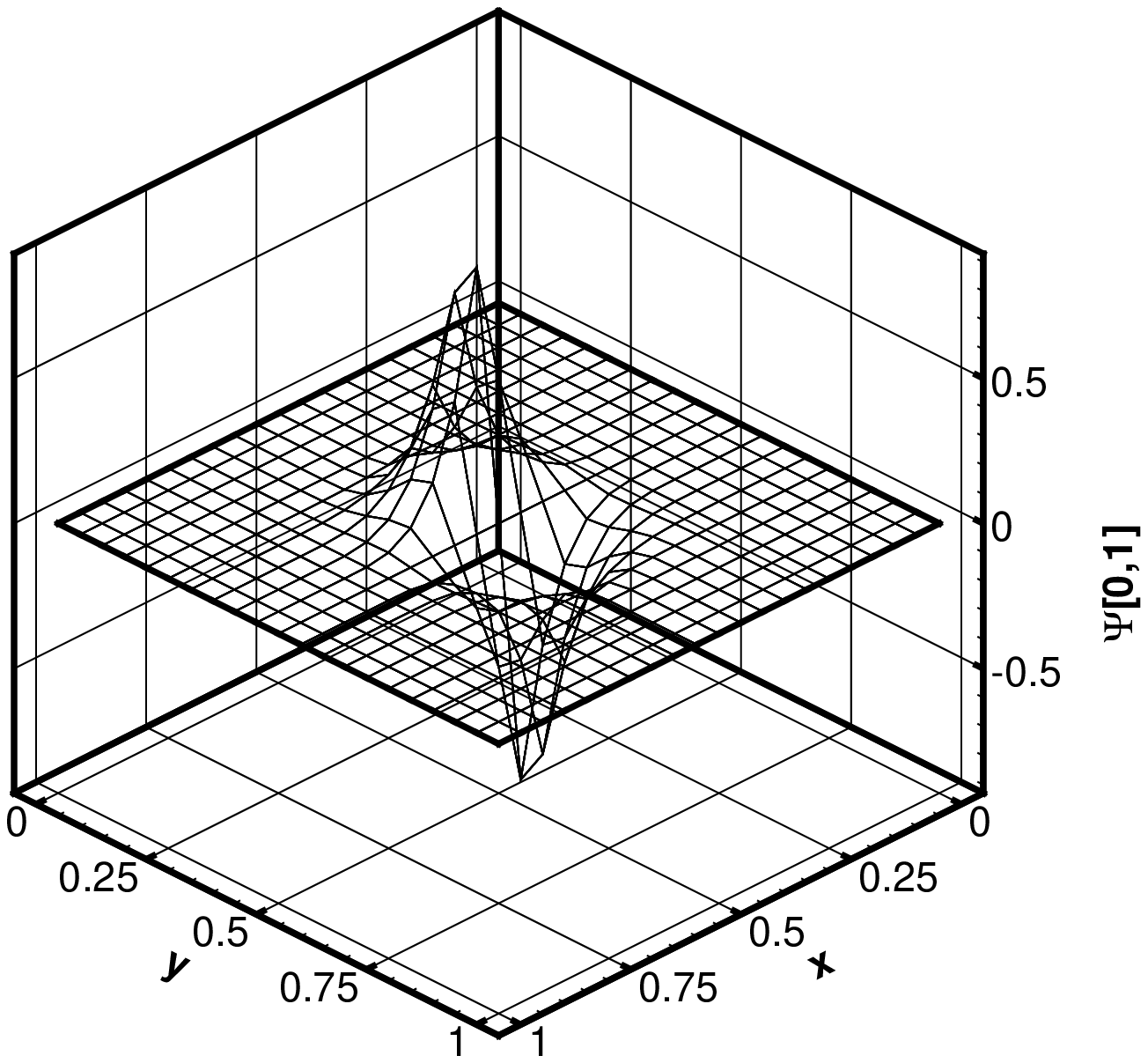}}
        }
\mbox{
    \subfigure[$\Psi^{[2,0]}$]{
        \includegraphics[height=3.5 cm]
        {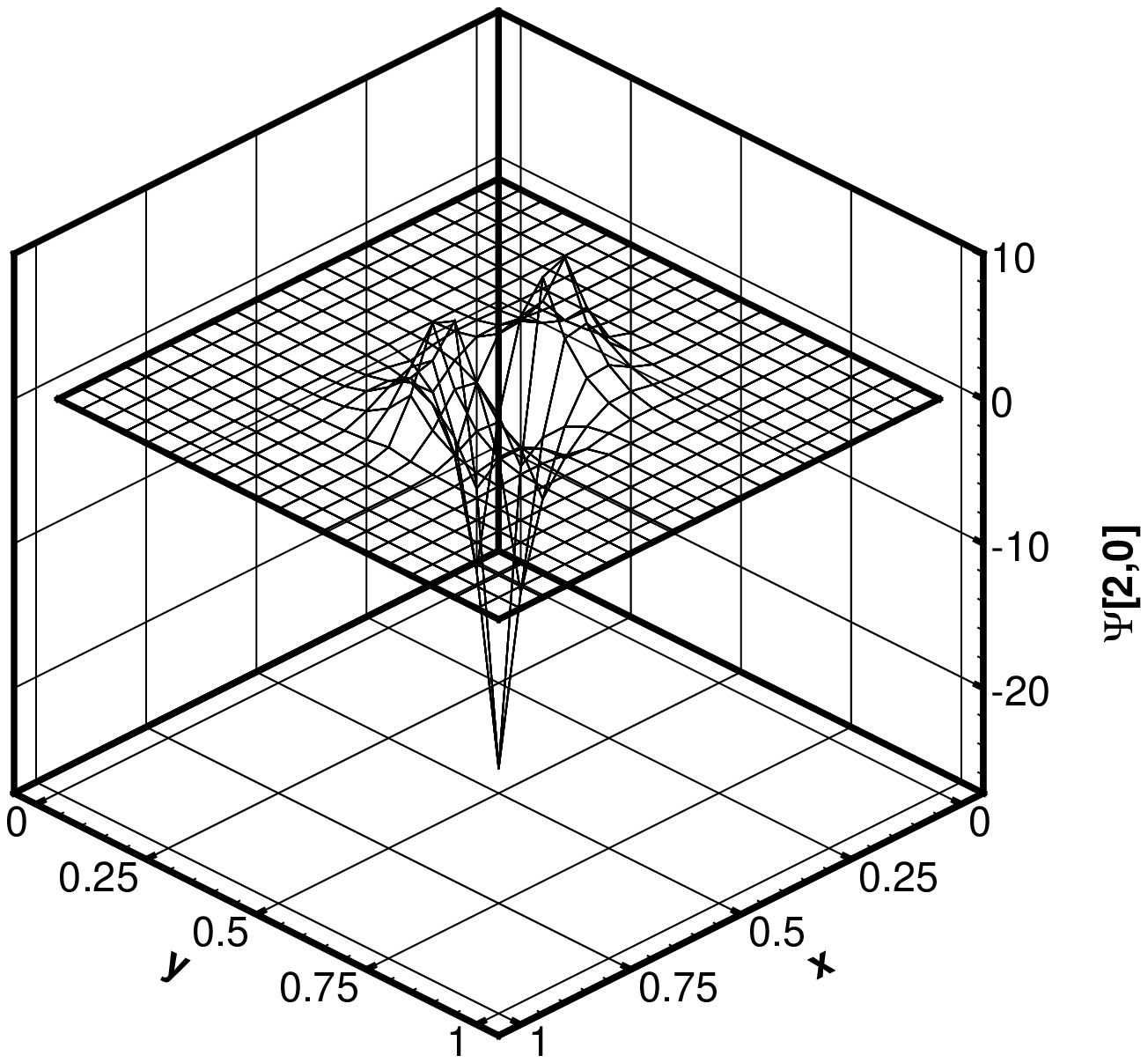}}
    \subfigure[$\Psi^{[1,1]}$]{
        \includegraphics[height=3.5 cm]
        {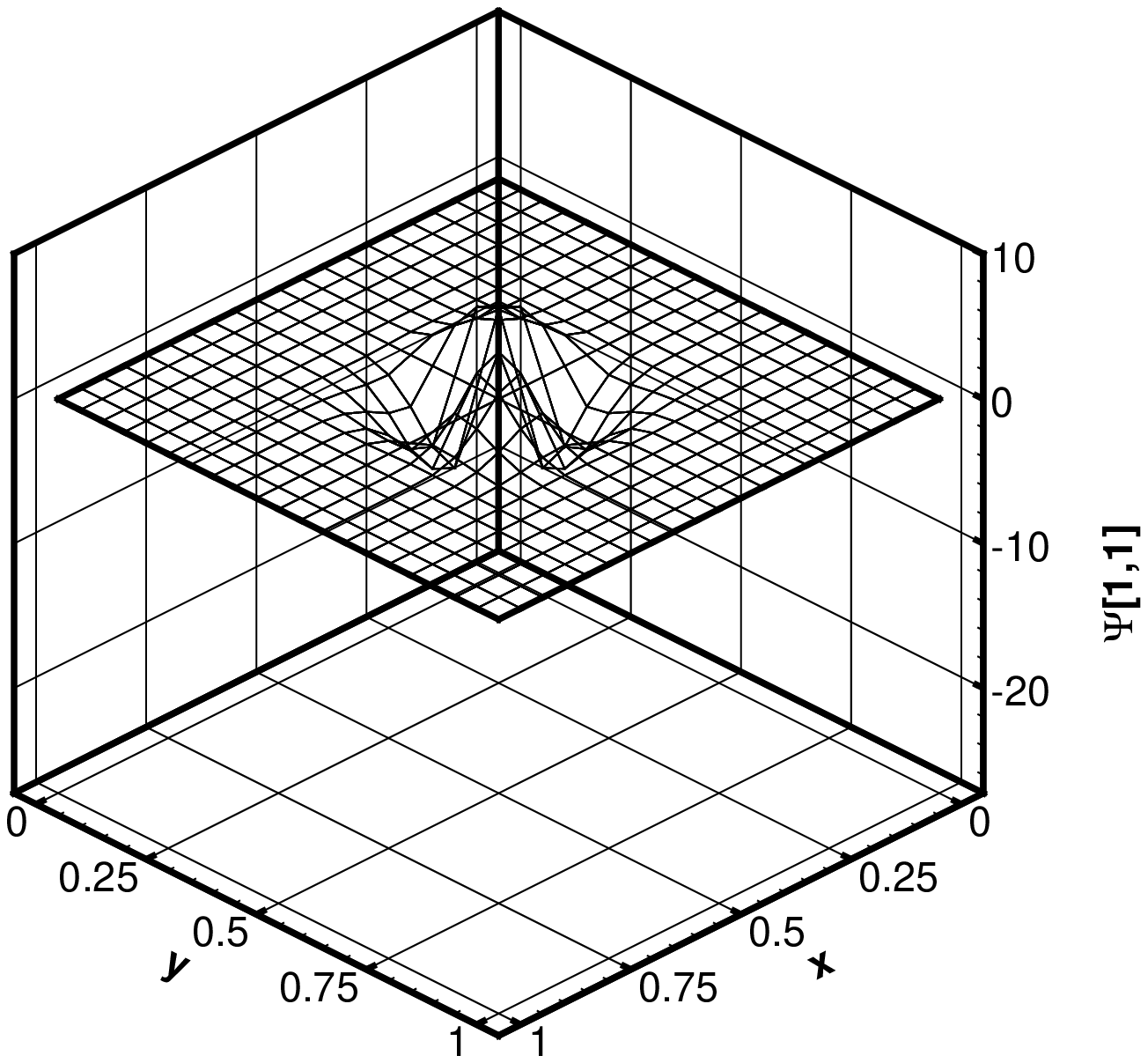}}
    \subfigure[$\Psi^{[0,2]}$]{
        \includegraphics[height=3.5 cm]
        {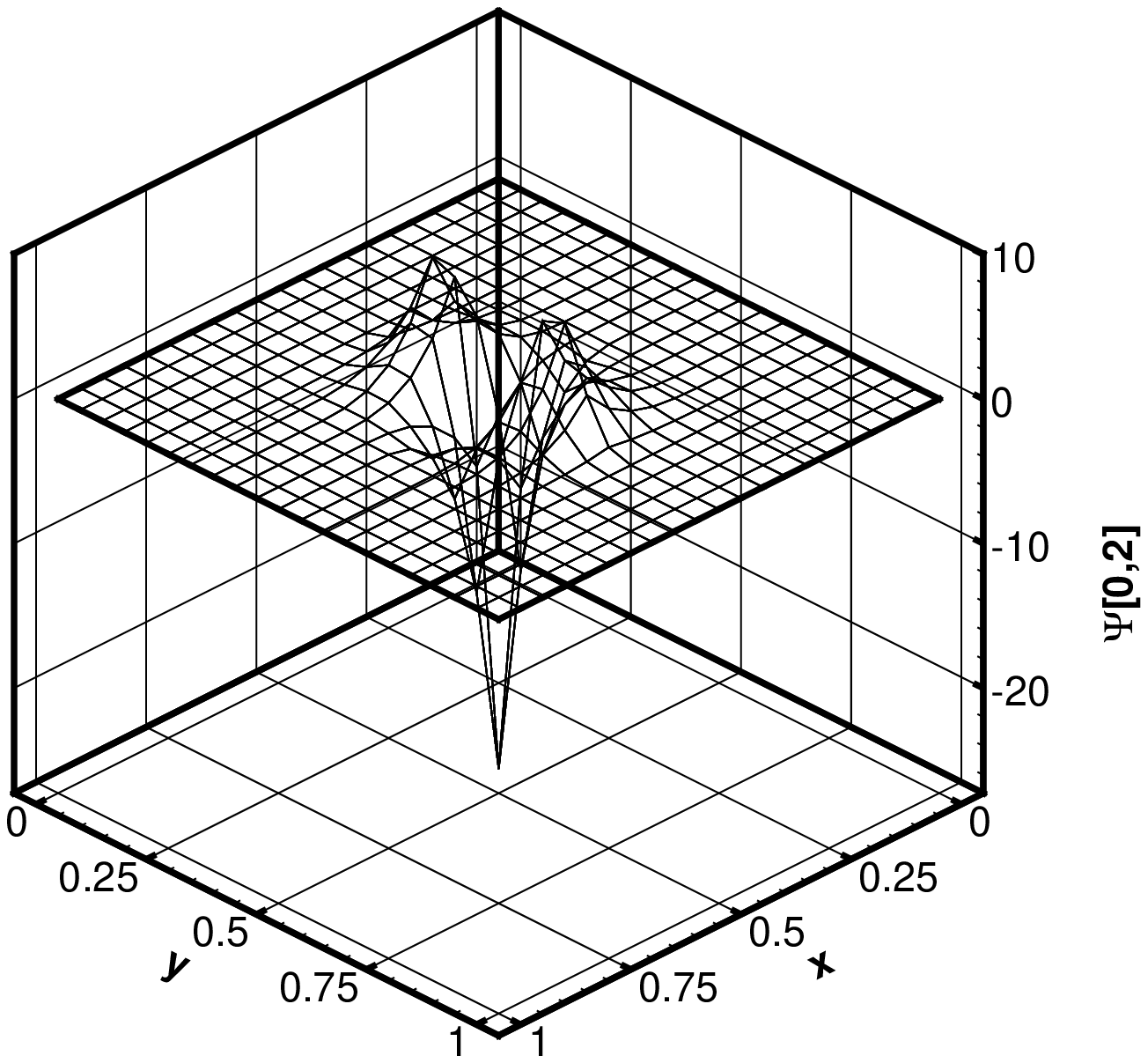}}}
\caption{ 2D $\alpha$-th shape functions with quadratic basis}
\label{fig_2D_MLS shape_quadratic}
\end{figure}

In the VIP method, the high order derivative approximations for
constructing node-wise difference equations are easily obtained(see,
Figure~\ref{fig_2D_MLS shape_quadratic}). Also, computational
enhancement is considerable since the present method directly
discretizes the strong forms of incompressible Navier-Stokes
equations without numerical integration.

The focus of this paper is laid on the contribution to a stable flow
computation without explicit structure of staggered grid. In our
method, we don't have to explicitly construct the staggered grid at
all. Instead, there exists only virtual interpolation point at each
computational node, which plays a key role in discretizing the
conservative quantities of the incompressible flow. In fact, it can
be regarded as an imaginary staggered structure, accordingly.
Particularly in our method, one set of nodes distributed on the flow
domain is needed due to the virtual interpolation point.

We emphasize that the mesh and grid generation are unnecessary in
our scheme and that we are convinced that our method is more
effective in higher dimensions such three-dimensional problem with a
complex geometry.

The reminder of the paper is organized as follows: Sections 2 and 3
present the time integration and the spatial approximation. In
section 4, we propose the stable second-order VIP method for solving
the incompressible Navier-Stokes equation. Various numerical results
are presented to show the accuracy, efficiency, stability, and
robustness and superiority of proposed scheme in Section 5. In
Section 6, conclusions are drawn.

\section{ Governing Equations and Time Integration}
\label{fractional-step} The use of a virtual staggered structure
arrangement, which stores all the variables at the same physical
location and employs only one set of nodes using virtual
interpolation points(see, Figure~\ref{fig:1}), reduces the
geometrical complexity. In the present study, the VIP scheme, as in
Park et al.~\cite{S.Park, H.J.Choe, S.Park2} is applied to
satisfy the continuity for the local stencil in complex domains.

The incompressible Navier-Stokes flows are represented with the
following governing equations,

\begin{align}
\mathbf{u}_t +  \nabla \cdot \left( \mathbf{u} \mathbf{u}^{T}
\right) + \nabla p -\frac{1}{Re}\nabla^2 \mathbf{u} &= 0\;\;\;\;
\textmd{on}
\; \Omega,\label{eq:1} \\
\nabla \cdot \mathbf{u} &= 0 \;\;\;\; \textmd{on} \; \Omega,
\label{eq:2}
\end{align}

where $\mathbf{u}$ and $p$ are the velocity components and pressure
of the flow. All the variables are nondimensionalized by the
characteristic velocity and length scales, and $Re$ is the Reynolds
number.

The time integration method used to solve Eqs.(\ref{eq:1}) and
(\ref{eq:2}) is based on a fractional step method where a
pseudo-pressure is used to correct the velocity field so that the
continuity equation is satisfied at each computational time step. In
this study, we use a second-order semi-implicit time advancement
scheme (a second-order Adams-Bashforth for the convection terms and
a second-order Crank-Nicolson method for the diffusion terms)

\begin{align}\label{eq:12}
&\mathbf{Au}^{\star} = \mathbf{r}^{n}, \\ \label{eq:13} &\Delta
\phi^{n+1} = \frac{1}{\triangle t} \nabla \cdot \mathbf{u}^{\star},
\\ \label{eq:14} &\mathbf{u}^{n+1} =
\mathbf{u}^{\star} - (\triangle t ) \nabla \phi^{n+1},  \\
\nonumber
\end{align}
\begin{align}\label{eq:15}
p^{n+1} &= p^{n} + \phi^{n+1} - \frac{\triangle t}{2Re} \Delta
\phi^{n+1},
\end{align}

where

\begin{align*}
\mathbf{A} &= \frac{1}{\triangle t} \left( \mathbf{I} -
 \frac{ \triangle t}{2 Re}
\Delta \right), \\
\mathbf{r}^n &= \frac{1}{\triangle t} \left( \mathbf{I} + \frac{
\triangle t}{2 Re} \Delta \right) \mathbf{u}^{n} - \left[ \nabla
\cdot
\left( \mathbf{uu}^{T} \right) \right]^{n+\frac{1}{2}}, \\
\end{align*}

$\mathbf{u}^{\star}$ is the intermediate velocity, and $\phi$ is the
pseudo-pressure. Also, $\triangle t$ and $\mathbf{I}$ are the
computational time step and the identity operator.

In the present study, the VIP method is applied to Cartesian
coordinate. The time-integration method is based on the method of
Kim and Moin~\cite{J.Kim} to enhance computational efficiency.

\section{Moving Least-Squares Approximation}

For spatial approximation of solutions, the MLS approximation in the
literature~\cite{D.W.Kim} is employed. We briefly explain the MLS
approximation.

For simplicity, we just consider 2-dimensional space and take $m=2$
but it can be extended to $n$-dimension. Multi-index notations are
adapted throughout the paper
\begin{align*}
\mathbf{x} = (x, y) \textmd{ and }
\partial^{\mathbf{(\alpha}, \beta)}_{\mathbf{x}} =
\partial_{x}^{\alpha} \partial_{y}^{\beta},
\end{align*}
where $\alpha$ and $\beta$ are non-negative integers. For a
continuous function $u(\mathbf{x})$ we can approximate this function
at a point $\bar{\mathbf{x}} = (\bar{x},\, \bar{y})$ in terms of
polynomials up to some order dependently on a neighborhood of
$\bar{\mathbf{x}}$, which is found by Weierstass.

Let $u_m(\mathbf{x}, {\bar{\mathbf{x}}})$ be a polynomial up to
degree $m$ which depends on the point $\bar{\mathbf{x}}$. Then for
some coefficient vectors $\mathbf{c}$, it can be assumed that
\begin{align*}
u_m(\mathbf{x}, {\bar{\mathbf{x}}}) &=
\mathbf{c}_{1}(\bar{\mathbf{x}})
+\left(\frac{x-\bar{x}}{\rho(\bar{\mathbf{x}})}\right)\mathbf{c}_{2}(\bar{\mathbf{x}})
+\left(\frac{y-\bar{y}}{\rho(\bar{\mathbf{x}})}\right)\mathbf{c}_{3}(\bar{\mathbf{x}})
\\
&+\left(\frac{x-\bar{x}}{\rho(\bar{\mathbf{x}})}\right)^{2}\mathbf{c}_{4}(\bar{\mathbf{x}})
+\left(\frac{x-\bar{x}}{\rho(\bar{\mathbf{x}})}\right)\left(\frac{y-\bar{y}}{\rho(\bar{\mathbf{x}})}\right)\mathbf{c}_{5}(\bar{\mathbf{x}})
+\left(\frac{y-\bar{y}}{\rho(\bar{\mathbf{x}})}\right)^{2}\mathbf{c}_{6}(\bar{\mathbf{x}}).
\end{align*}
The dilation function $\rho(\bar{\mathbf{x}})$ can be regarded as
 the size of a neighborhood at $\bar{\mathbf{x}}$ for approximation.

In order to find the best approximation $u_m(\mathbf{x},
{\bar{\mathbf{x}}})$ with the coefficient $
\mathbf{c}(\bar{\mathbf{x}})$, we define the locally weighted square
functional of the form on a given set of nodes, $\Lambda = \{
\mathbf{x}_i \in \bar{\Omega} | i = 1, 2, \cdots, N \}$.

\begin{align*}
J(\mathbf{c}(\bar{\mathbf{x}})) = \sum_{i=1}^{N}
\Phi\left(\frac{\mathbf{x}-\bar{\mathbf{x}}}{\rho(\bar{\mathbf{x}})}\right)|{
u_m(\mathbf{x}_i, {\bar{\mathbf{x}}}) - u(\mathbf{x}_i)}|^2
\end{align*}
where the weight function is taken as the following form,

\begin{equation*}
\Phi\left(\mathbf{y}\right)=
\left\{%
\begin{array}{ll}
    \left(1-\|\mathbf{y}\|^{\left(1-\|\mathbf{y}\|\right)^{2}}\,\right)^4, & \quad if \quad \| \mathbf{y} \|< 1  \\
    0\quad\quad\quad\quad\quad\quad\quad\quad, &  \quad otherwise. \\
\end{array}%
\right.
\end{equation*}

Minimizing the functional $J$, the local approximation
$u_m(\mathbf{x}, {\bar{\mathbf{x}}})$ is determined with coefficient
$\mathbf{c}(\bar{\mathbf{x}})$. In fact, it is the best
approximation, partially near $\bar{\mathbf{x}}$. Moreover, $u_m$ is
a polynomial is $\mathbf{x}$, so that we can differentiate it as
many times as we want. It is also natural that the derivative of
$u_m(\mathbf{x}, {\bar{\mathbf{x}}})$ w.r.t. $\mathbf{x}$ are good
approximate of derivatives. Therefore, we pay our attention to the
values of $u_m(\mathbf{x}, {\bar{\mathbf{x}}})$ and its derivatives
at $\bar{\mathbf{x}}$. This observation produces the following
approximates,

\begin{align*}
\left( \mathcal{D}^{[\alpha, \beta]} u\right)(\bar{\mathbf{x}})
\equiv \lim_{\mathbf{x}\rightarrow\bar{\mathbf{x}}}
\partial^{(\alpha, \beta)}_{\mathbf{x}}u_m(\mathbf{x},
{\bar{\mathbf{x}}}).
\end{align*}

Following the above procedure, we finally have the representation
formula for approximated derivatives,
\begin{align}\label{fmlsrk:2}
\left( \mathcal{D}^{[\alpha, \beta]} u\right)(\bar{\mathbf{x}})
=\sum_{i=1}^{N} u(\mathbf{x}_i) \Psi^{[\alpha,
\beta]}_{i}(\bar{\mathbf{x}}),
\end{align}
in which we call $\Psi^{[\alpha, \beta]}_{i}(\bar{\mathbf{x}})$ the
$[\alpha, \beta]$-th approximative of a shape function at
$\bar{\mathbf{x}}_i$. For detailed description, see the
reference~\cite{D.W.Kim}.

\section{Implementation using Virtual Interpolation Point}

Using the conventional MLS approximations only, we have empirically
experienced that the fractional step method becomes unstable. This
is why we elaborate the VIP scheme on a virtual local stencil for
stability.

The VIP method is developed for the solution of computational fluid
dynamics problems that does not require the use of staggered grid
systems. This implementation of this scheme is performed on only one
set of nodes for both velocities and pressure. It makes a key
difference from others. In this implementation, the set of nodes for
computation can be distributed arbitrarily in principle and hence
the proposed method can be applied to the flow problems on
complicated geometry.

\subsection{Numerical Flux using the VIPs on a Local Stencil}
The key idea of VIP scheme is that conservative variables are
obtained by the conventional MLS approximation at the auxiliary
virtual interpolation points, which are not necessary nodes. In the
present method, the choice of an accurate interpolation scheme
satisfying the spatial approximation in the complex domain is
important because there is the virtual staggered structure for
computation of the velocities and pressure but there is no explicit
staggered structure for stability. In the proposed method, the high
order derivative approximations for constructing node-wise
difference equations are easily obtained.

We first introduce the approximations of the identity and Laplacian
operators,
\begin{align}\label{vip:1}
\mathbf{I} = \mathcal{D}^{[0,0]} \textrm{ and } \Delta =
\mathcal{D}^{[2,0]} + \mathcal{D}^{[0,2]},
\end{align}
using the MLS approximations in~(\ref{fmlsrk:2}). The Laplacians
in~(\ref{eq:12}),~(\ref{eq:13}),and~(\ref{eq:15}) are all replaced
with the operator in~(\ref{vip:1}). In addition, every variable
without differential operators is approximated through the identity
operator.

$\left[ \nabla \cdot \left( \mathbf{u} \mathbf{u}^T \right)
\right]^{n+1/2}$ is the non-linear convection term at node
$\mathbf{x}$. Implicitly handling viscous term eliminates the
numerical instability due to the CFL restriction. The term $\left[
\nabla \cdot \left( \mathbf{u} \mathbf{u}^T \right) \right]^{n+1/2}$
is approximated in a second-order temporal approximation for the
convective derivative term at time level $t_{n+1/2}$ which is
usually called Adams-Bashforth,
\begin{align*}
\left[ \nabla \cdot \left( \mathbf{u} \mathbf{u}^T \right)
\right]^{n+1/2} = \frac{1}{2} \left[ 3 \nabla \cdot \left(
\mathbf{u}^n \left(\mathbf{u}^n\right)^{T} \right) - \nabla \cdot
\left( \mathbf{u}^{n-1} \left(\mathbf{u}^{n-1}\right)^{T} \right)
\right] + O(\Delta t^2).
\end{align*}

\begin{figure}[H]
    \centering
    \includegraphics[width=13cm]{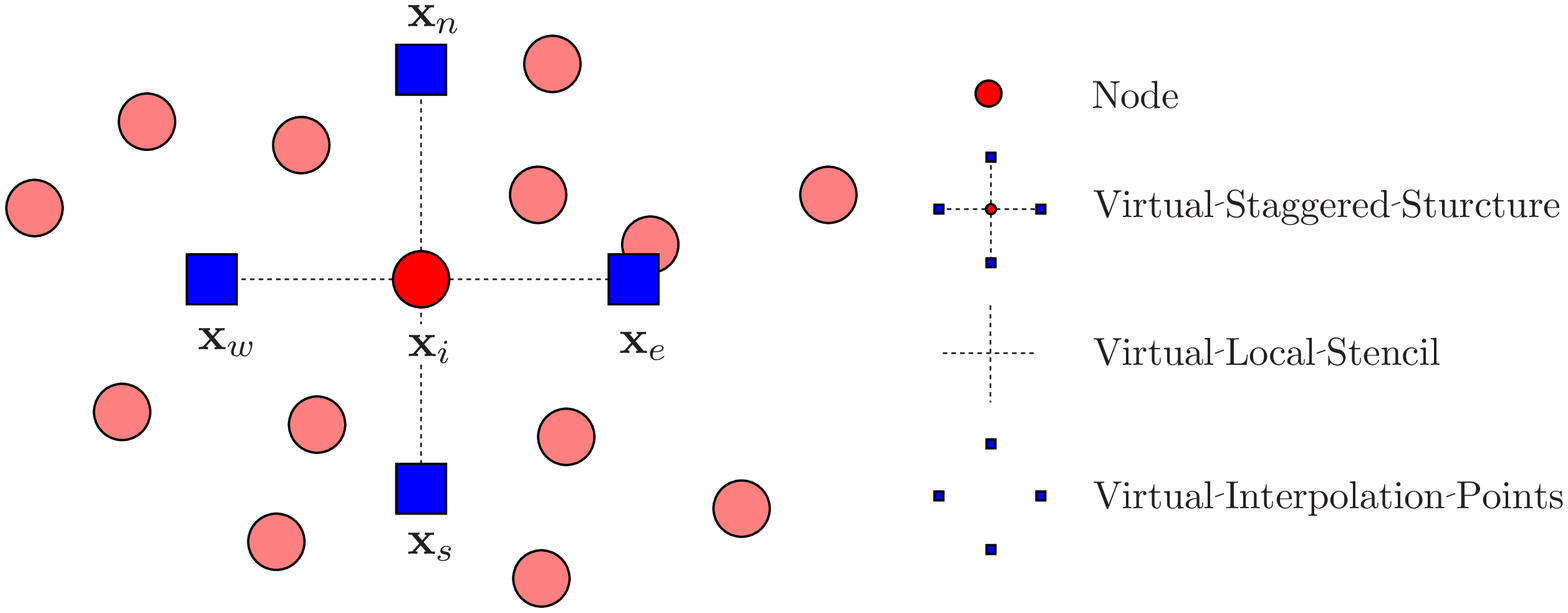}
    \caption{Schematic diagram for virtual interpolation
points on the local stencil at a node $\mathbf{x}_i \in
\Lambda$.}\label{fig:1}
\end{figure}

Second, instead of directly applying the derivative approximations
in MLS approximation to the convective terms, we simply take the
direct difference for implementing the divergence operator as in the
finite difference method. Let $ {\mathbf{x}}_{e}, {\mathbf{x}}_{w},
{\mathbf{x}}_{s},$ and ${\mathbf{x}}_{n}$ denote east, west, south,
and north points from the node $\mathbf{x}_i$ on the local stencil
in Fig. \ref{fig:1}. We call these points the virtual interpolation
points of $\mathbf{x}_i$. When $\mathbf{x}_i$ is far away from the
boundary, it is not difficult to choose the virtual interpolation
points around $\mathbf{x}_i$. However, technical problem can happen
in case where the node $\mathbf{x}_i$ is close to the boundary of
the computational domain.

Three conservative terms under consideration are, $\nabla \cdot
\left( \mathbf{uu}^{T}\right)$ appearing in $\mathbf{r}^n$
in~(\ref{eq:12}), $\nabla \cdot \mathbf{u}$ in~(\ref{eq:13}), and,
$\nabla \phi$ in~(\ref{eq:14}). At each interior node $\mathbf{x}_i
\in \Omega \cap \Lambda$, the following discretizations are used for the local
conservation over a local stencil consisting to the virtual
interpolation points at node $\mathbf{x}_i \in \Omega \cap \Lambda$;

\begin{enumerate}
\item {Convection Term using the VIPs on the Virtual Local Stencils}
\begin{align*}
\nabla \cdot \left( \mathbf{u} \mathbf{u}^T
\right)\bigg|_{\textbf{x}_i}  &\approx \left(
  \begin{array}{c}
    \frac{u^2(
{\mathbf{x}}_{e}) - u^2({\mathbf{x}}_{w})}{\| {\mathbf{x}}_{e} -
{\mathbf{x}}_{w}\|}
    + \frac{u(
{\mathbf{x}}_n)\,v( {\mathbf{x}}_n) - u({\mathbf{x}}_s)\,v(
{\mathbf{x}}_s)}{\| {\mathbf{x}}_n -
{\mathbf{x}}_s\|}   \\
    \frac{u(
{\mathbf{x}}_e)\,v( {\mathbf{x}}_e) - u({\mathbf{x}}_w)\,v(
{\mathbf{x}}_w)}{\| {\mathbf{x}}_e - {\mathbf{x}}_w\|}
    + \frac{ v^2( {\mathbf{x}}_n) - v^2(
{\mathbf{x}}_s)}{\| {\mathbf{x}}_n -
{\mathbf{x}}_s\|}\\
  \end{array}
\right), \\
\end{align*}

where
\begin{align*}
  \mathbf{u}(\mathbf{x}_{_{vip}})
  =  \sum_{i=1}^{N} \Psi^{[0,0]}_{i}(\mathbf{x}_{_{vip}})\mathbf{u}_i, \quad vip = e,w,n,s,
\end{align*}

\item{Divergence Term using the VIPs on the Virtual Local Stencils}
\begin{align*}
\nabla \cdot \left( \mathbf{u}^{\star} \right)\bigg|_{\textbf{x}_i}
&\approx
    \begin{array}{c}
      \frac{u^{\star}( {\mathbf{x}}_e) -
u^{\star}({\mathbf{x}}_w)}{\| {\mathbf{x}}_e - {\mathbf{x}}_w\|}
    + \frac{v^{\star}( {\mathbf{x}}_n) -
v^{\star}( {\mathbf{x}}_s)}{\| {\mathbf{x}}_n -
 {\mathbf{x}}_s\|}, \\
    \end{array}
  \\
\end{align*}
where
\begin{align*}
  \mathbf{u}^{\star}(\mathbf{x}_{_{vip}})
  =  \sum_{i=1}^{N} \Psi^{[0,0]}_{i}(\mathbf{x}_{_{vip}})\mathbf{u}^{\star}_i, \quad vip = e,w,n,s,
\end{align*}

\item{Gradient Term using the VIPs on the Virtual Local Stencils}
\begin{align*}
\nabla \phi \bigg|_{\textbf{x}_i} &\approx\left(
    \begin{array}{c}
      \frac{\phi({\mathbf{x}}_e) -
\phi( {\mathbf{x}}_w)}{\|
{\mathbf{x}}_e -  {\mathbf{x}}_w\|} \\
     \frac{\phi({\mathbf{x}}_n) -
\phi({\mathbf{x}}_s)}{\| {\mathbf{x}}_n -
{\mathbf{x}}_s\|} \\
    \end{array}
  \right),\\
\end{align*}
where
\begin{align*}
  \phi(\mathbf{x}_{_{vip}})
  =  \sum_{i=1}^{N} \Psi^{[0,0]}_{i}(\mathbf{x}_{_{vip}})\phi_i, \quad vip = e,w,n,s.
\end{align*}

\end{enumerate}

\section{Numerical experiments}
In this section, several different flow problems (Decaying vortices,
lid-driven cavity flow, triangular cavity flow, flow over a circular
cylinder and a bumpy circular cylinder) are simulated using the VIP
method proposed in this study and the results agree very well
previous numerical and experimental results, verifying the accuracy
of the present method.

\subsection{Taylor decaying vortices}
The temporal and spatial accuracy of the VIP method is verified by
simulating the two-dimensional unsteady flows such as
\begin{align*}
    u(x, y, t) &= -\cos x \sin y e^{-2 t},    \\
    v(x, y, t) &= \sin x \cos y e^{-2 t}, \\
    p(x, y, t) &= -\frac{1}{4}\left[\cos 2x + \cos y\right] e^{-4  t}.
\end{align*}
We consider the Taylor decaying vortices on the domain $\Omega = (0,
\pi) \times (0, \pi)$, which is discretized with regular nodes.

\begin{figure}[H]
\mbox{
        \subfigure[{$\Delta t = 10^{-4}.$}]
        {
            \includegraphics[width=0.46\textwidth]{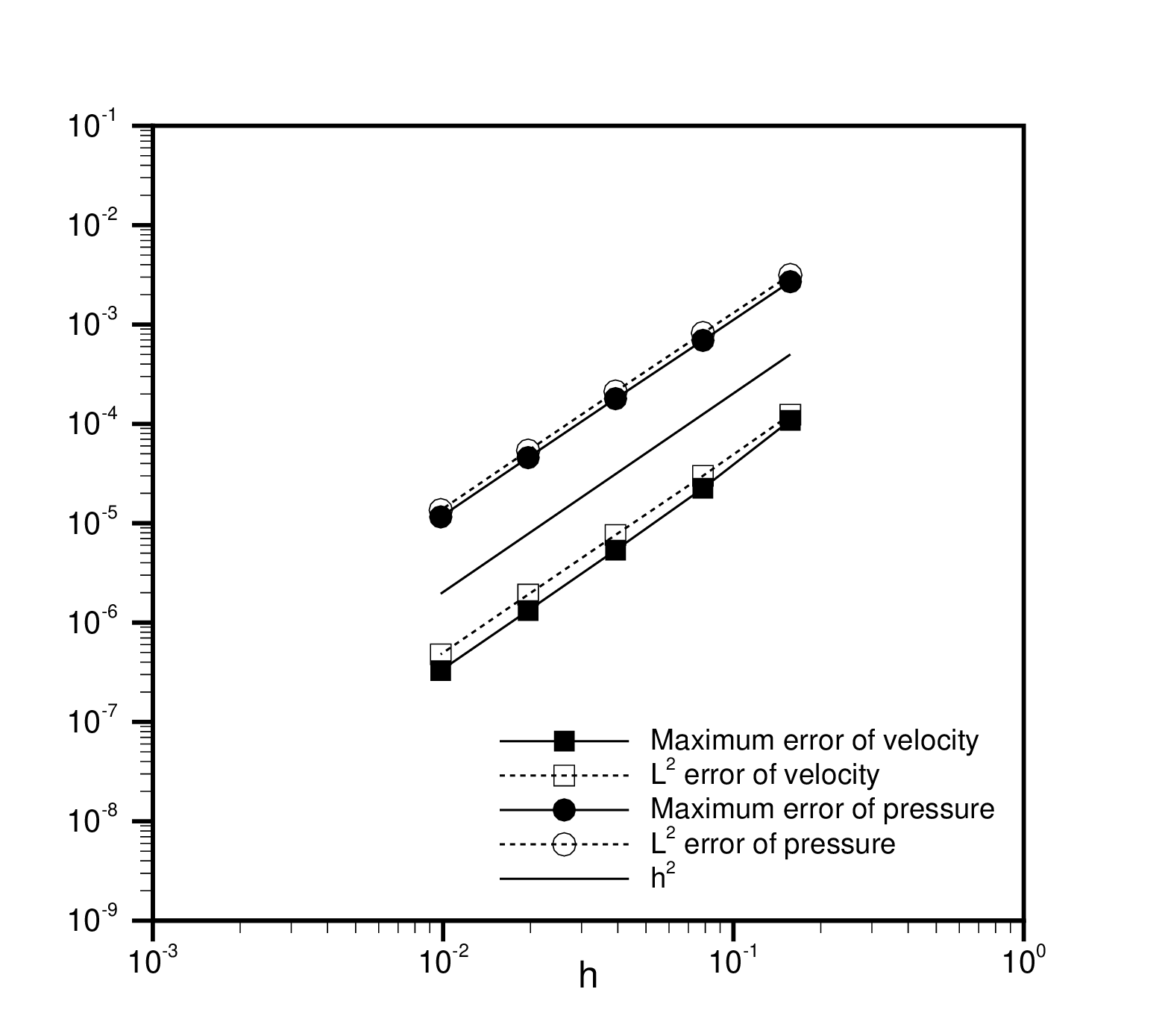}
        }
}
\mbox{
        \subfigure[{$h = \pi / 320$.}]
        {
                \includegraphics[width=0.46\textwidth]{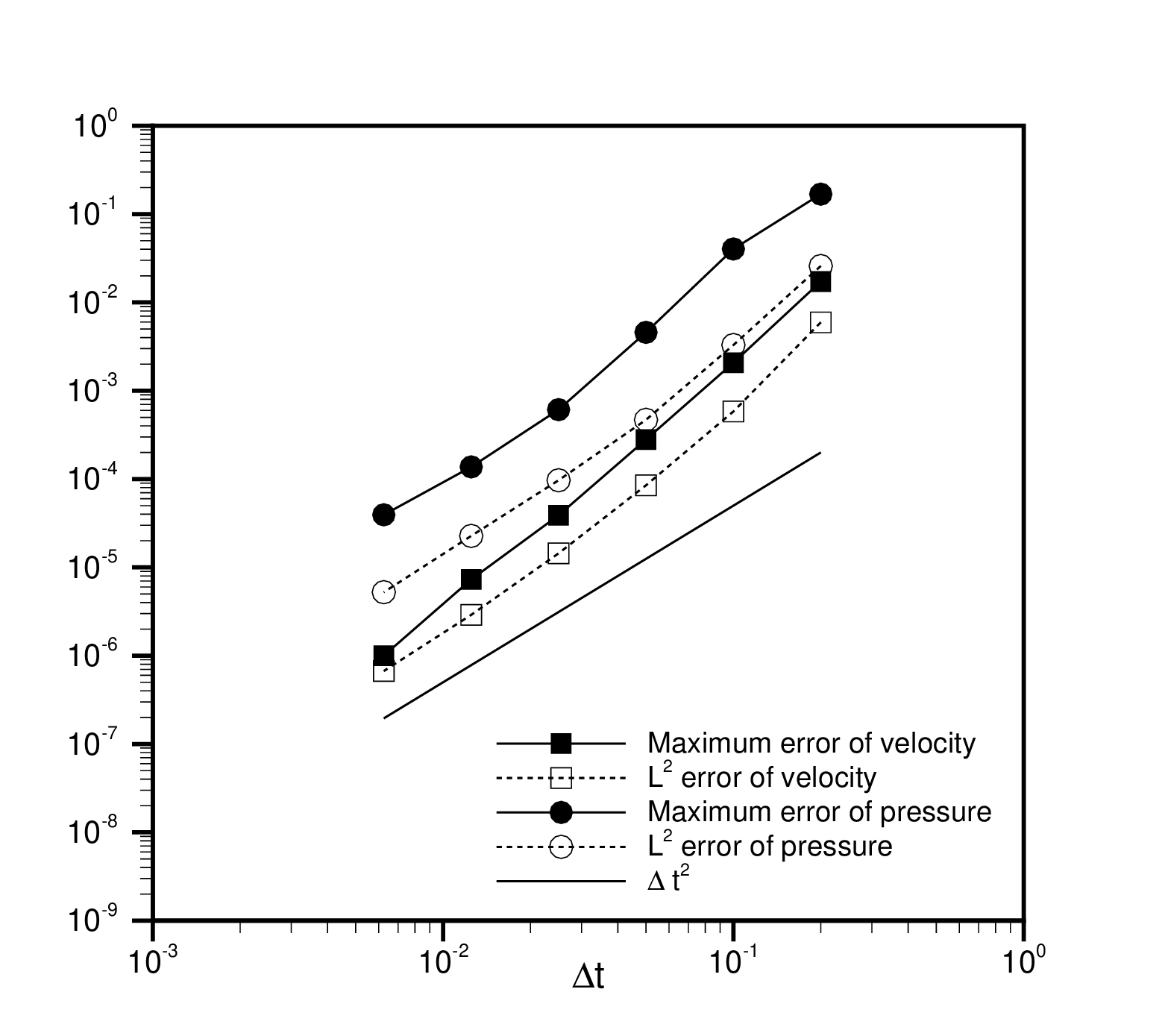}
        }
}
 \caption{Error convergence plots of Taylor decaying vortex: (a) Spatial of the VIP method evaluated using the unsteady Navier-Stokes Taylor vortex analytical
    solution at t = 0.1; (b) Temporal accuracy of the VIP method evaluated using the unsteady Navier-Stokes Taylor vortex analytical solution at t = 1.0.}\label{fig:Decaying Vortices}
    \end{figure}

As shown in the Fig.~\ref{fig:Decaying Vortices}, we have obtained
the convergence results of $O(h^2)$ and $O(\Delta t^2)$ for the
temporal and spatial on uniform nodes, respectively.

\subsection{Lid-driven cavity flow}

\begin{figure}[H]
        \centering
        \subfigure[22,925 nodes]
        {
            \includegraphics[height=4.7cm]{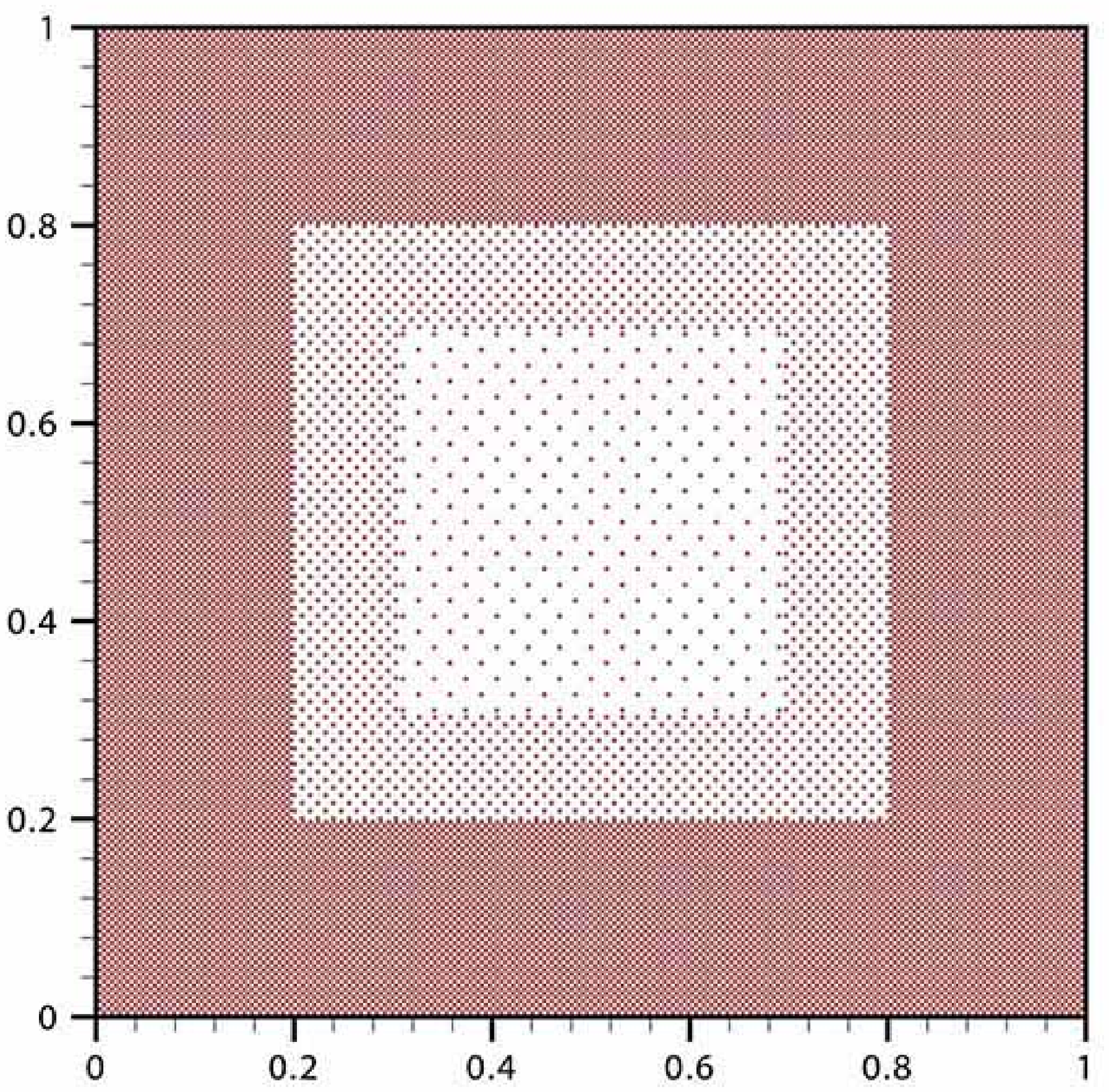}
        }
        \subfigure[12,048 nodes]
        {
            \includegraphics[height=4.7cm]{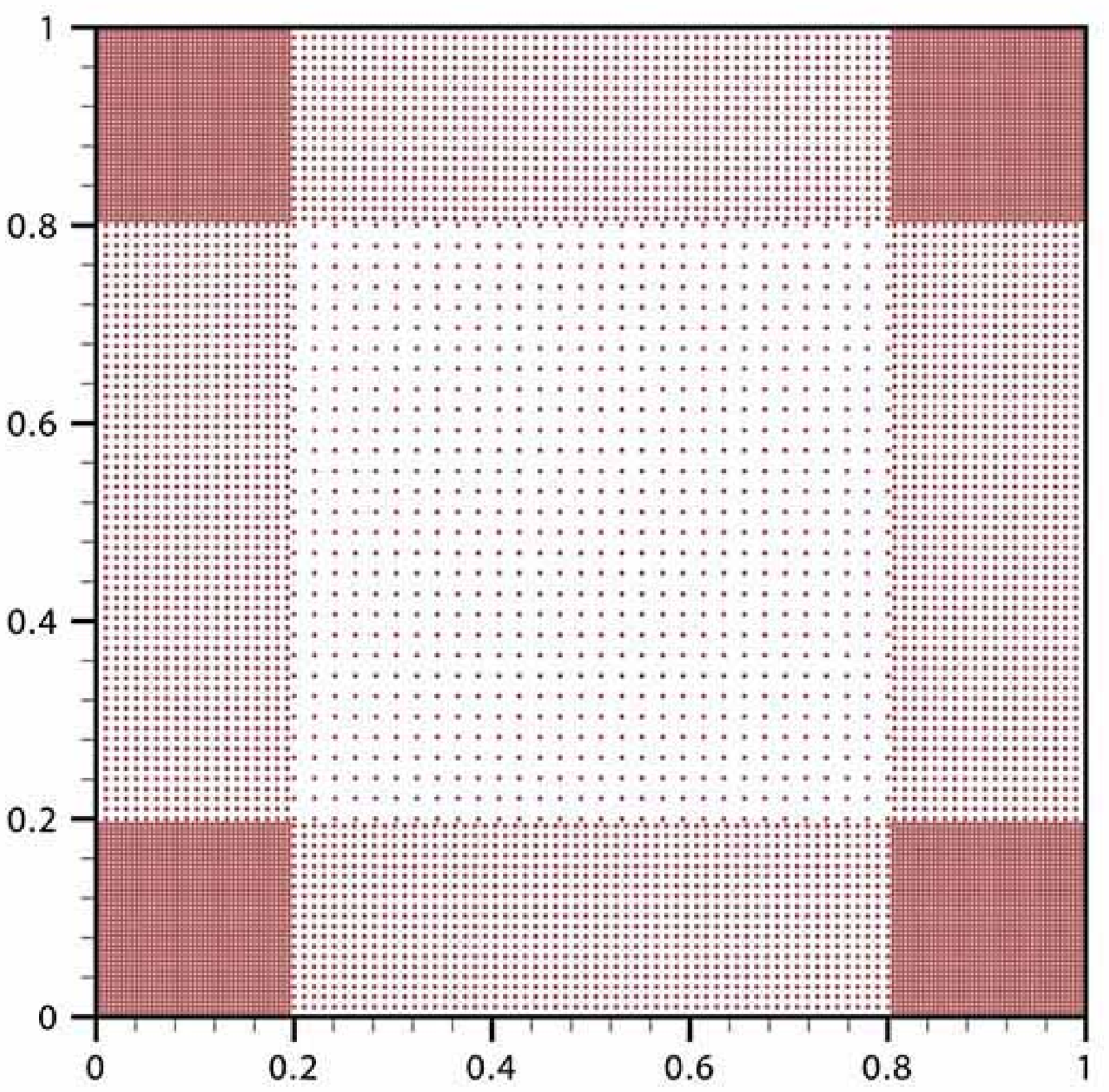}
        }\\
        \subfigure[7,897 nodes]
        {
            \includegraphics[height=4.7cm]{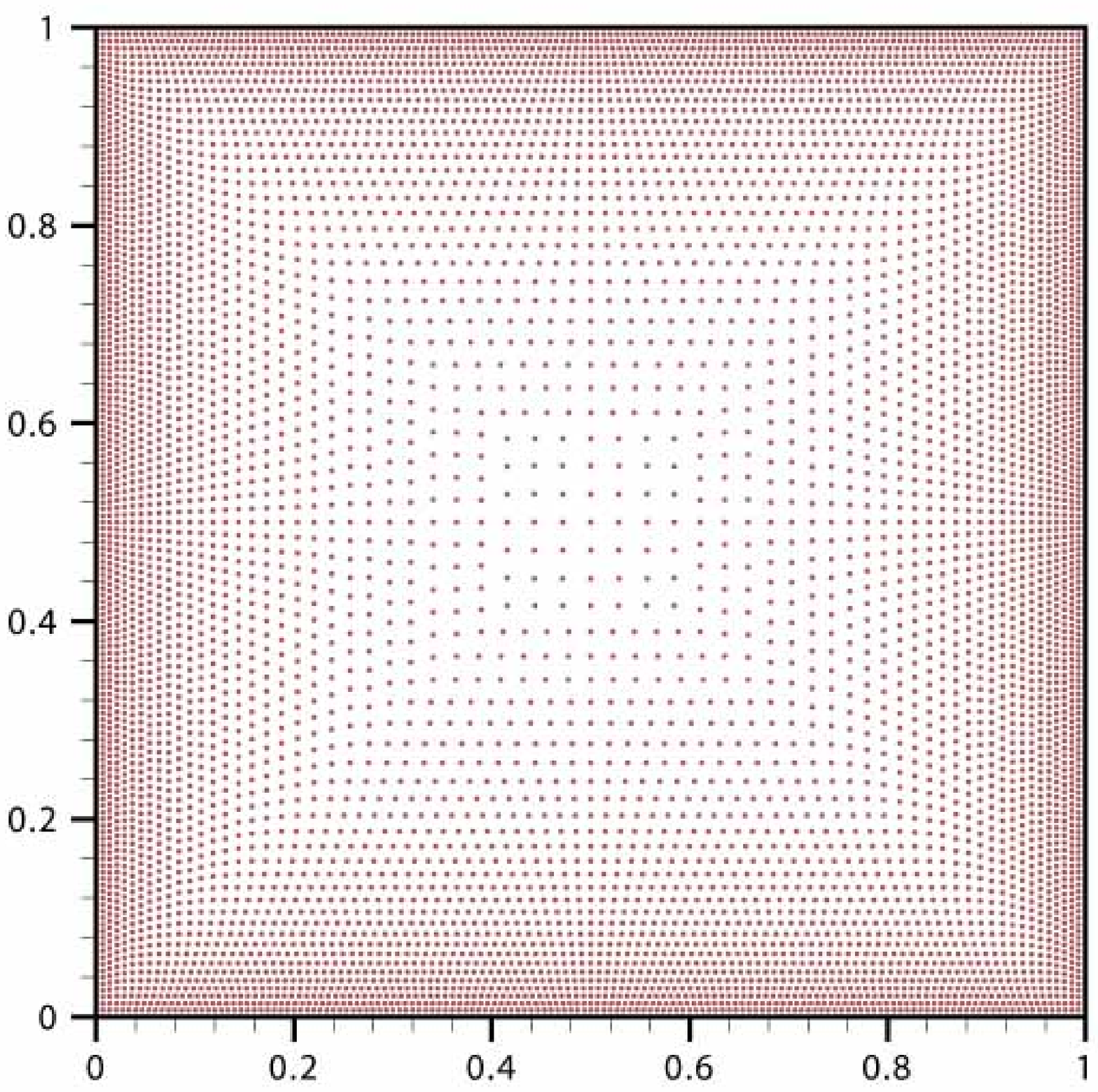}
        }
        \subfigure[]
        {
            \includegraphics[height=4.7cm]{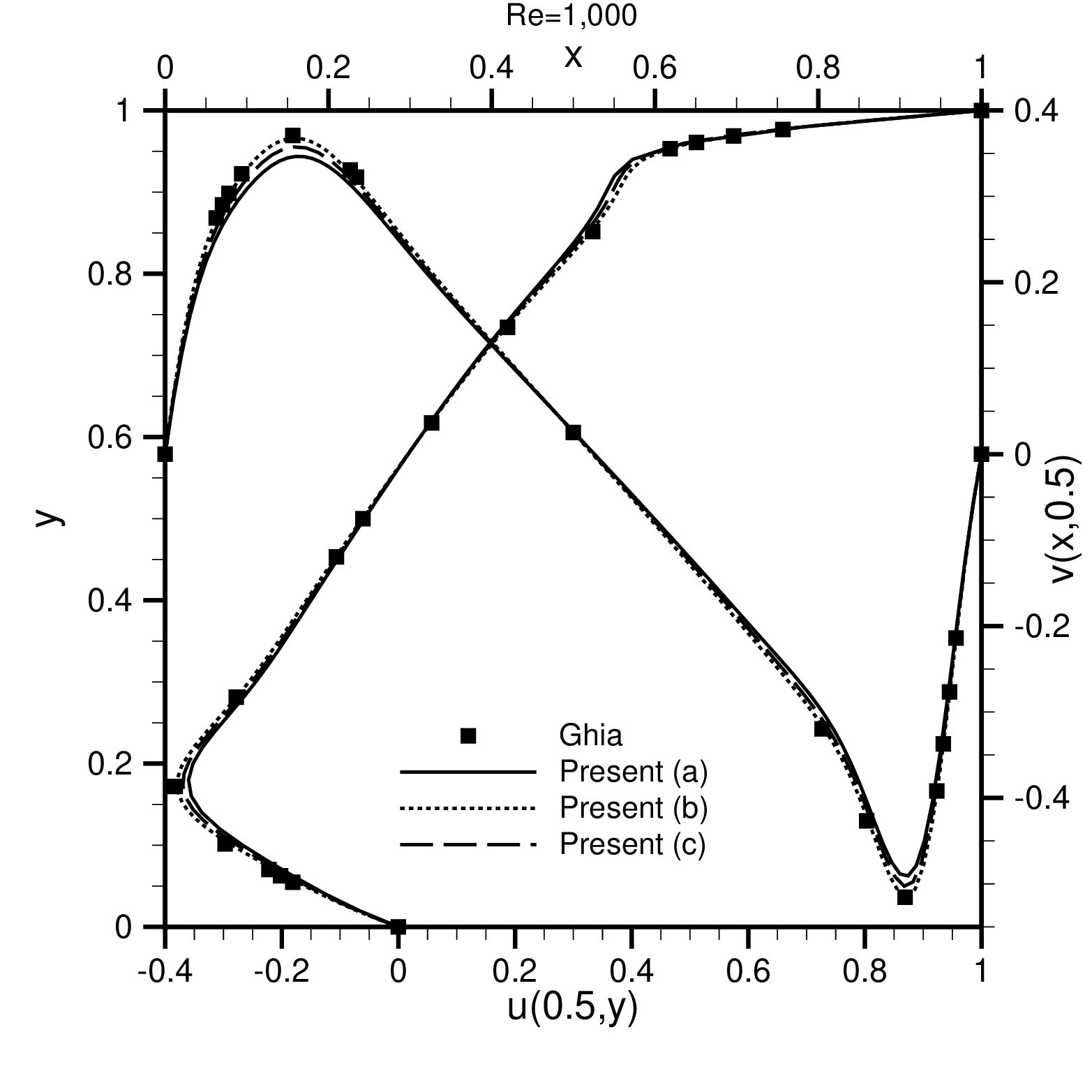}
        }
        \caption{Cavity flow: (a), (b), and (c) are various regular distributed
        nodes(a-c); (d) velocity profiles along middle sections for
        the various regular nodes of the square cavity problem
        $Re$=1000.}
        \label{fig:cavity}
\end{figure}

This classical problem has become a standard benchmark for assessing
the performance of algorithms for the incompressible Navier-Stokes
equations. For a typical example of the interior flow with corner
singularity, many researchers have extensively studied the square
cavity flow on a unit square domain to access the accuracy of the
numerical solution. The $u$-velocity on the vertical center line
$x=0.5$ and the $v$-velocity on the horizontal center line $y=0.5$
are given in Fig.~\ref{fig:cavity}d. For each sectional velocity,
typical regular(see, Figures~\ref{fig:cavity}a -~\ref{fig:cavity}c)
distributed nodes are, respectively, employed for comparison
purpose. It is shown that all the results obtained from the VIP
method are in good agreement with the data by Giha et
al.~\cite{ghia} which ha been a widely accepted reference for the
validation.

\subsection{Triangular cavity flow}
The two-dimensional steady incompressible flow inside a triangular
driven cavity is also an interesting subject like the square driven
cavity flow. This flow was studied analytically by
Moffatt~\cite{moffatt} in the Stokes regime. Moffatt showed that the
intensities of eddies and the distance of eddy centers from the
corner, follow a geometric sequence.

\begin{figure}[H]
        \centering
        \subfigure[]
        {
            \includegraphics[height=3.5cm]{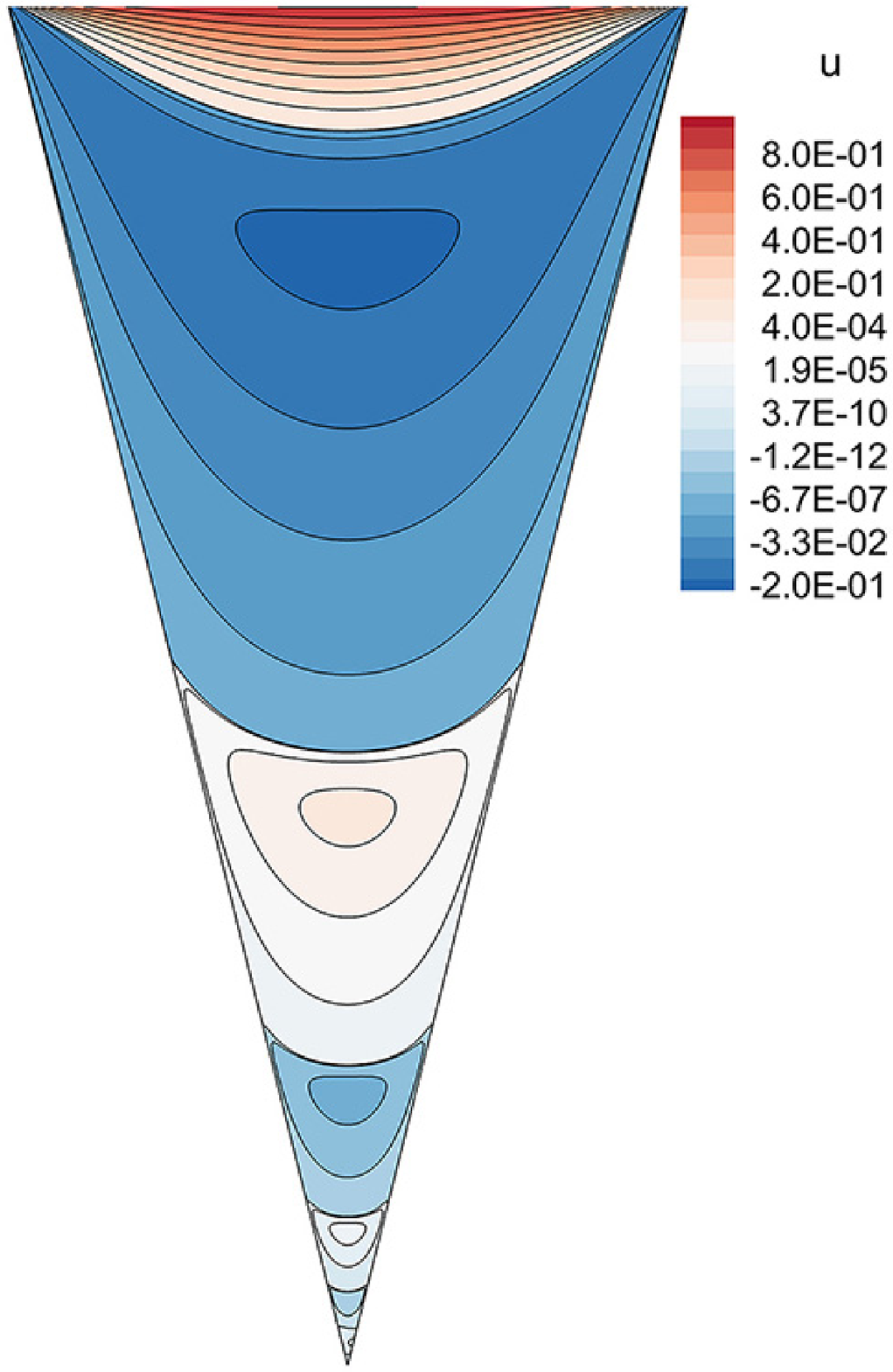}
        }
        \subfigure[]
        {
            \includegraphics[height=3.5cm]{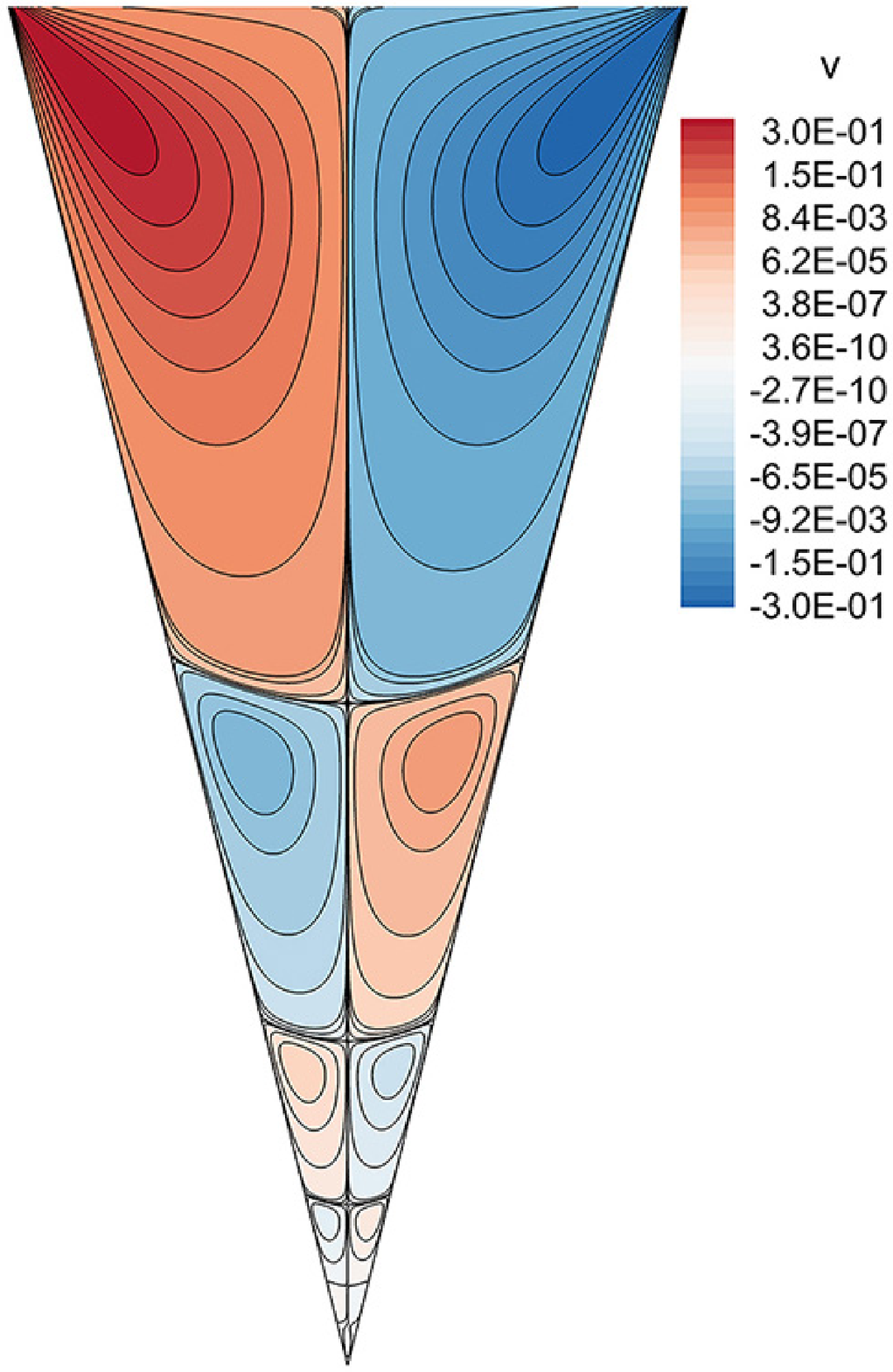}
        }
        \subfigure[]
        {
            \includegraphics[height=3.5cm]{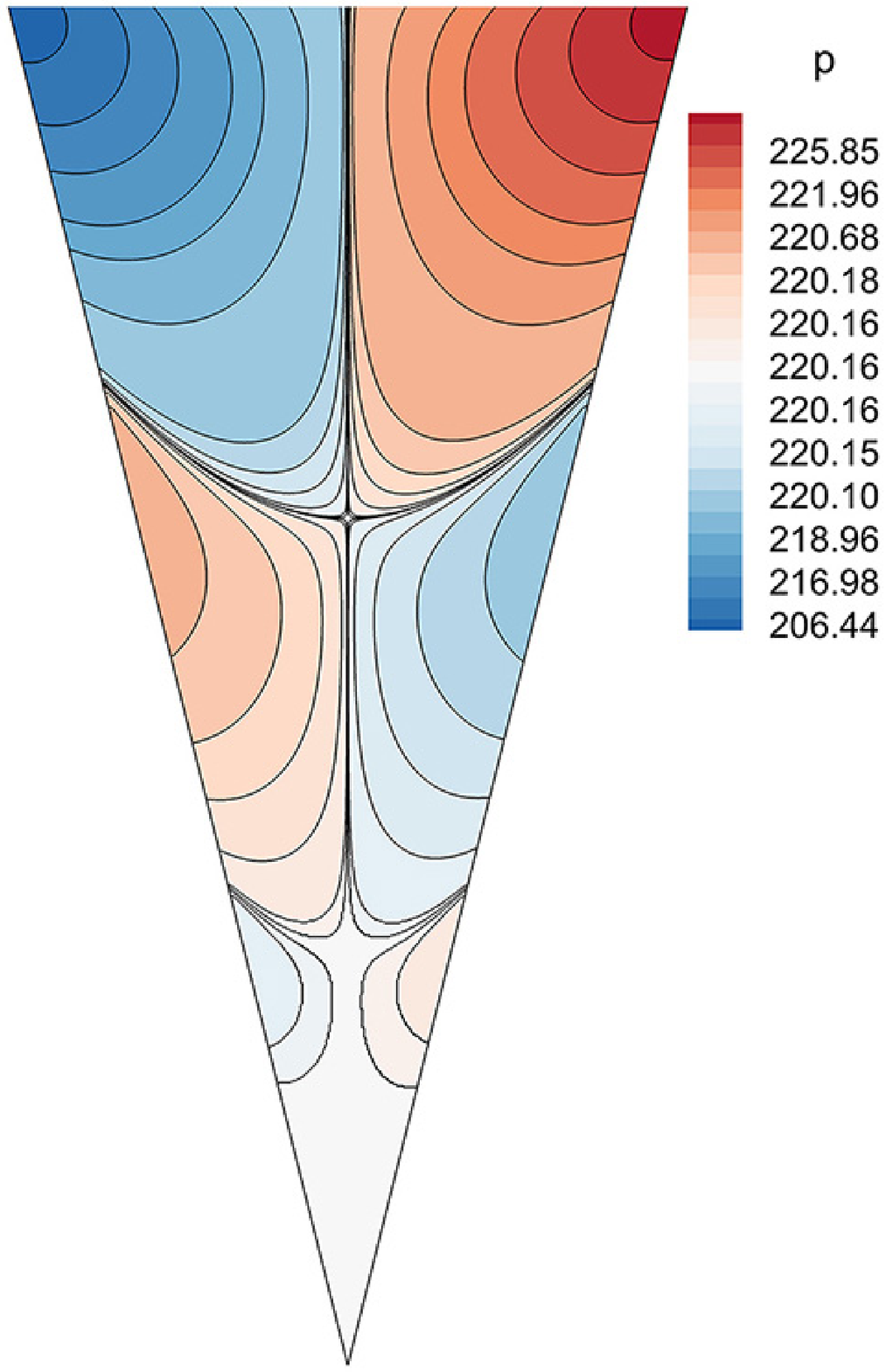}
        }
        \subfigure[]
        {
            \includegraphics[height=3.5cm]{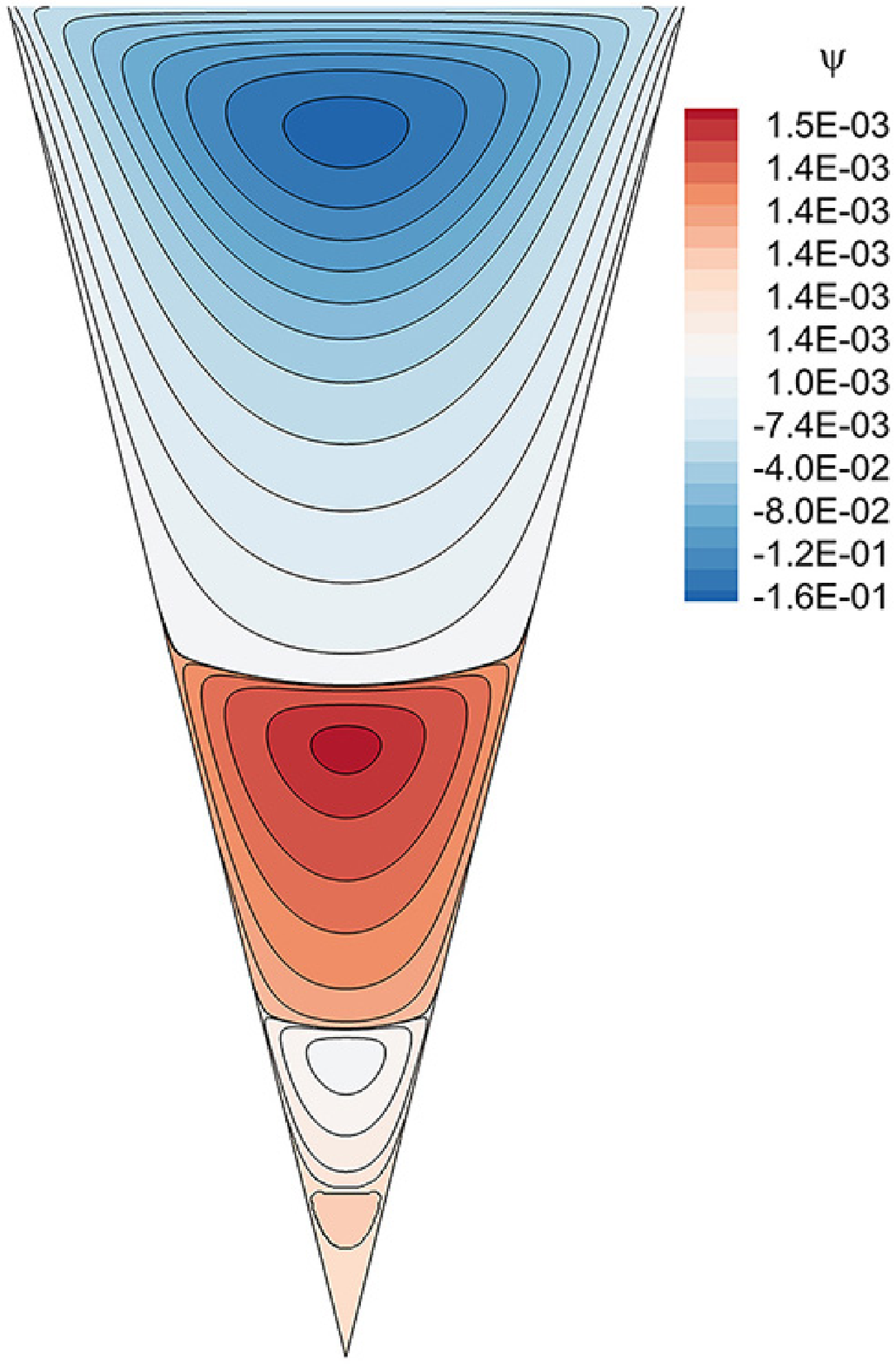}
        }
        \subfigure[]
        {
            \centering
            \includegraphics[height=3.5cm]{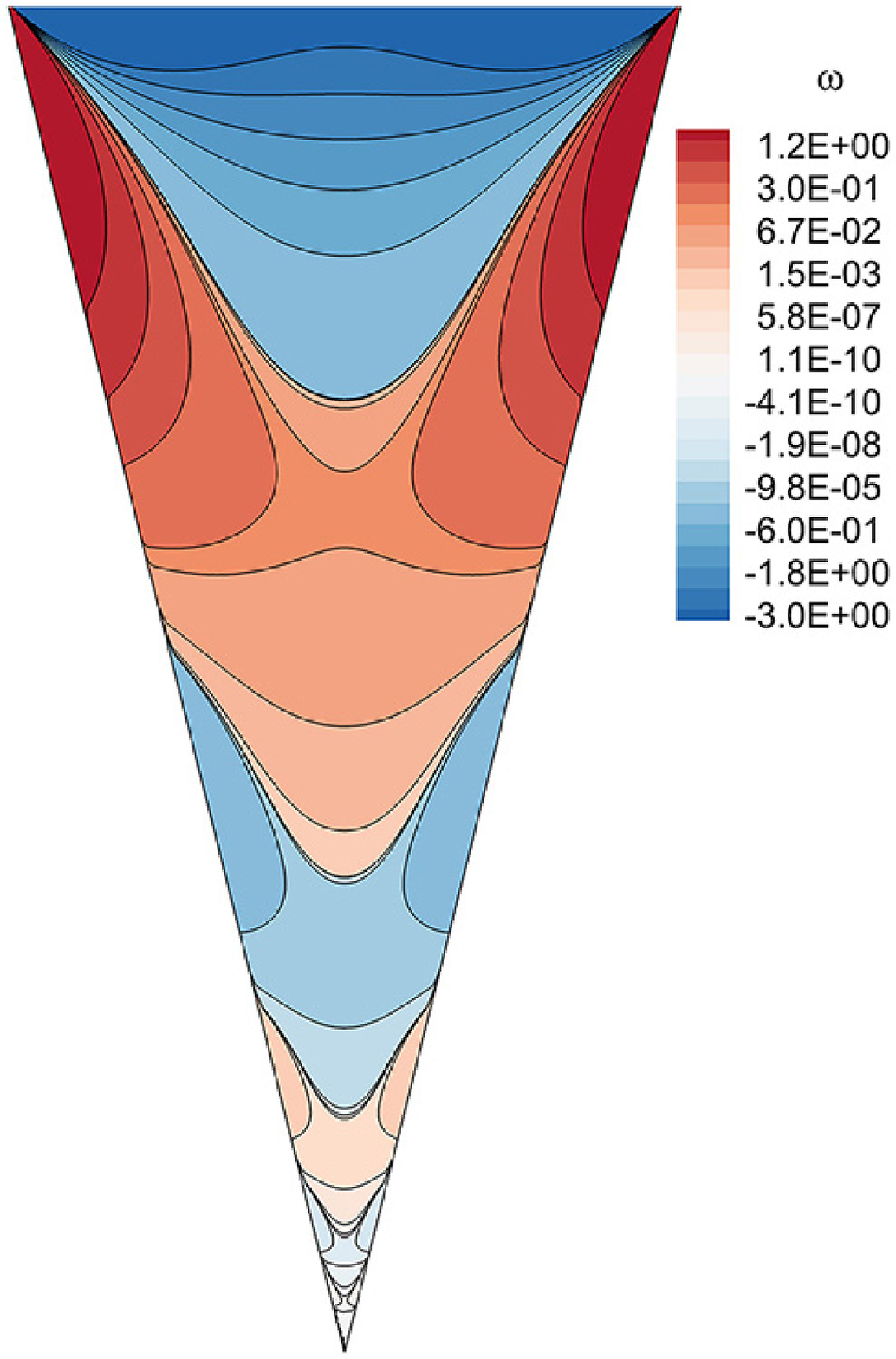}
        }
        \caption{Moffatt eddies toward the vertex of an isosceles triangle with $\theta = 28.072^\circ$ : (a) $u$-velocity; (b) $v$-velocity; (c) pressure; (d) stream function; (e) vorticity. }\label{fig:triangular cavity}
    \end{figure}

We apply the VIP method to this Moffatt eddy simulation for the flow
in a wedge-shaped domain(see, Figure~\ref{fig:triangular cavity
nodes}). In the velocity profile, the points where the $u$-velocity
has local maximum correspond to the $u$-velocity at the dividing
streamline between the eddies, where Moffatt have used these
velocities as a measure of the intensity of consecutive eddies.
Table~\ref{table:r} and ~\ref{table:I} tabulate the calculated
ratios of $r_n/r_{n+1}$ and $I_n/I_{n+1}$ for isosceles triangle
with $\theta = 28.072^\circ$ along with analytical predictions of
Moffatt and the agreement is good. The contour lines for u-velocity,
v-velocity, pressure, stream function, and vorticity are shown in
Fig.~\ref{fig:triangular cavity} where the sequence of eddies is
well presented. The changing signs of the velocity components are
properly illustrated toward the vertex of the wedge, which causes
the small eddies.

\begin{table}[H]
\centering
\tabcolsep=0.3cm
\begin{tabular}{cccccc}
\hline
 & $r_1/r_2$ & $r_2/r_3$ & $r_3/r_4$ & $r_4/r_5$ & $r_5/r_6$ \\
\hline\hline
VIP method($400 \times 800/2$) & 1.99 & 2.01 & 2.01 & 2.00 & 1.96 \\
Moffatt~\cite{moffatt} & \multicolumn{5}{c}{$r_n / r_{n+1} = 2.01$ } \\
\hline
\end{tabular}
\caption{Relative eddy center locations $r_n / r_{n+1}$ for
isosceles triangle with $\theta = 28.072^\circ$}\label{table:r}
\end{table}

\begin{table}[H]
\centering
\tabcolsep=0.3cm
\begin{tabular}{cccccc}
\hline
 & $I_1/I_2$ & $I_2/I_3$ & $I_3/I_4$ & $I_4/I_5$ & $I_5/I_6$ \\
\hline\hline
VIP method($400 \times 800/2$) & 385.7 & 406.0 & 402.3 & 388.1 & 380.9 \\
Moffatt~\cite{moffatt} & \multicolumn{5}{c}{$I_n / I_{n+1} = 407$} \\
\hline
\end{tabular}
\caption{Relative eddy center intensities $I_n / I_{n+1}$ for
isosceles triangle with $\theta = 28.072^\circ$}\label{table:I}
\end{table}

\begin{figure}[H]
        \centering
        \subfigure[]
        {
            \includegraphics[height=6cm]{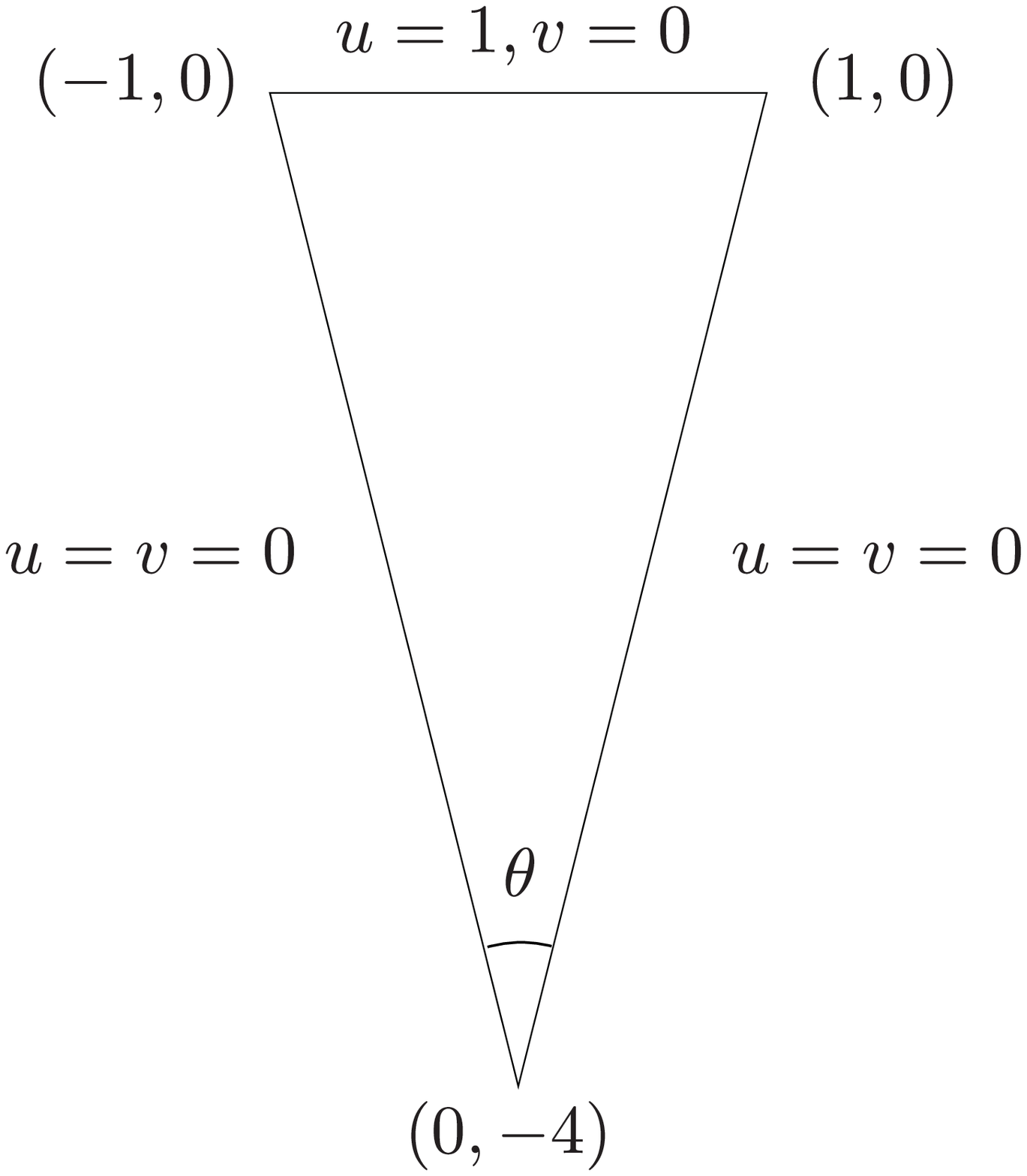}
        }
        \qquad
        \subfigure[]
        {
            \includegraphics[height=6cm]{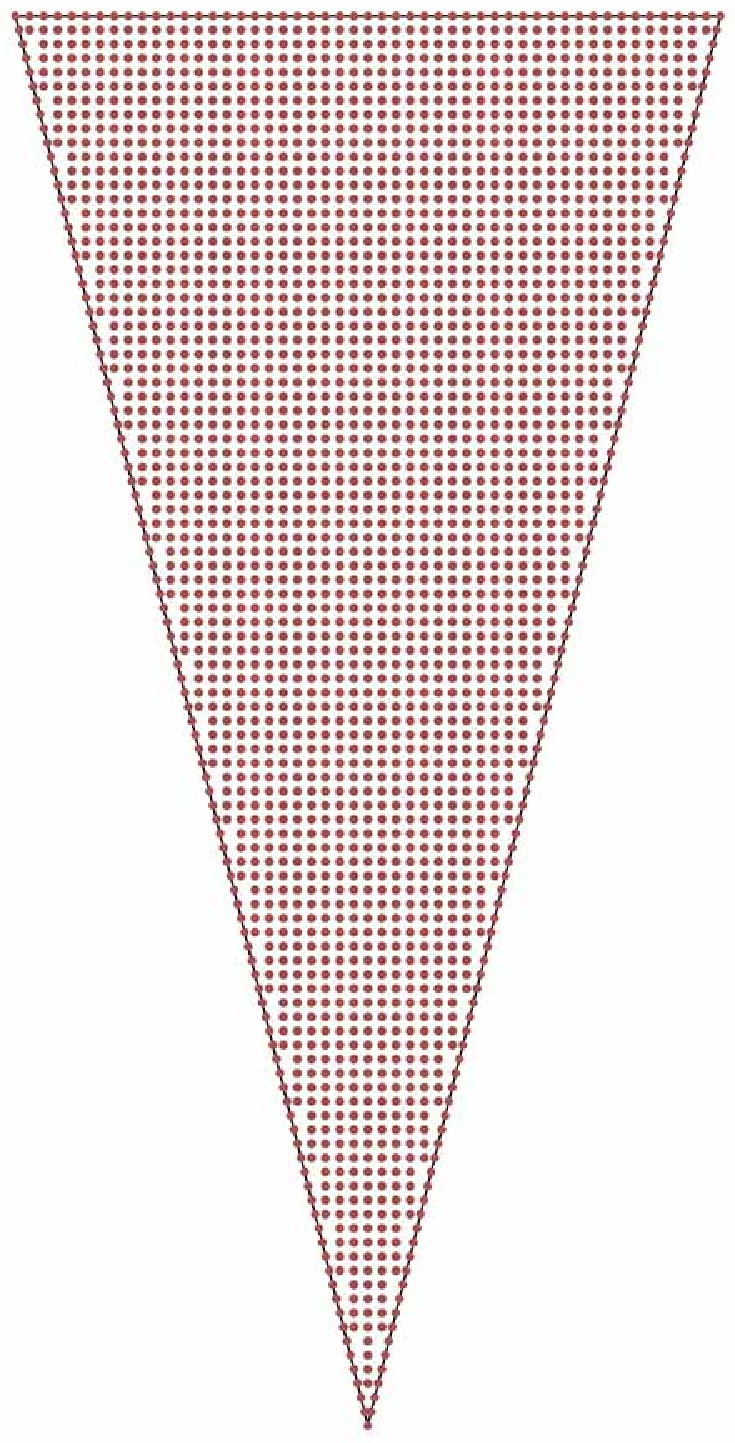}
        }
        \caption{ Isosceles triangle cavity flow with $\theta = 28.072^\circ$ : (a) Problem
        description; (b) regular distributed nodes.}\label{fig:triangular cavity nodes}
    \end{figure}

\subsection{Flow past a circular cylinder: Steady and Unsteady}
Flow past a circular cylinder is one of the classical problems of
fluid mechanics. For lower value of Reynolds number, the flow is
steady and symmetric. And as the Reynolds number is increased, the
flow past a circular cylinder is a problem unsteady in nature and,
therefore, good numerical accuracy is required in order to capture
the different phenomena present in the evolving solution.

\begin{figure}[H]
\centering \mbox{
        \subfigure[]
        {
            \includegraphics[height=5cm]{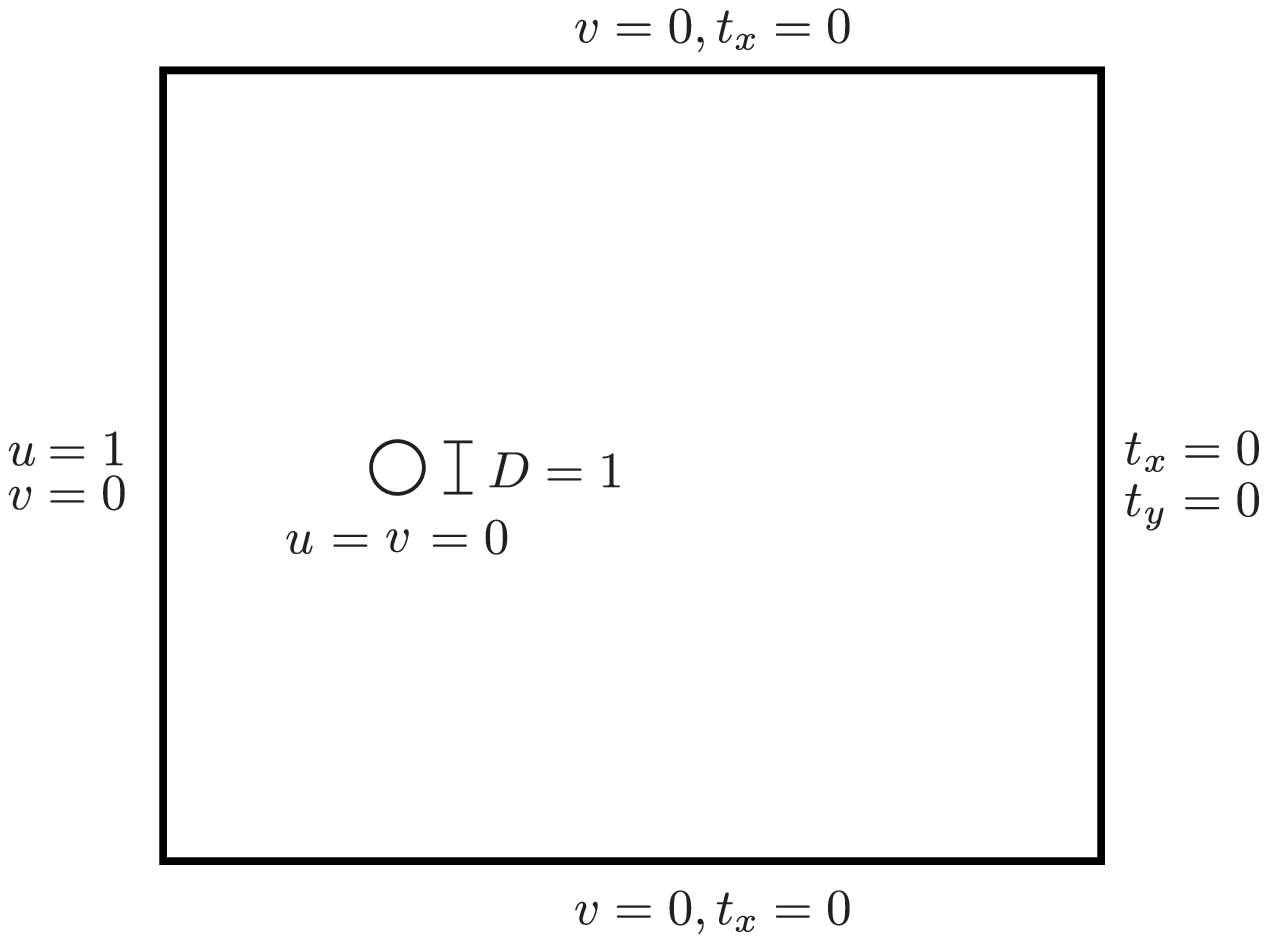}
        }
        \subfigure[]
        {
            \includegraphics[height=5cm]{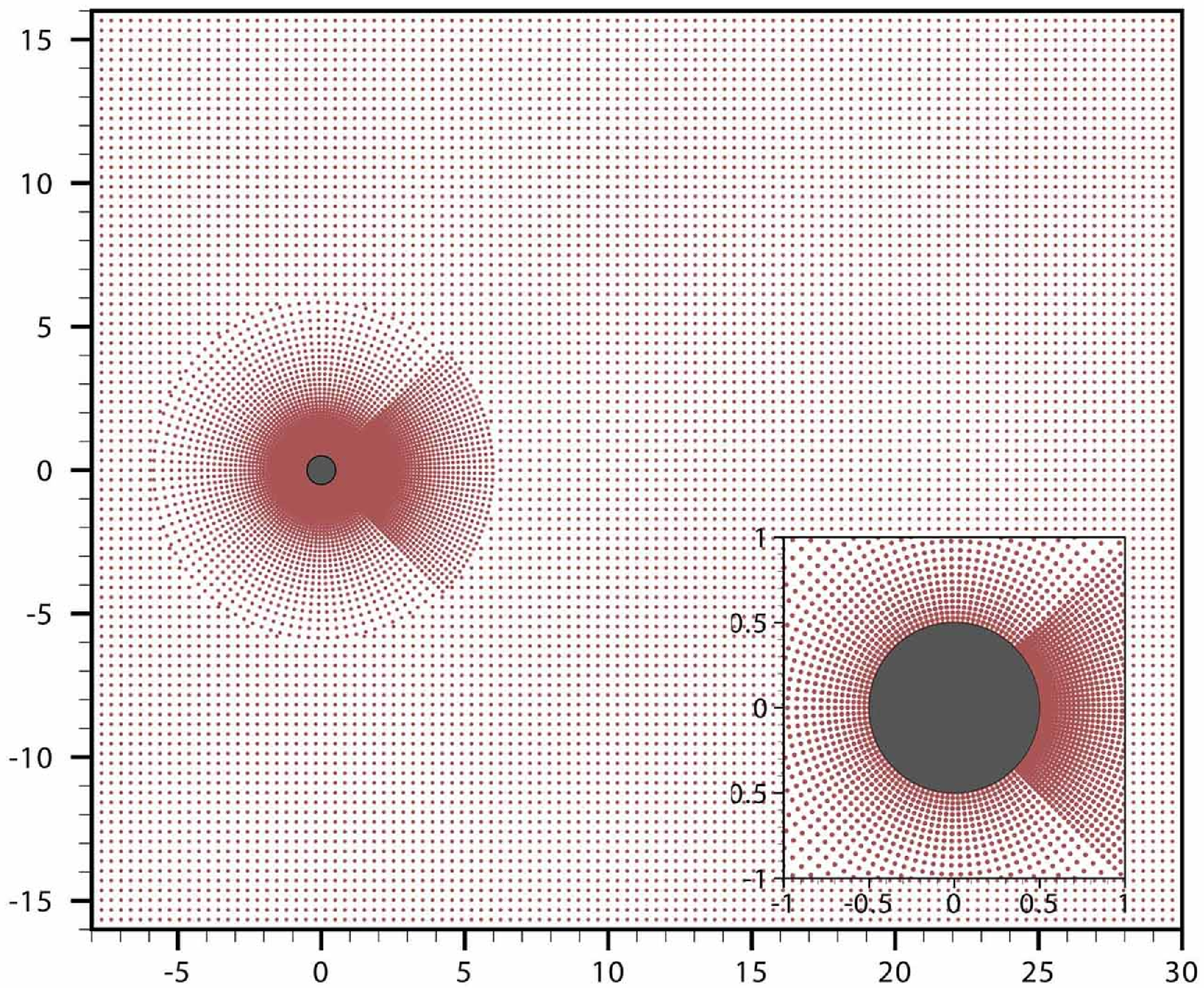}
        }
} \caption{Flow past a circular cylinder : (a) Problem
        description; (b) regular distributed nodes and close-up of the central block containing the
cylinder. }\label{fig:circular_cylinder}
\end{figure}

To validate the VIP scheme, the numerical simulation of the steady
and unsteady flows past a circular cylinder is carried out. In the
problem under investigation, depicted in
Fig.~\ref{fig:circular_cylinder} along with the computational domain
and the distributed nodes. The Reynolds number in the this flow is
defined as $ Re = u_{\infty} D / \nu $ where $D$ is the diameter of
the cylinder. We impose $u$ = 1 and $v$ = 0 for inlet, and traction
free condition, $(t_x, t_y)$ = 0 for outlet, and $v$ = 0 and $t_x$ =
0 for top and bottom boundary. The traction vector
$\bf{t}$ is defined by $\mathbf{t} = -p \mathbf{n} + \nu
\partial \mathbf{u}/\partial \mathbf{n}$ with $\mathbf{n}$ denoting
the outward normal. At low Reynolds numbers, the flow develops two
symmetric wakes past the cylinder. This solution becomes unstable
for Reynolds numbers over $40$, and periodic vortex shedding
appears. These vortices are transported by the flow, creating what
is known in the literature as Von Karman vortex street.

\begin{figure}[H]
 \centering  \mbox{
        \subfigure[]
        {
            \includegraphics[height=3.5cm]{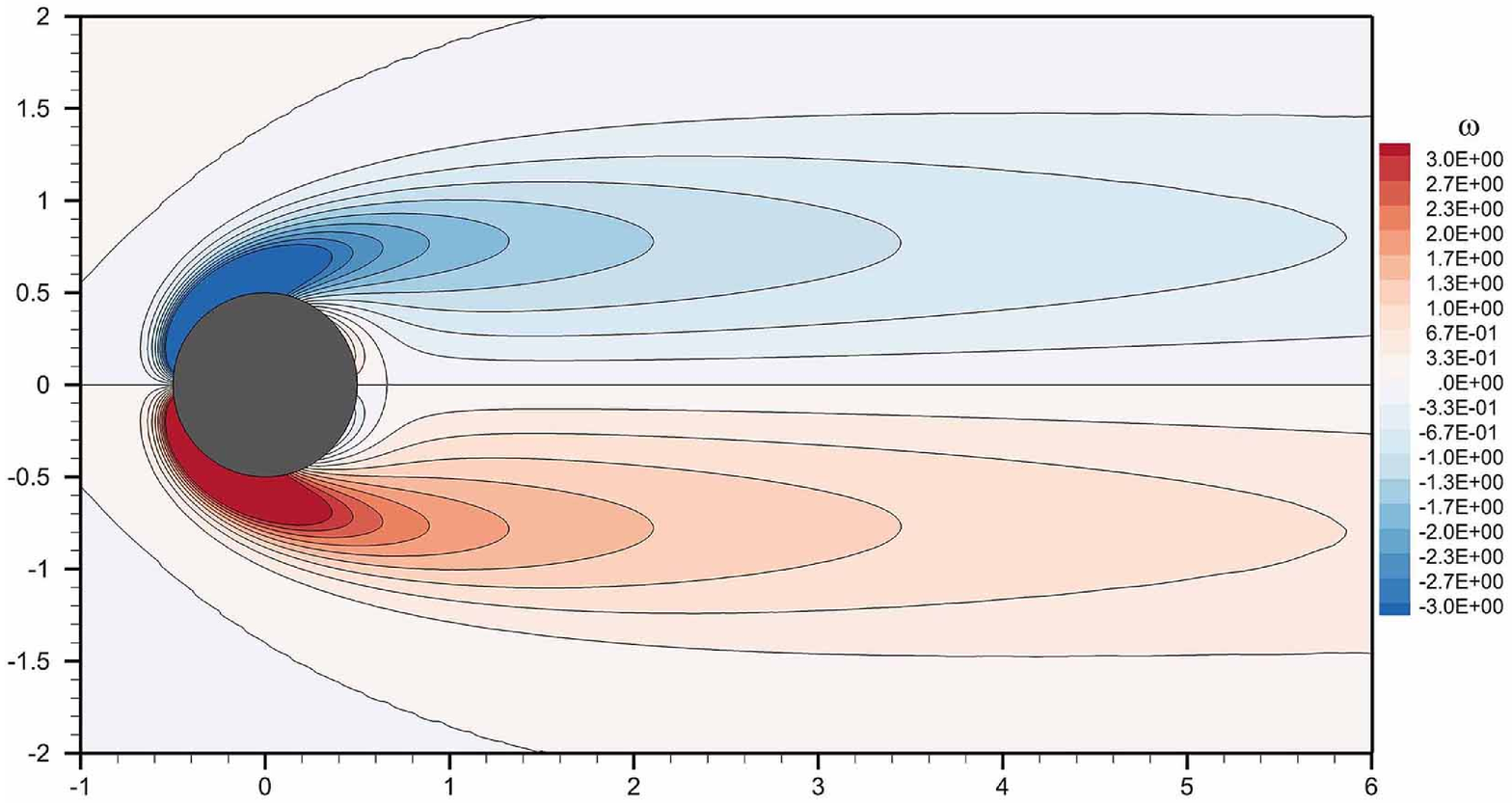}
        }
        \subfigure[]
        {
            \includegraphics[height=3.5cm]{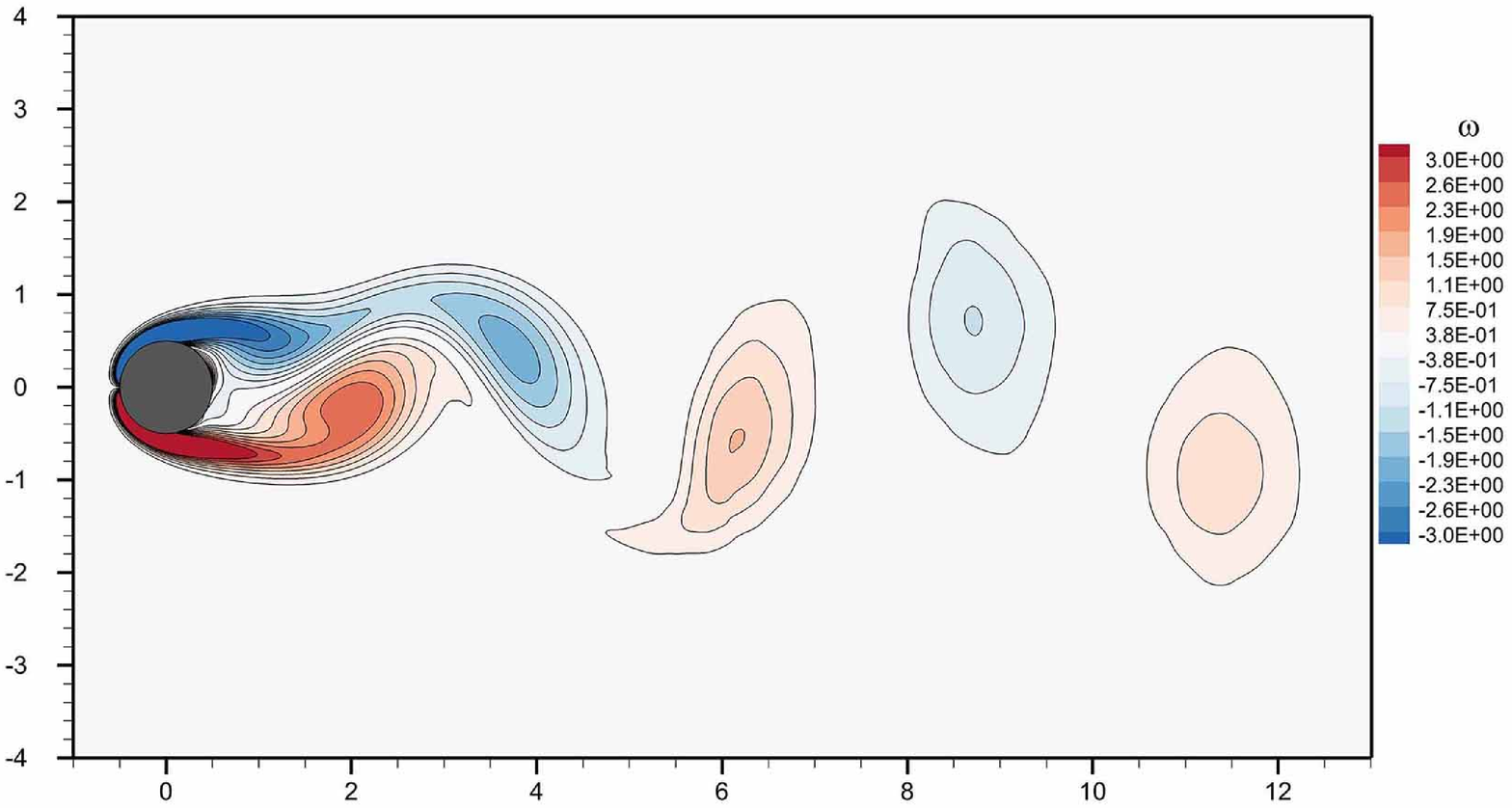}
        }
}
        \caption{Spanwise vorticity contours near a circular cylinder : (a)
$Re$ = 40; (b) $Re$ = 100, instantaneous vorticity contours are
    drawn. }\label{fig:circular_cylinder vorticity}
\end{figure}

\begin{figure}[H]
        \centering
        \subfigure[]
        {
            \includegraphics[height=4.5cm]{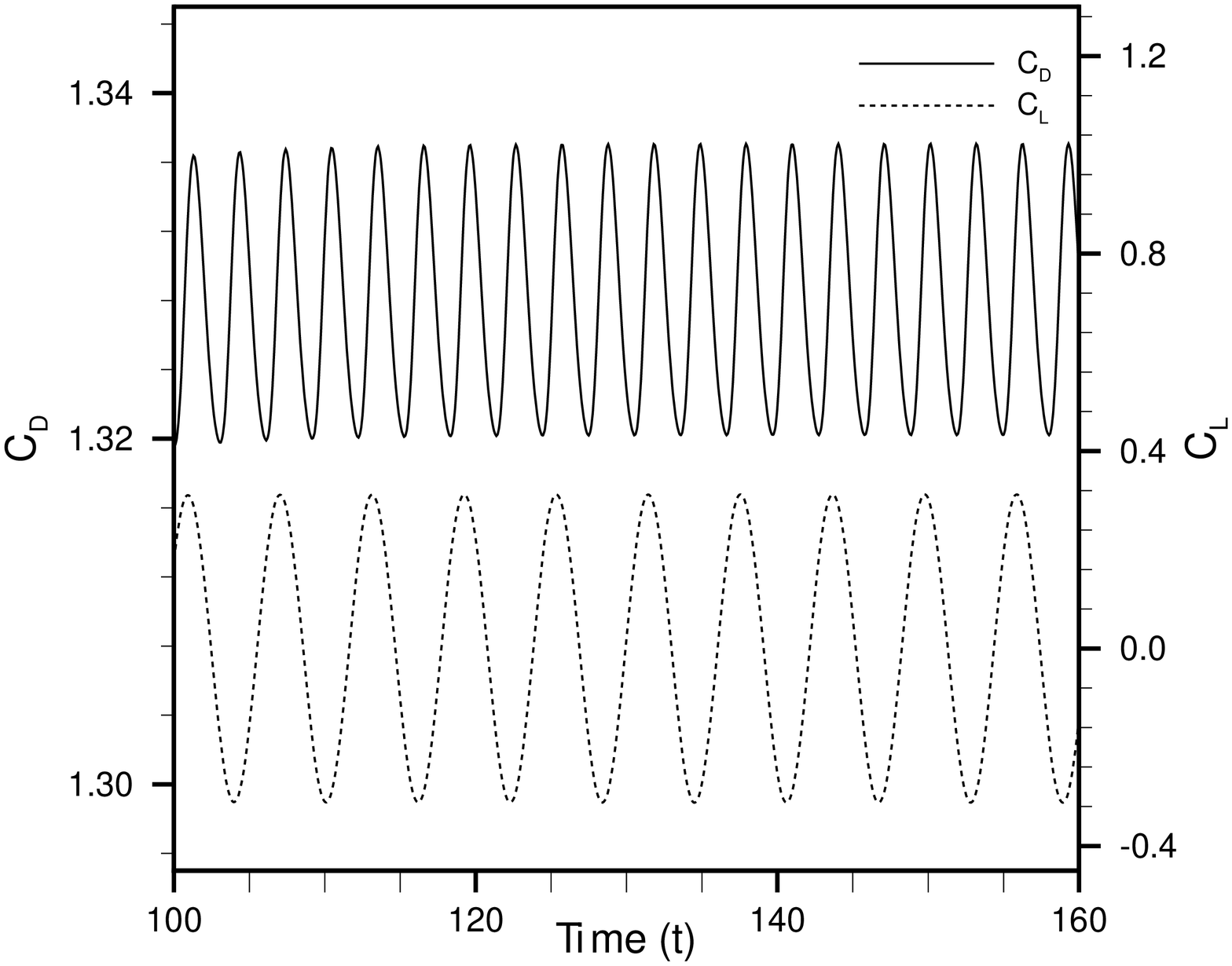}
        }
        \subfigure[]
        {
            \includegraphics[height=4.5cm]{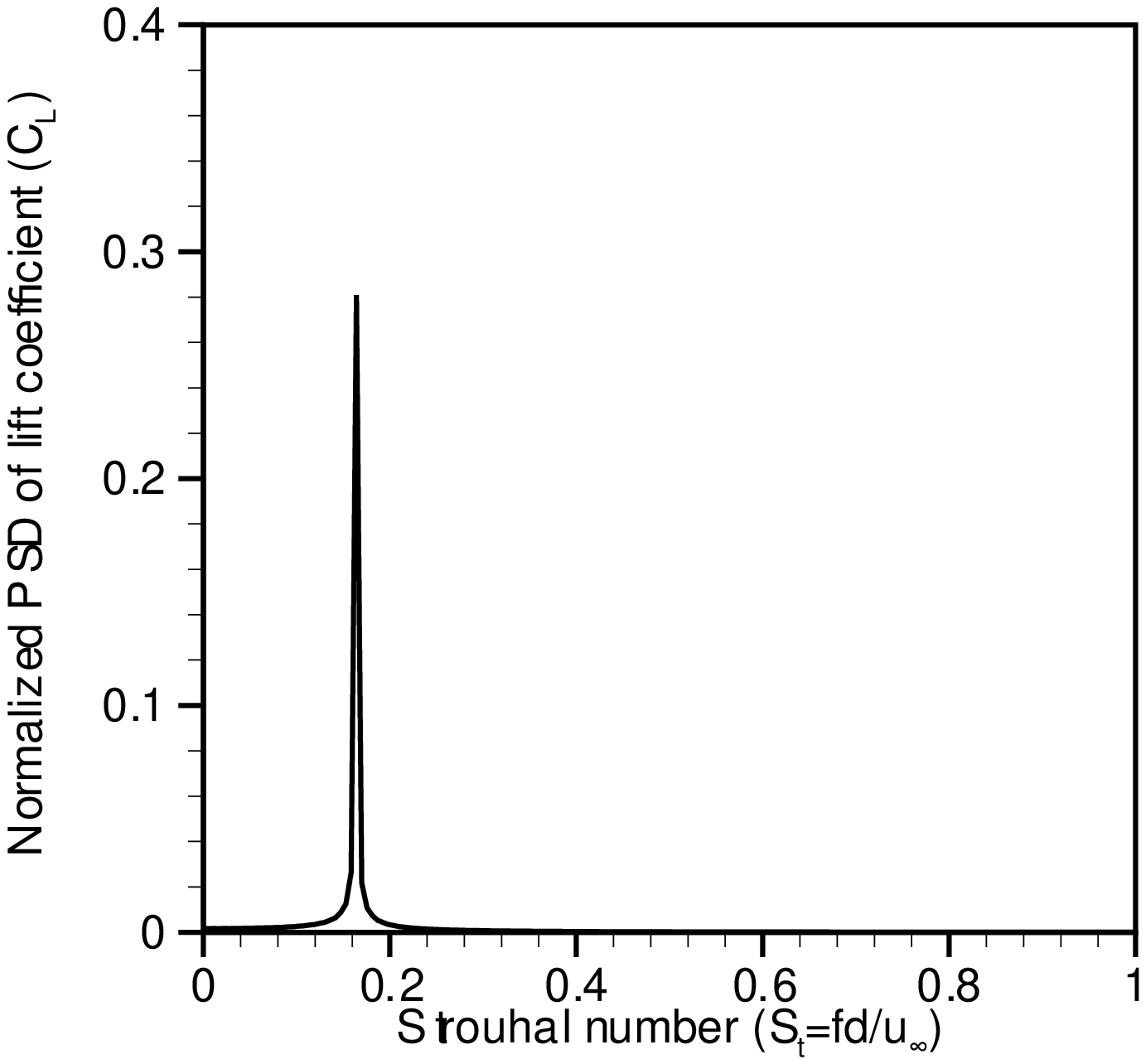}
        }
        \caption{Results of flow past a circular cylinder with $Re=100$: (a) Temporal evolution of lift coefficient and drag coefficient; (b) Normalized power spectral density (PSD) for lift coefficient. The Strouhal number ($S_t$) is indicated by the summit of the curve ($S_t = 0.164)$). } \label{fig:circular_cylinder CdClSt}
\end{figure}

Figure~\ref{fig:circular_cylinder vorticity} shows the spanwise
vorticity contours for $Re=40$ and $Re=100$.
Table~\ref{table:Re40100} and Figure~\ref{fig:circular_cylinder
CdClSt} show the results of simulations together with the previous
numerical results of Takami $\textit{et al.}$~\cite{takami} and H.
Ding~\cite{ding}, where $C_D$ is the drag coefficient(time-averaged
value in case of $Re$=100) and $C^{\prime}_L$ is the amplitude of
lift-coefficient fluctuations (maximum deviation from the
time-averaged value) at $Re$=100. The Strouhal number
($St=fd/u_{\infty}$) is also in excellent agreement with numerical
results(see, Table 3 and Figure~\ref{fig:circular_cylinder CdClSt}).

\begin{table}[H]\begin{center}
    {\tabcolsep=0.2cm
        \begin{tabular}{llllllllllllllllll}
            \hline
            & & $Re$ &  & $C_D$ & &  $C^{'}_L$ & & $St$ \\
            \hline\hline
VIP method   & & 40  & & 1.536 &&  &&  \\
             & & 100 & & 1.328 && $\pm$ 0.31 && 0.164 \\
\\
Takami $\textit{et al.}$~\cite{takami} & & 40 && 1.536 &&  &&  \\
H. Ding~\cite{ding}  & & 100     && 1.325 && $\pm$ 0.28 && 0.164 \\
\hline
        \end{tabular}}
        \caption{Simulation results for flow over a circular
        cylinder.}\label{table:Re40100}
        \end{center}
    \end{table}

\subsection{Flow past a bumpy circular cylinder}
Analytic solutions are rarely available, and conventional numerical
computations are usually out of reach since the rapidly varying
wrinkles and the domain have different length scales. The
traditional remedy is to pose special boundary conditions on a
mollified domain to capture the geometrical influence of the
wrinkles. The development of such conditions is cumbersome in
general, and modeling error estimates can be out of reach.

Nevertheless, many problems involving oscillating boundaries or
interfaces arise in many fields of physics and engineering sciences,
such as the scattering of acoustic waves on small periodic
obstacles, the free vibrations of strongly nonhomogeneous elastic
bodies, the behavior of fluids over rough walls.

Interesting example in the fluid mechanics is the flow field around
golf balls, in which the wrinkles associated to the curvature
decrease the gap between the air-pressure behind and in front of the
ball.

\begin{figure}[H]
        \centering
        \subfigure[]
        {
            \includegraphics[width=0.55\textwidth]{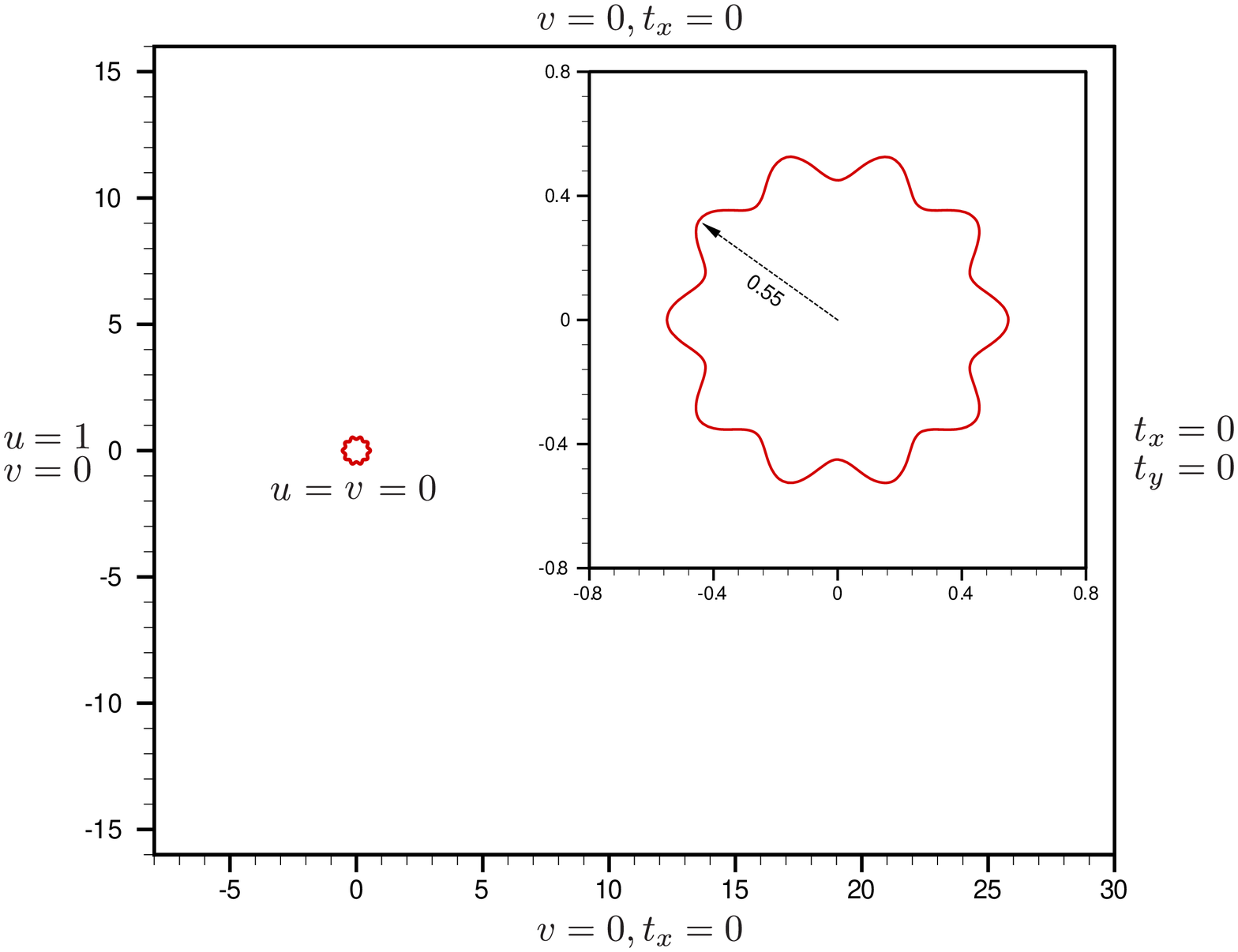}
        }\\
        \subfigure[$m_b=10$.]
        {
            \includegraphics[width=0.136\textwidth]{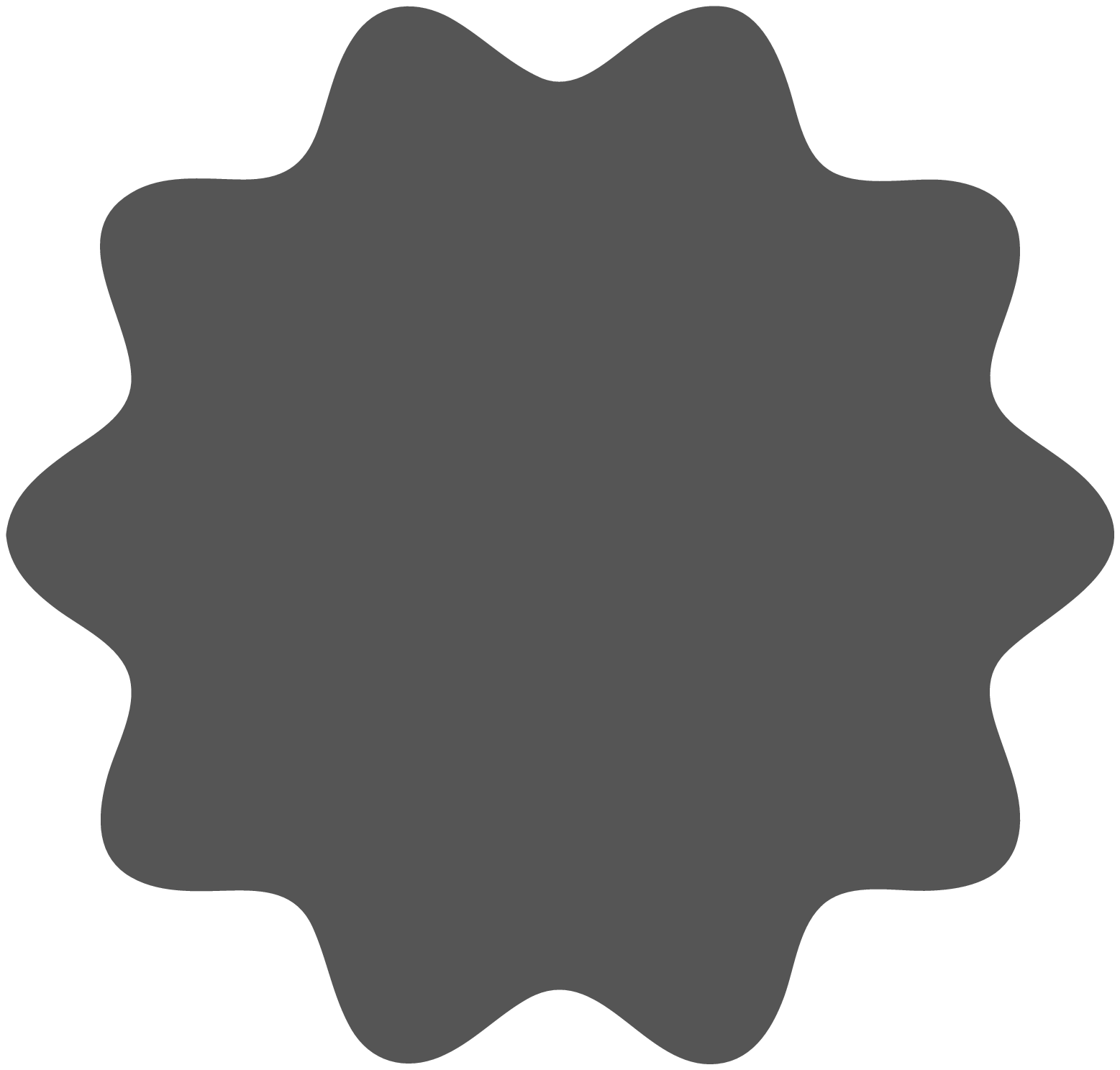}
        }\qquad
        \subfigure[$m_b=20$.]
        {
            \includegraphics[width=0.136\textwidth]{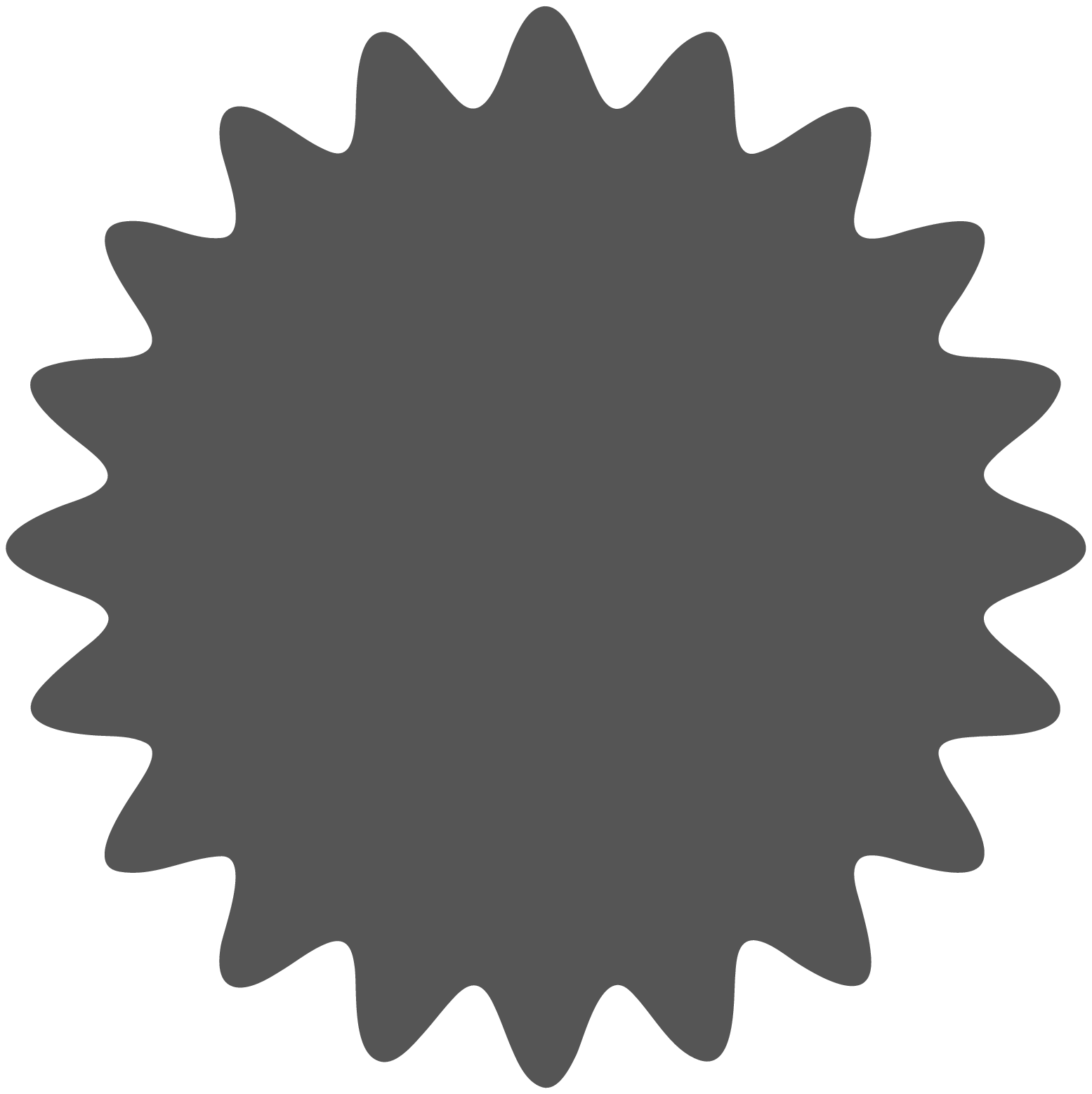}
        }\qquad
        \subfigure[$m_b=30$.]
        {
            \includegraphics[width=0.136\textwidth]{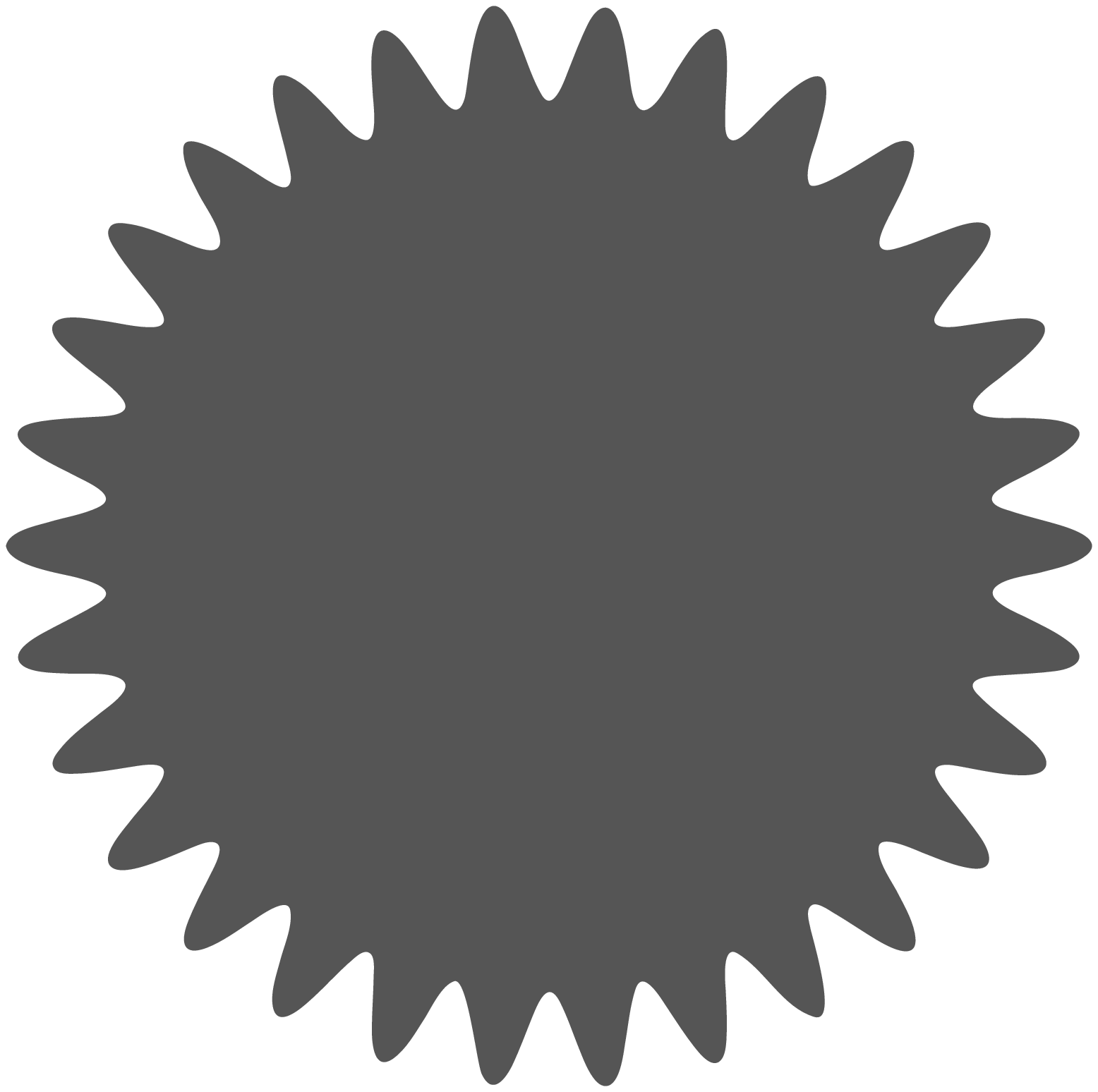}
        }\\
        \subfigure[$m_b=40$.]
        {
            \includegraphics[width=0.136\textwidth]{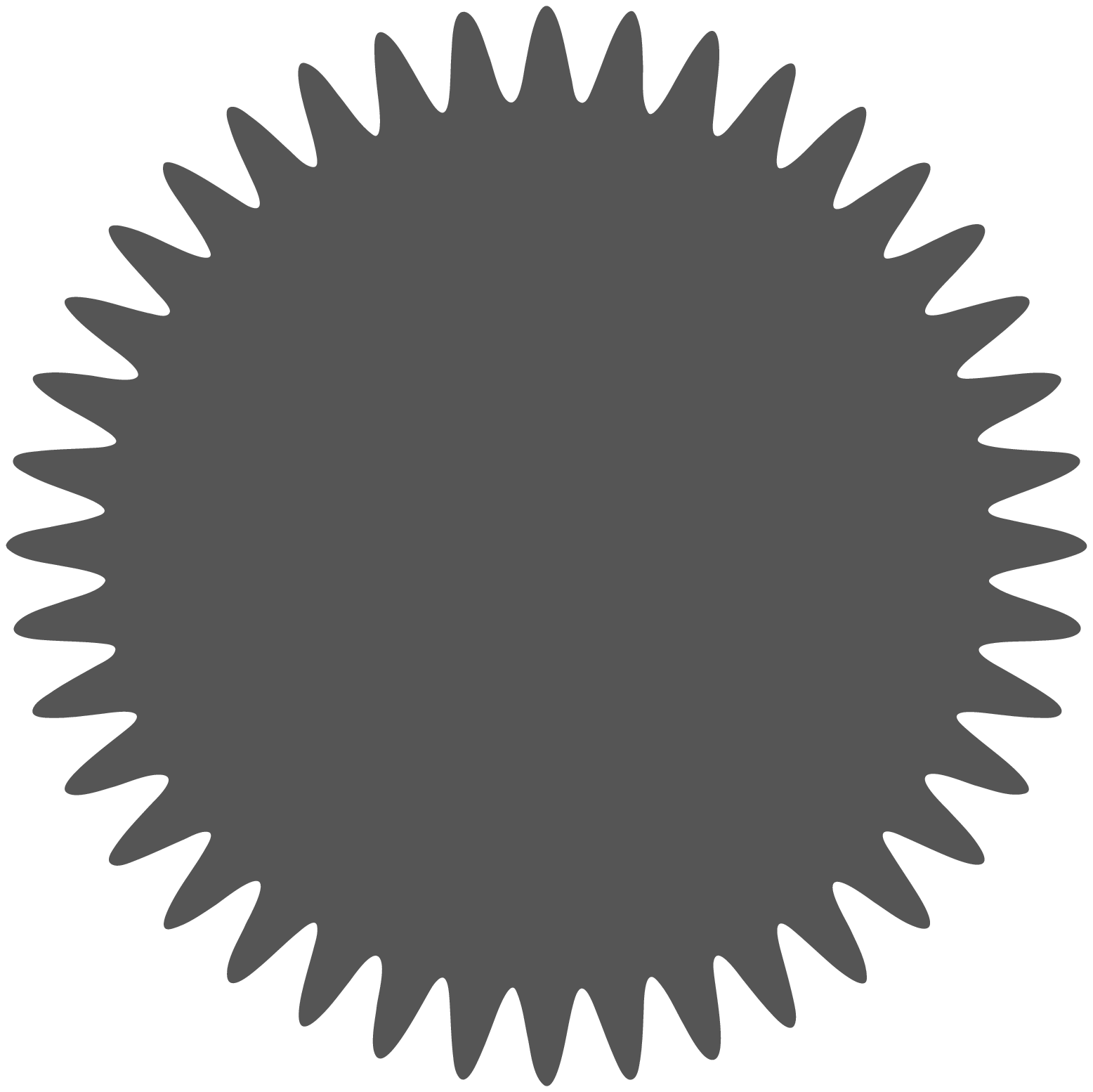}
        }\qquad
        \subfigure[$m_b=50$.]
        {
            \includegraphics[width=0.136\textwidth]{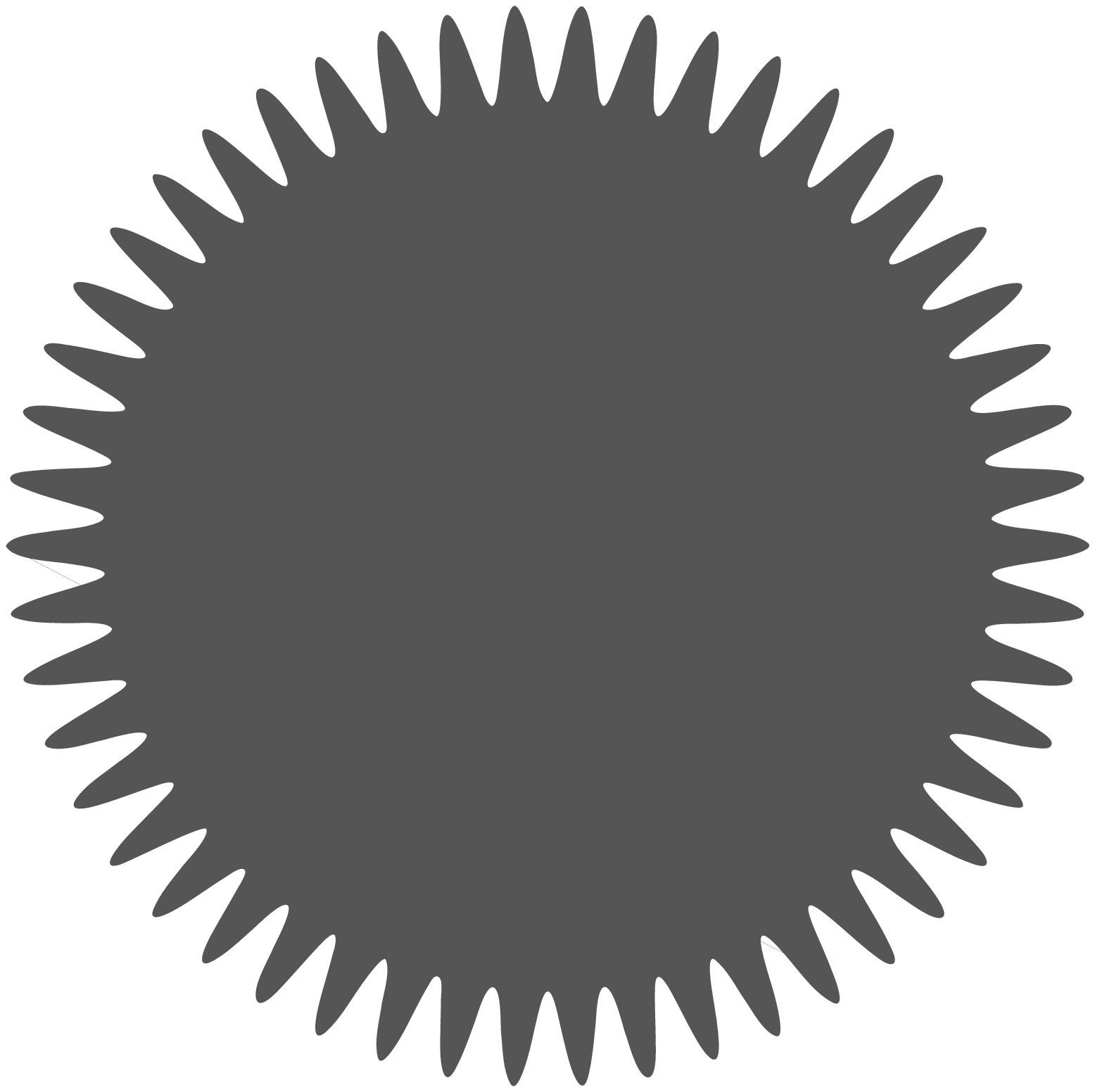}
        }\qquad
        \subfigure[$m_b=60$.]
        {
            \includegraphics[width=0.136\textwidth]{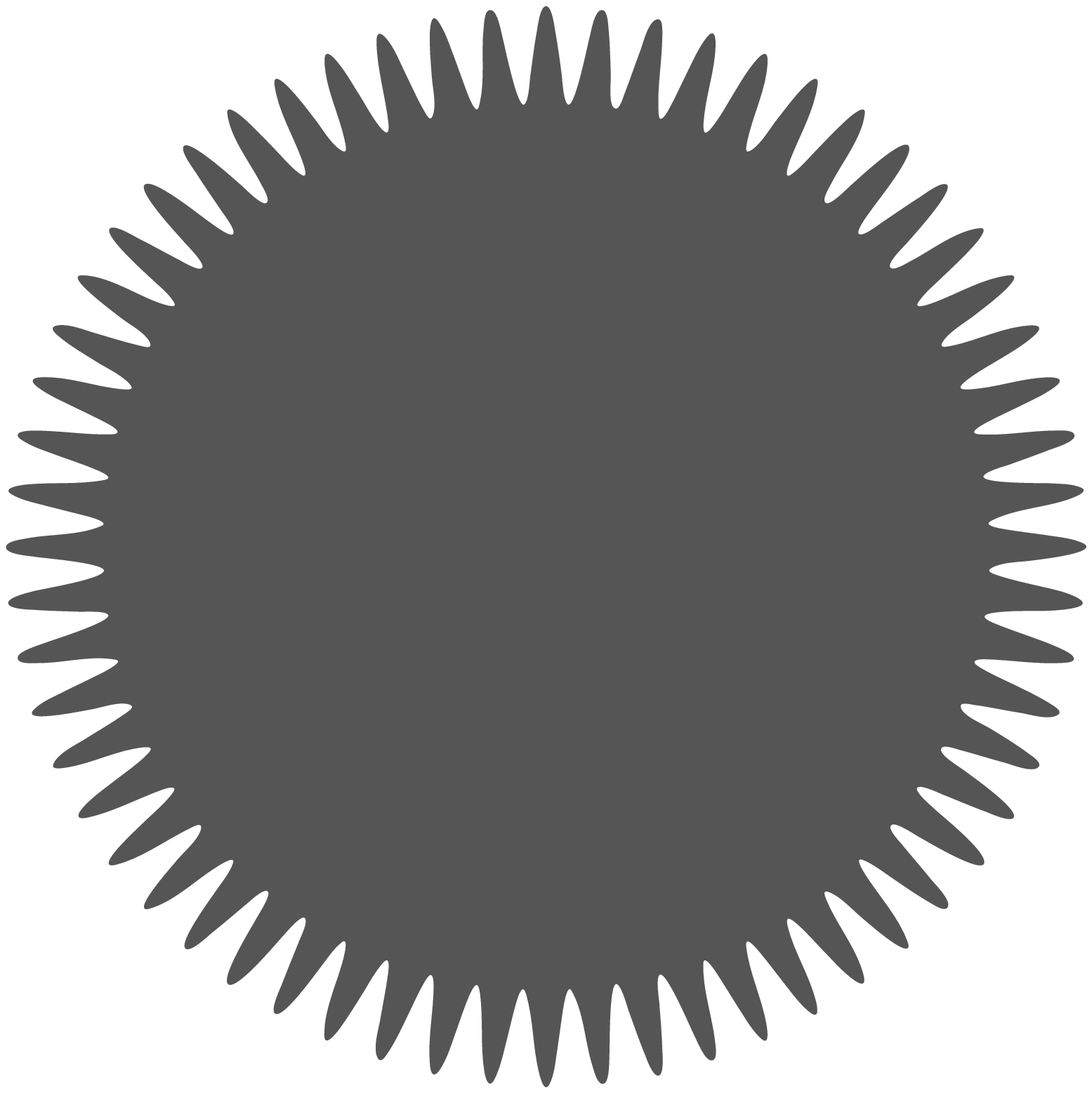}
        }\\
        \subfigure[$m_b=70$.]
        {
            \includegraphics[width=0.136\textwidth]{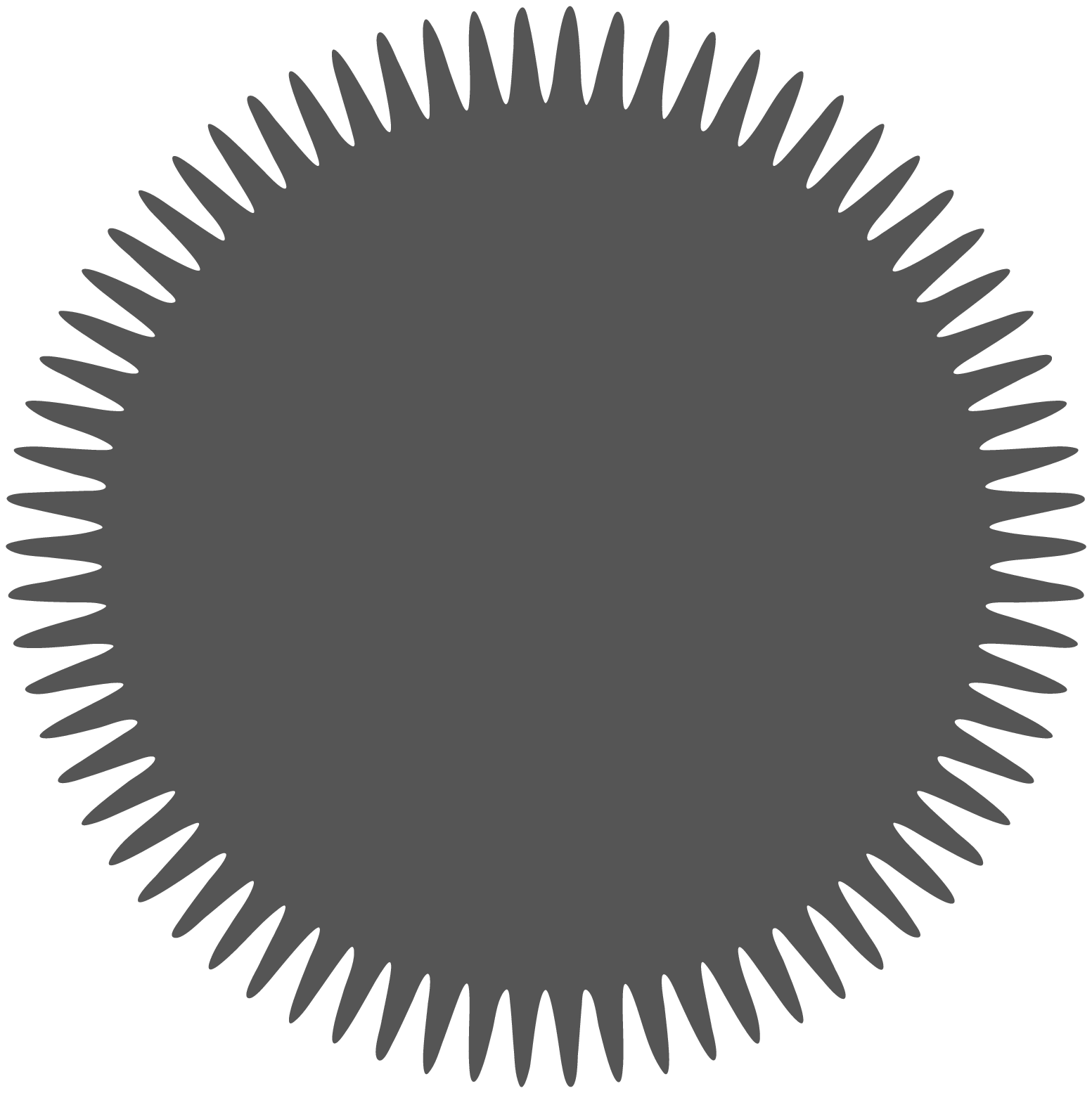}
        }\qquad
        \subfigure[$m_b=80$.]
        {
            \includegraphics[width=0.136\textwidth]{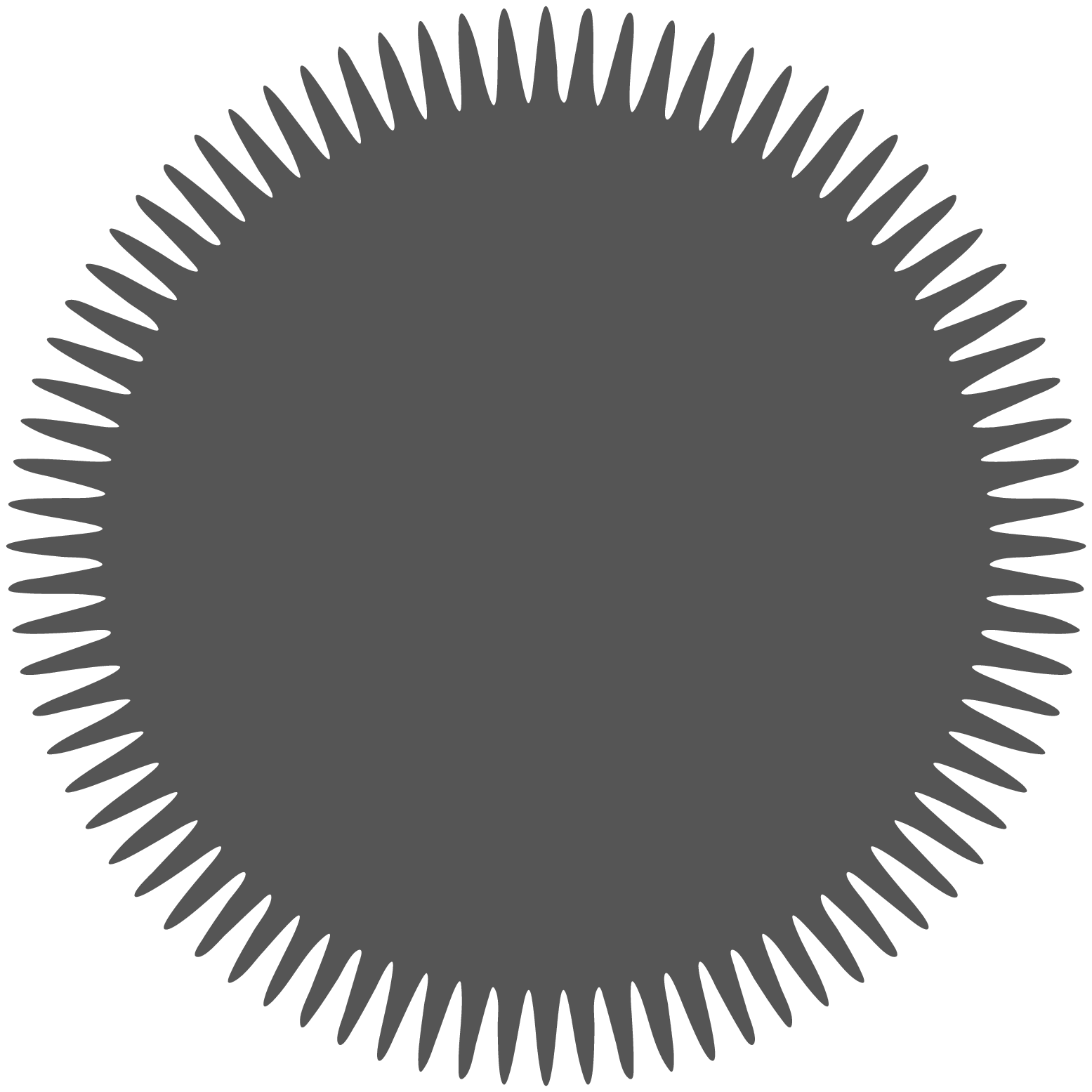}
        }\qquad
        \subfigure[$m_b=90$.]
        {
            \includegraphics[width=0.136\textwidth]{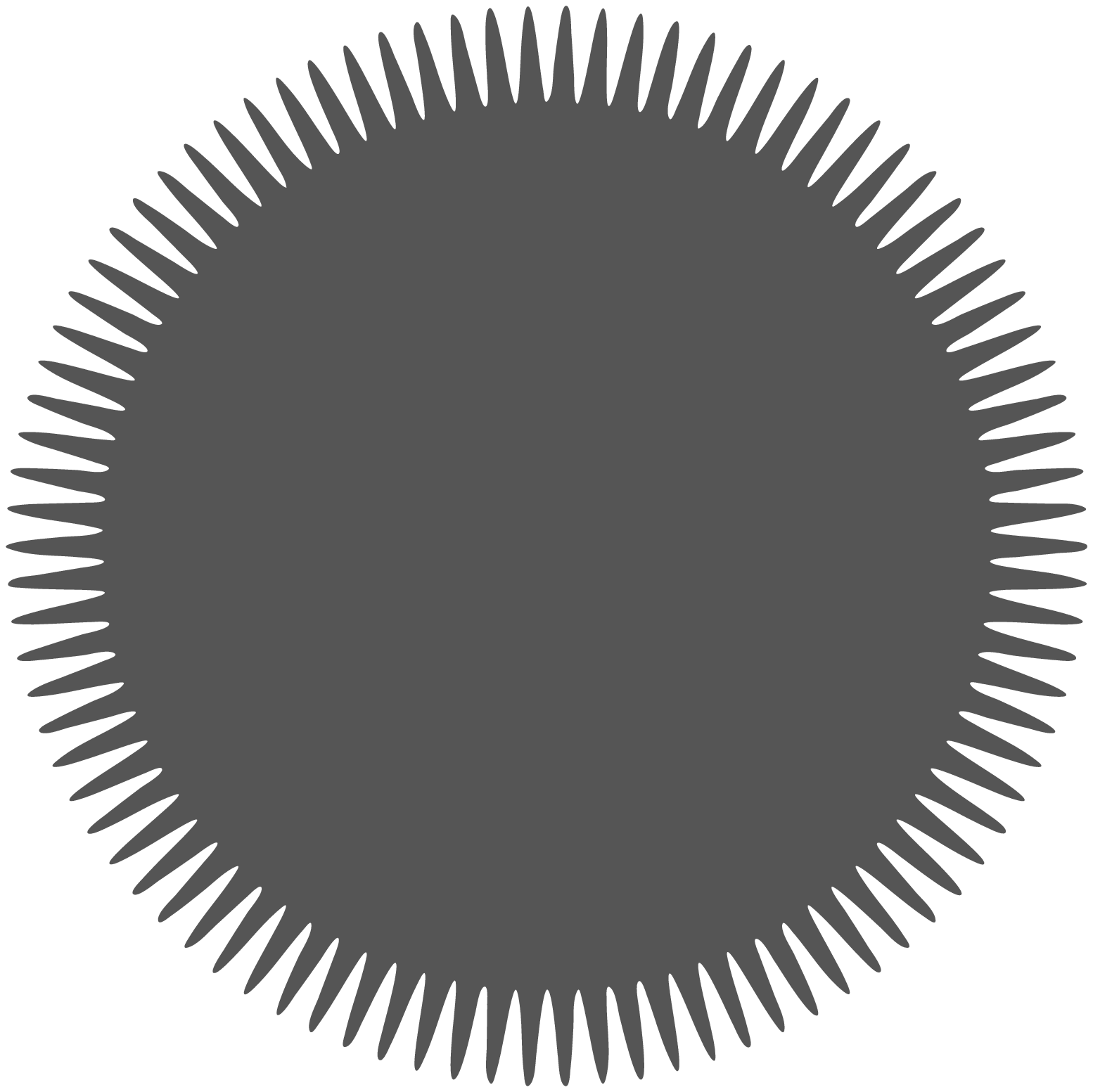}
        }\\
        \subfigure[$m_b=100$.]
        {
            \includegraphics[width=0.136\textwidth]{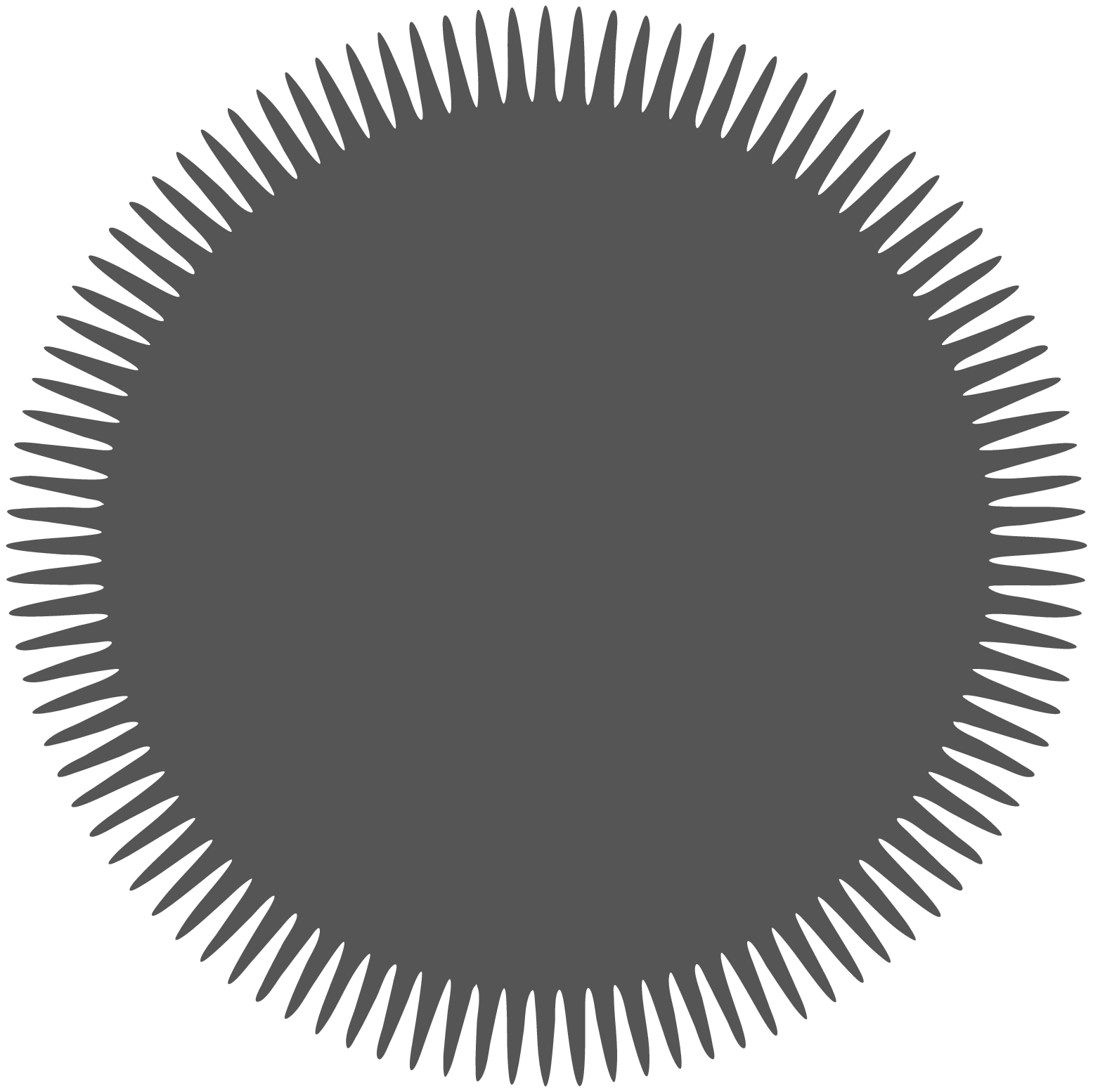}
        }
        \caption{The numerical setting for flow past various bump
        cylinders: (a) the computational domain with the external dimension of $(-8,30) \times (-16,16)$; (b)-(k) are various shape of bumpy circles. }\label{fig:bumpy cylinder}
    \end{figure}

\begin{figure}[H]
        \centering
        \subfigure[]
        {
            \includegraphics[width=0.65\textwidth]{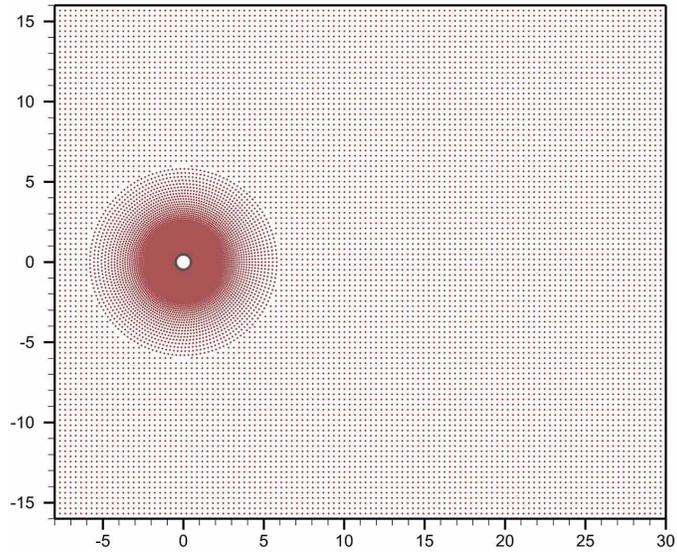}
        }\\
        \subfigure[$m_b=100$]
        {
            \includegraphics[width=0.65\textwidth]{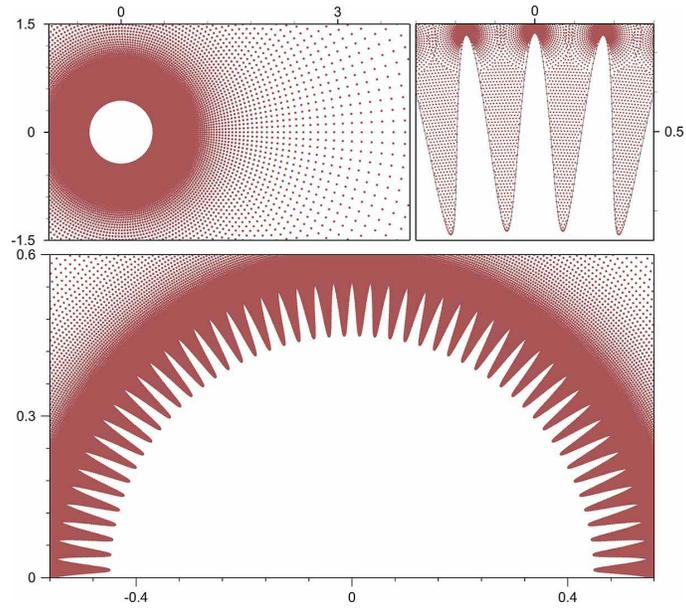}
        }
        \caption{The numerical feather for flow past a bumpy
        cylinder($m_b$=100): (a) the node distribution($236,604$ nodes) of the full domain for the computations of the external flow;
        (b) the close up and the zoom-in views of the bumpy circle and the radial stepwise uniform node distribution with the radius of R=0.55}\label{fig:bumpy cylinder distributed nodes}
    \end{figure}

As depicted below, we represent the surface by a bumpy circle
$\textbf{x}(\theta)$ of which the radial deviation is introduced by
the following sinusoidal curves:

\begin{align*}
        \textbf{x}(\theta) &= r\left[1+\gamma \cos (m_b \theta)\right] (\cos\theta, \sin\theta),
\end{align*}
where $r$ is a dimensionless average radius of the bumpy circle,
$\gamma$ denotes the amplitude ratio of a bump, and $m_b$ is the
total number of bumps along the circumference. $r$=0.5, $\gamma$=0.1
and $m_b=10,20,30,\cdots,90,100$ are used in our simulations(see,
Figure~\ref{fig:bumpy cylinder}). The present method never employs
any meshes, grids, or even integration cells, i.e., being entirely
free from connectivity data. Thus it has an advantage over other
numerical methods that are based on subdivisions in modelling this
complex geometry of bumpy circular cylinder, as can be noticed from
(b) in Fig.~\ref{fig:bumpy cylinder distributed nodes}. Stepwise
adaptive node distribution is also employed similarly as introduced
in the previous example of flow around a circular cylinder. For all
cases, we use a fixed value of $Re = 40$. Stepwise adaptive node
distribution is also employed similarly as introduced in the
previous example of flow around a circular cylinder.

The aim of this example is to investigate the correlation between
the drag coefficient and the shape of bumpy circles. With the
increase of bumpy numbers, $m_b$, we calculate the drag coefficients
given by
\begin{align*}
        C_D = \frac{F_D}{U^2 R},
\end{align*}
where $F_D$ is the drag force, $U$ is the characteristic
velocity ($U$=1), and the maximum radius $R$ is the characteristic
length($R$=0.55).

Drag reduction via altering the no-slip condition, often achieved
through microgeometries or micro-patterning which trap the fluid, is
a topic that has received attention over the years. This reduction,
known as the roller bearing effect, is due to the formation of
embedded vortices within the bumps(or triangular cavities).

\begin{figure}[H]
        \centering
        \subfigure[]
        {
            \includegraphics[width=0.6\textwidth]{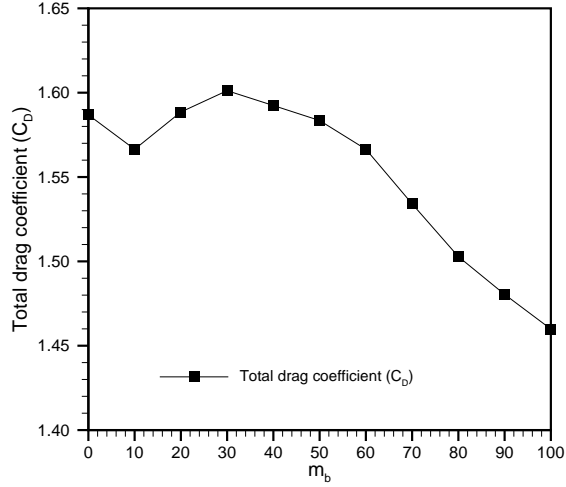}
        }\qquad
        \subfigure[]
        {
            \includegraphics[width=0.6\textwidth]{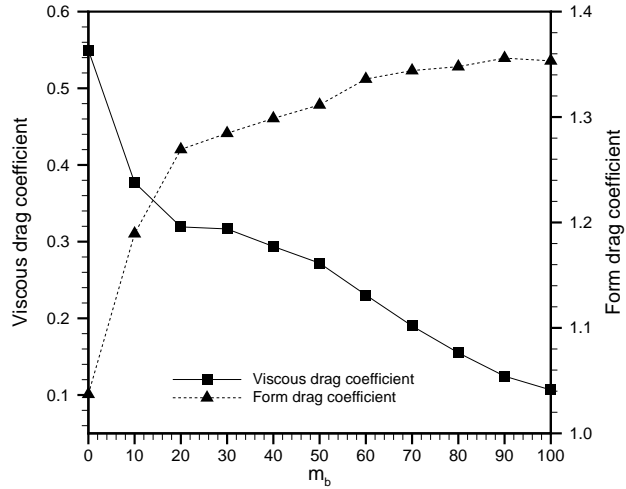}
        }
        \caption{Dependence of drag coefficients on the number of bumps with $Re=40$: (a) total Drag; (b) viscous drag and form drag.
        Three values on the vertical axis are the drag coefficients calculated in the reference case where the regular smooth circle of the radius $R=0.55$ is employed.}\label{fig:bumpy Cd}
    \end{figure}

Several interesting trends are found in our results. First of all,
the changes of drag coefficients are noticeable as the number of
bumps increases. Drag data was then compared to that obtained over a
regular smooth circular cylinder. The result of drag
coefficients($C_D$) are presented in Fig.~\ref{fig:bumpy Cd}(a). The
coefficient of total drag increases until $m_b$ reaches 30. Then,
after $m_b$=30, the total drag coefficient decreases back. Exceeding
$m_b$ = 100, the coefficient becomes far less than the reference
value that is the case of the regular smooth circle of the radius
$R$ and is plotted on the vertical axis in the figure. Results show
that for $m_b > 30$ an appreciable drag reduction of greater than
7.9\% is obtained.

We also show the dependence of form drag and viscous drag on the
bump number in Fig.~\ref{fig:bumpy Cd}(b). It is as well of
interesting note that the viscous drag coefficient decreases as the
number of bumps increases. It keeps decreasing down to less than
80\% lesser value than that of the regular circle ($m_b = 0$).

\begin{figure}[H]
        \centering
        \subfigure[$m_b=10$]
        {
            \includegraphics[width=0.4\textwidth]{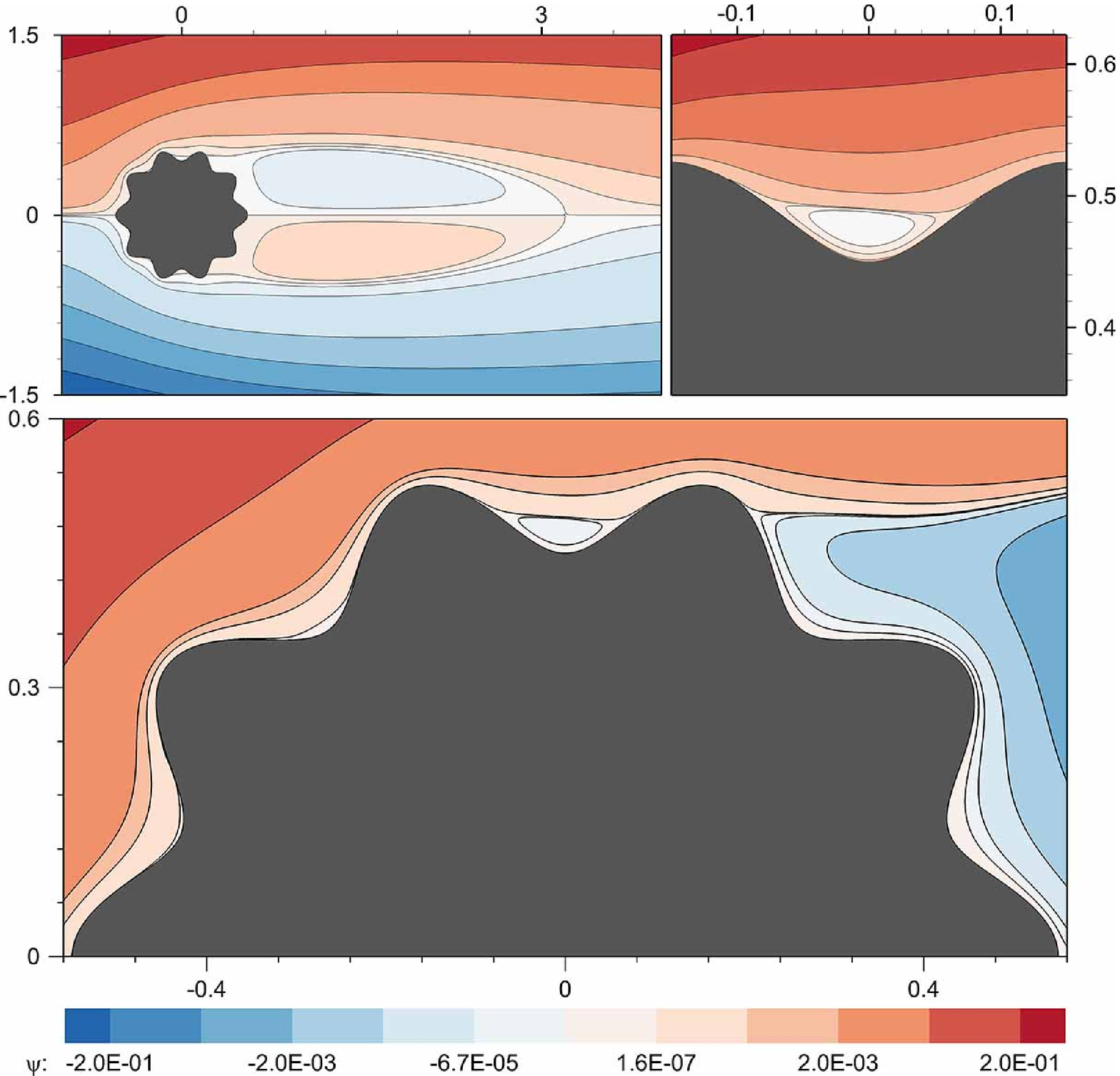}
        }\qquad
        \subfigure[$m_b=30$]
        {
            \includegraphics[width=0.4\textwidth]{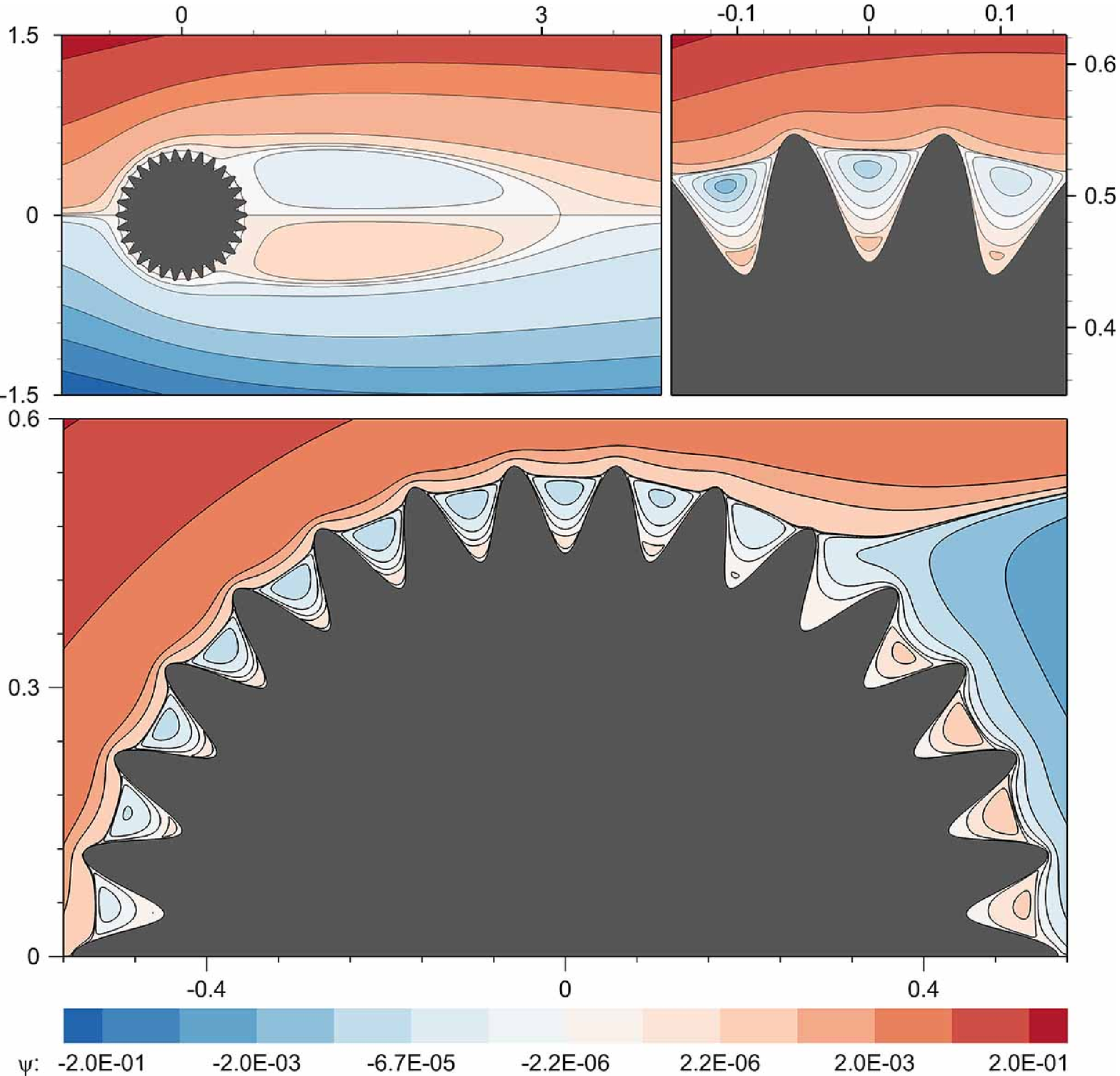}
        }\\
        \subfigure[$m_b=50$]
        {
            \includegraphics[width=0.4\textwidth]{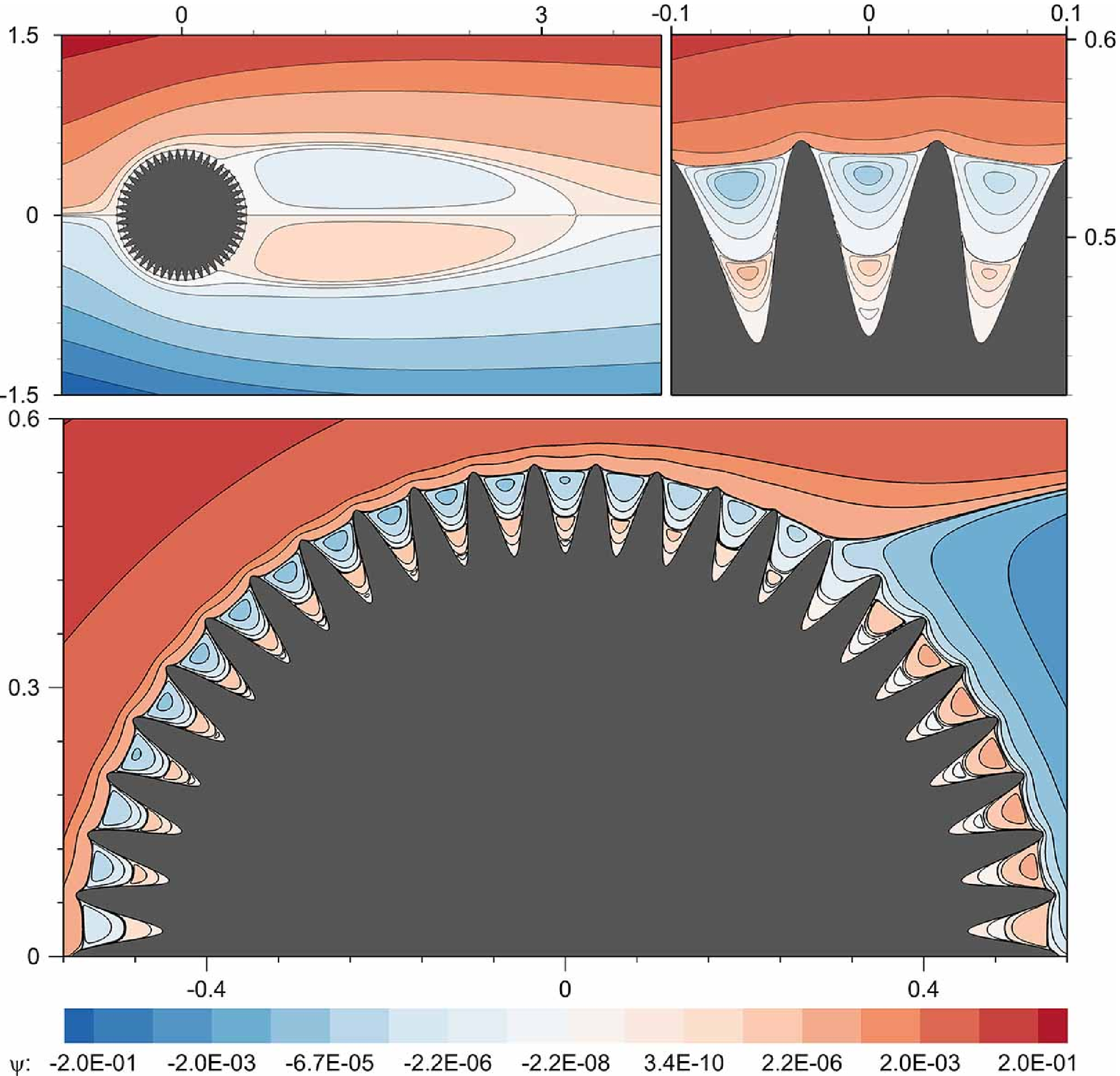}
        }\qquad
        \subfigure[$m_b=70$]
        {
            \includegraphics[width=0.4\textwidth]{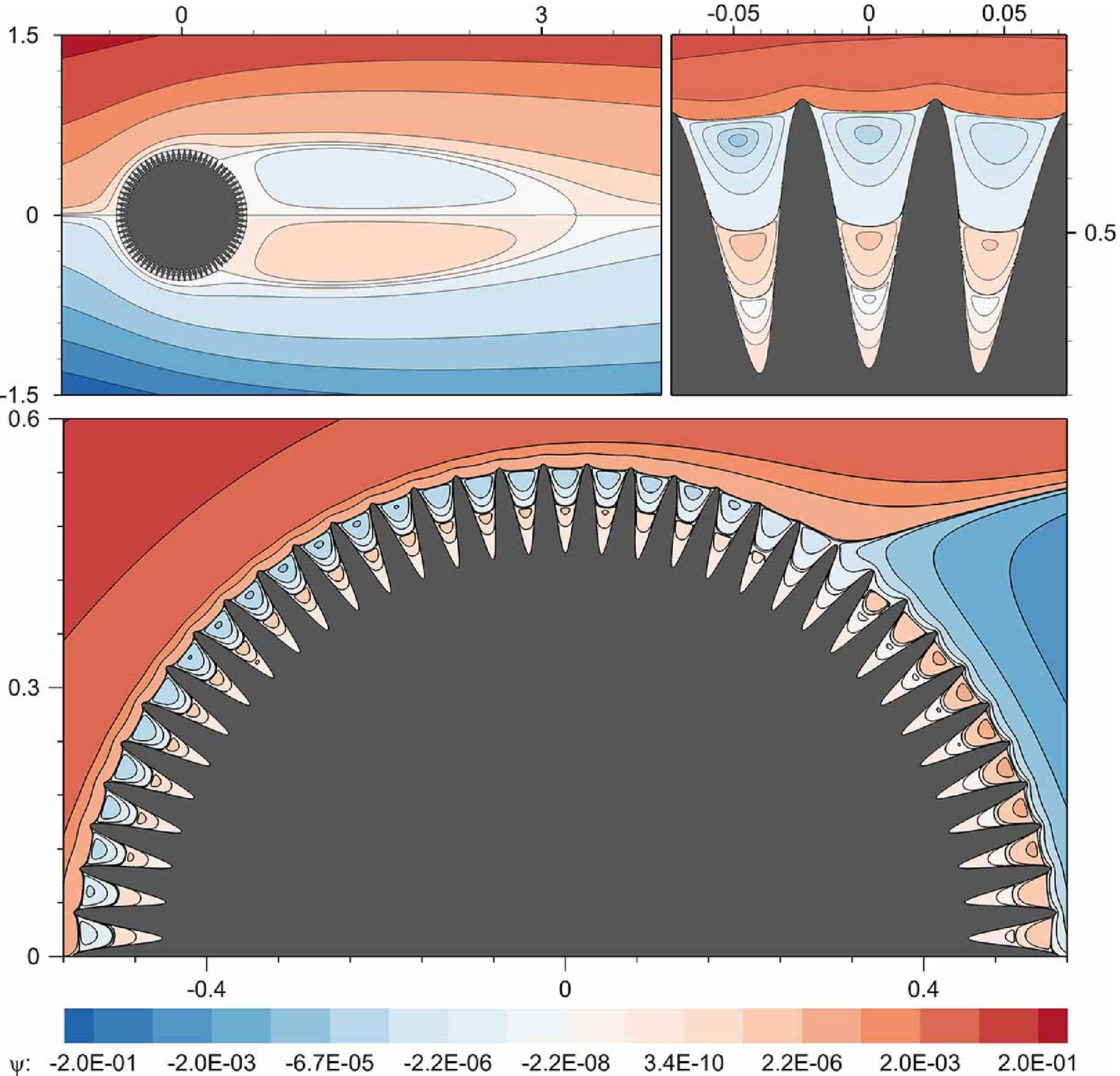}
        }\\
        \subfigure[$m_b=90$]
        {
            \includegraphics[width=0.4\textwidth]{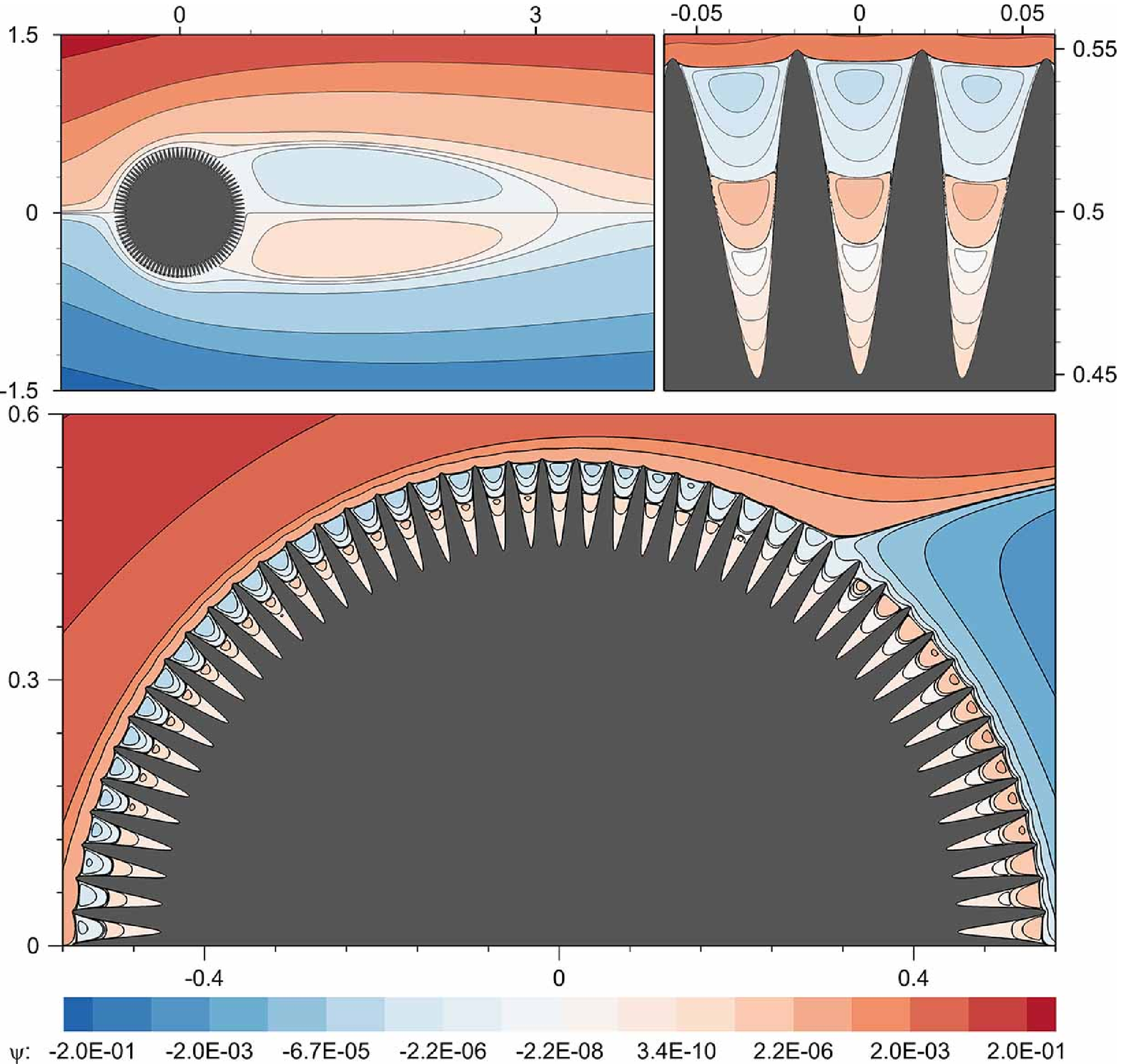}
        }\qquad
        \subfigure[$m_b=100$]
        {
            \includegraphics[width=0.4\textwidth]{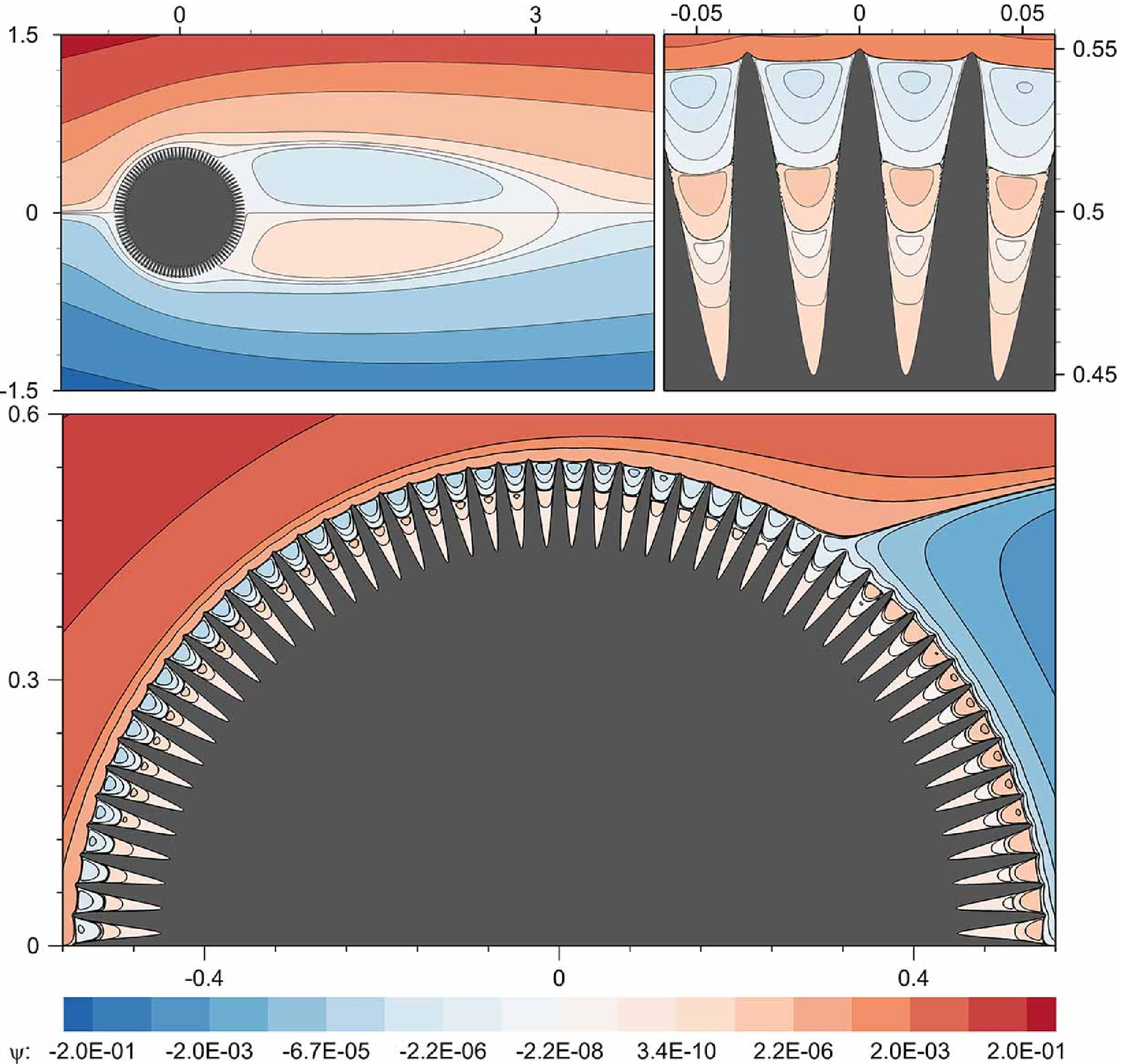}
        }
        \caption{Flow pattern around various bumpy circular cylinders with $Re=40$.
        The stream function contours around the body  for six different numbers of bumps, $m_b$ = 10, 30, 50, 70, 90, and 100.}\label{fig:bumpy stream functions}
    \end{figure}

\begin{figure}[H]
        \centering
        \subfigure[$m_b=10$]
        {
            \includegraphics[width=0.4\textwidth]{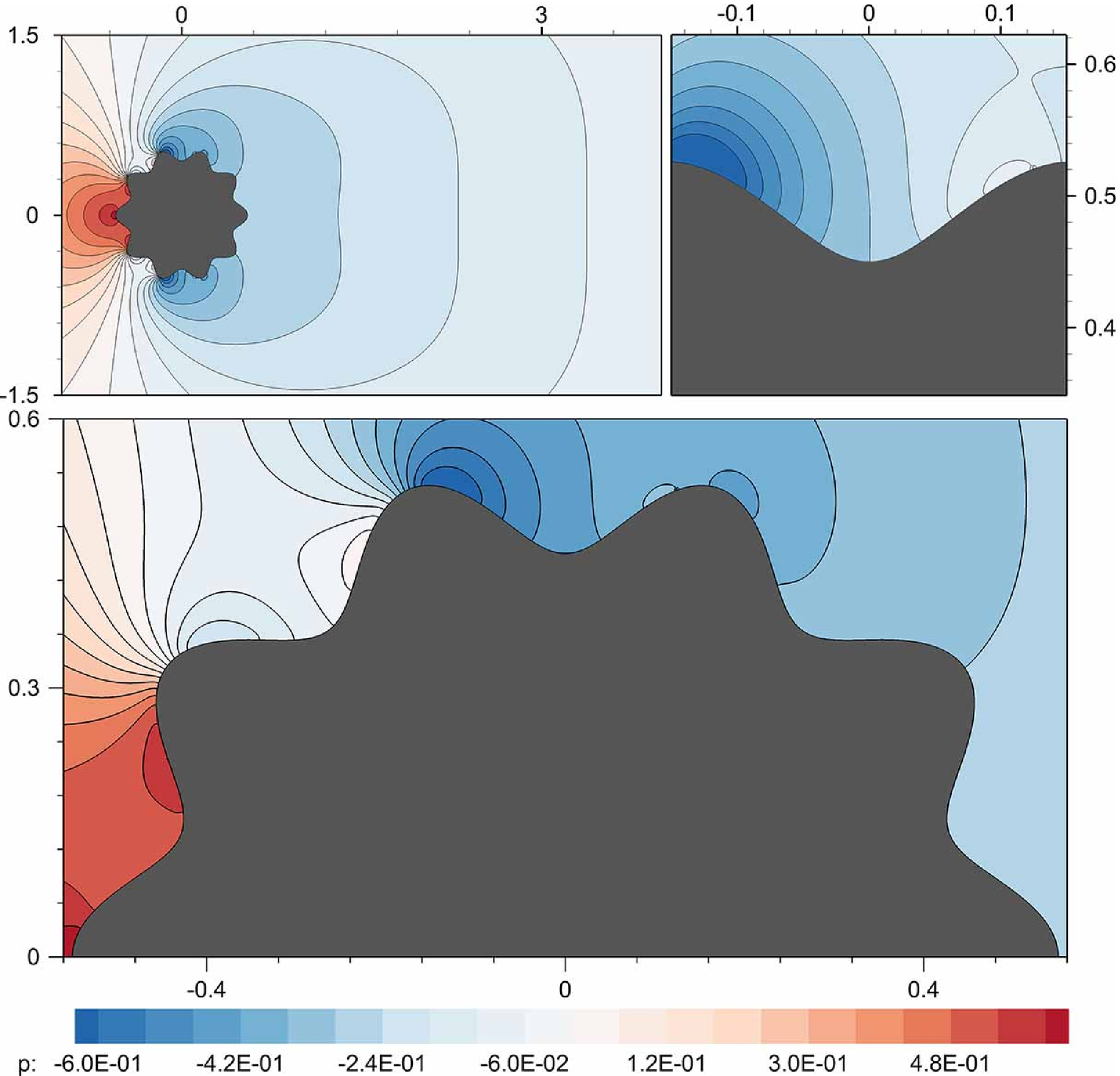}
        }\qquad
        \subfigure[$m_b=30$]
        {
            \includegraphics[width=0.4\textwidth]{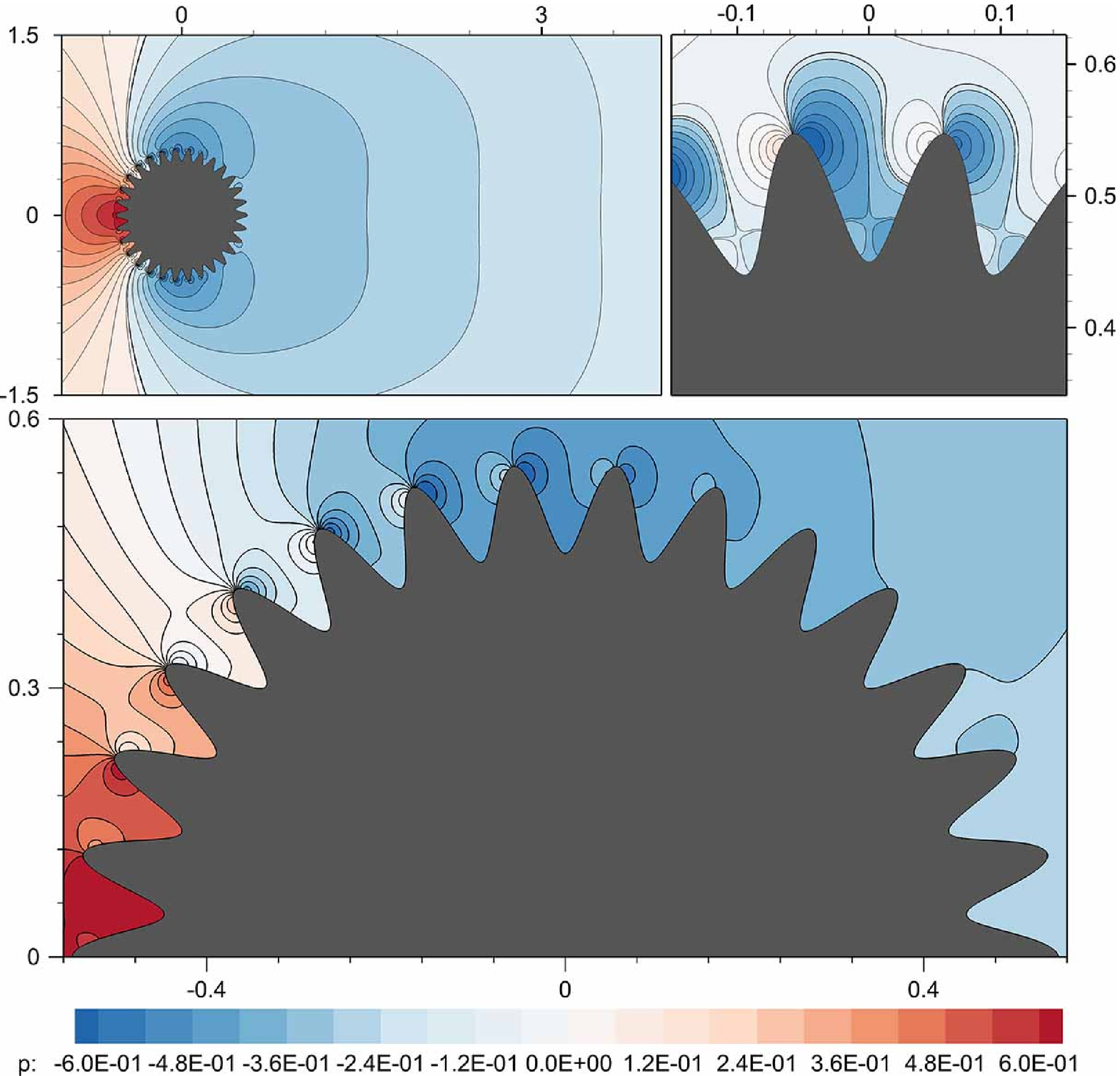}
        }\\
        \subfigure[$m_b=50$]
        {
            \includegraphics[width=0.4\textwidth]{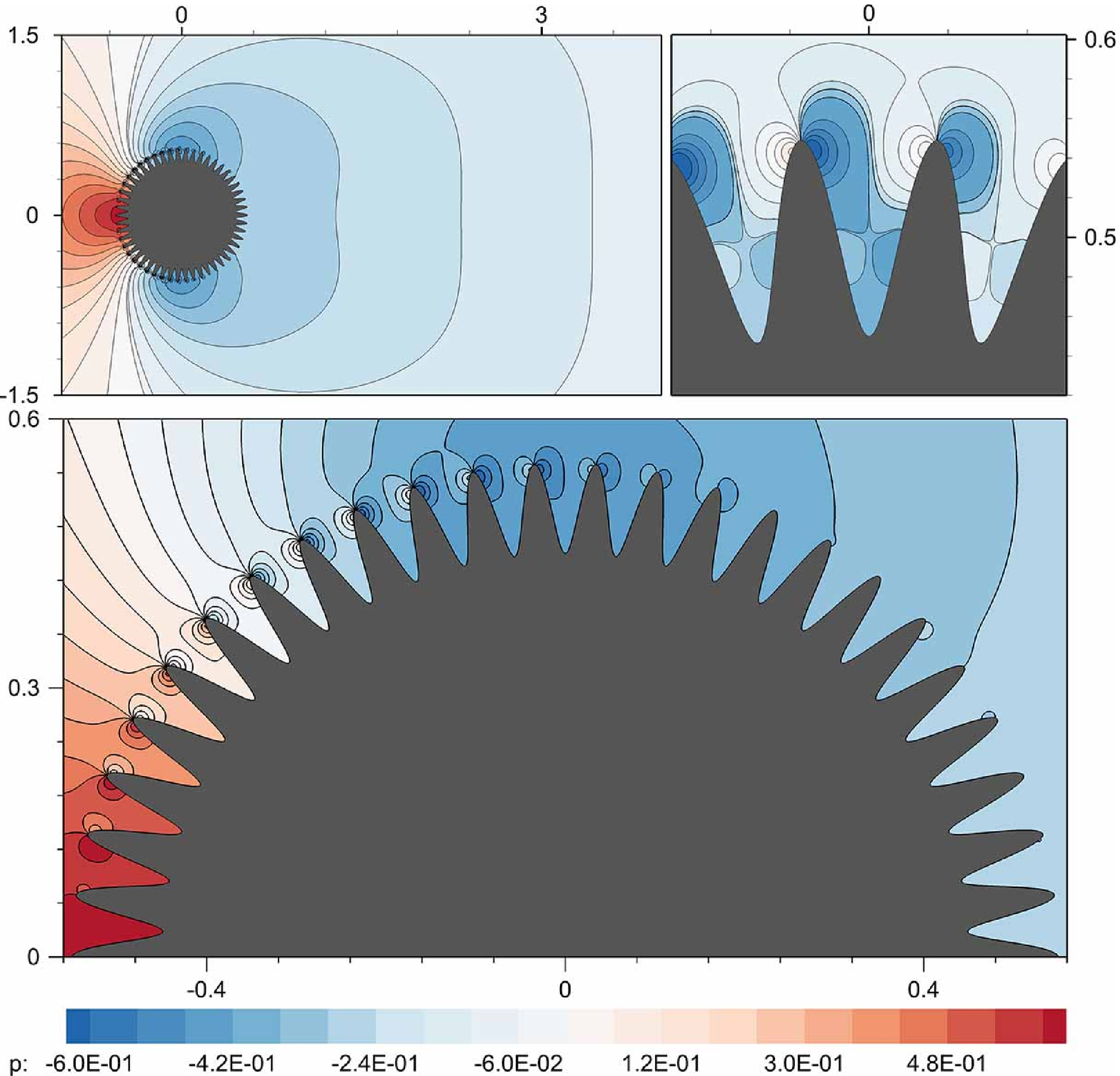}
        }\qquad
        \subfigure[$m_b=70$]
        {
            \includegraphics[width=0.4\textwidth]{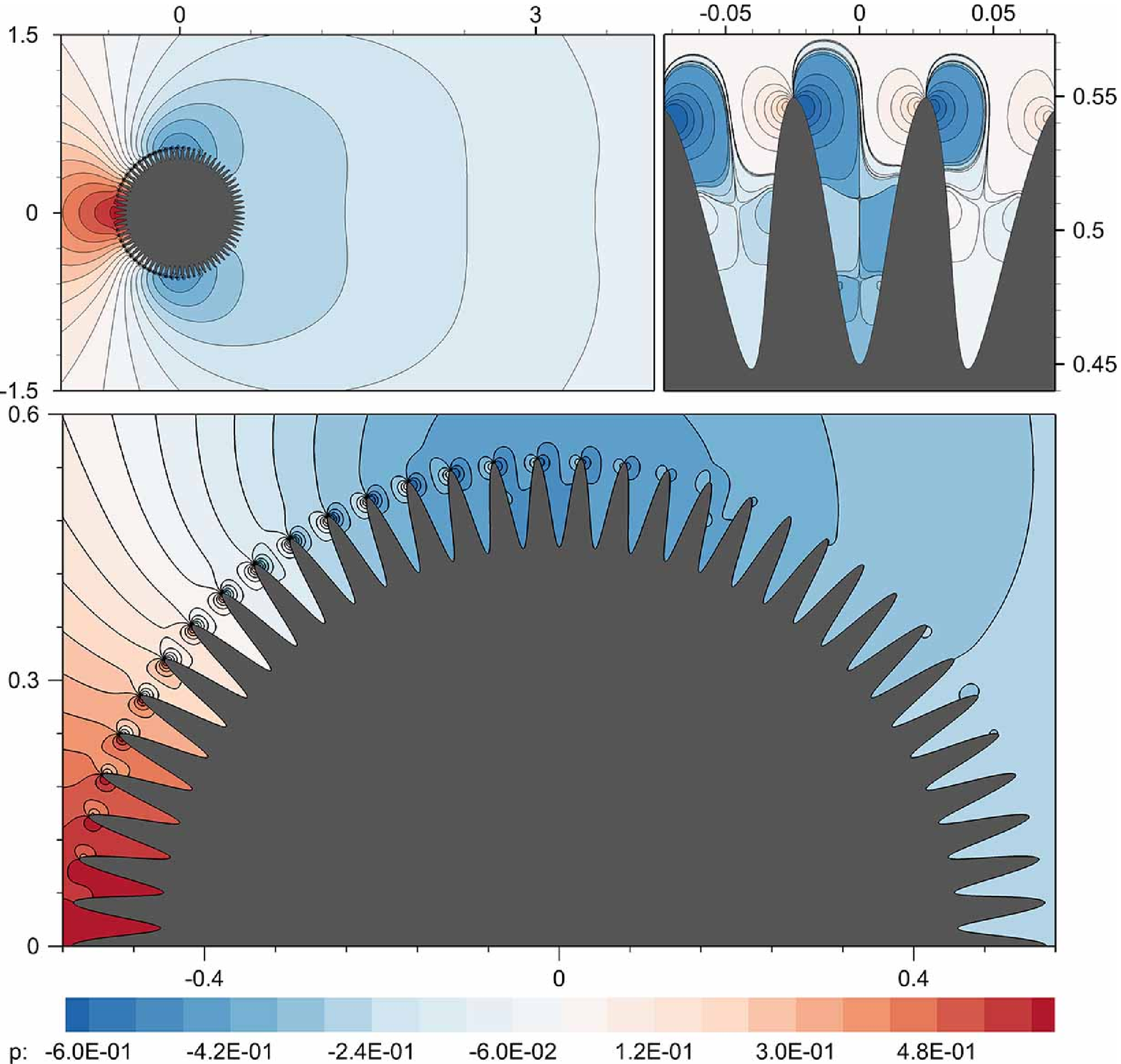}
        }\\
        \subfigure[$m_b=90$]
        {
            \includegraphics[width=0.4\textwidth]{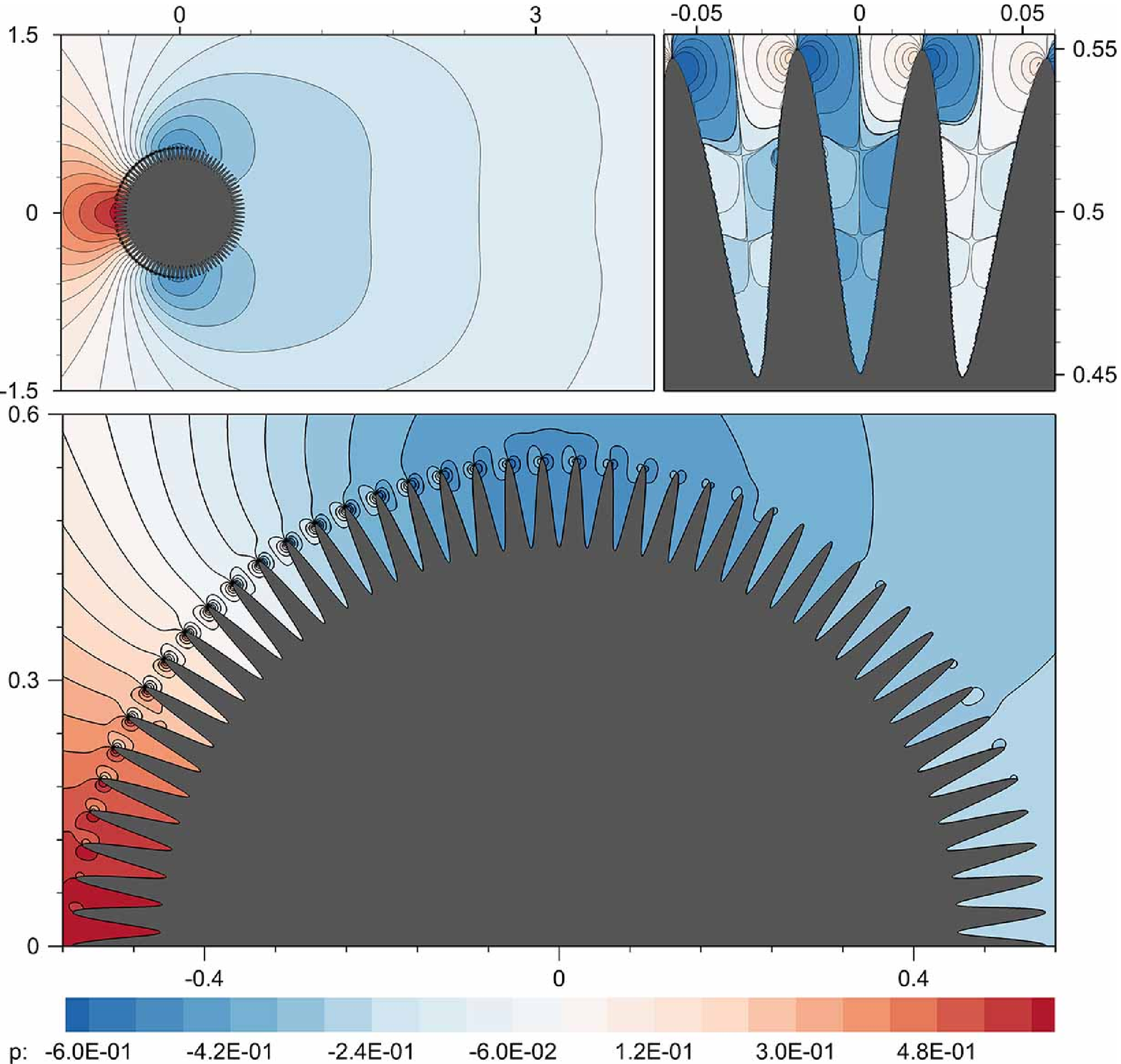}
        }\qquad
        \subfigure[$m_b=100$]
        {
            \includegraphics[width=0.4\textwidth]{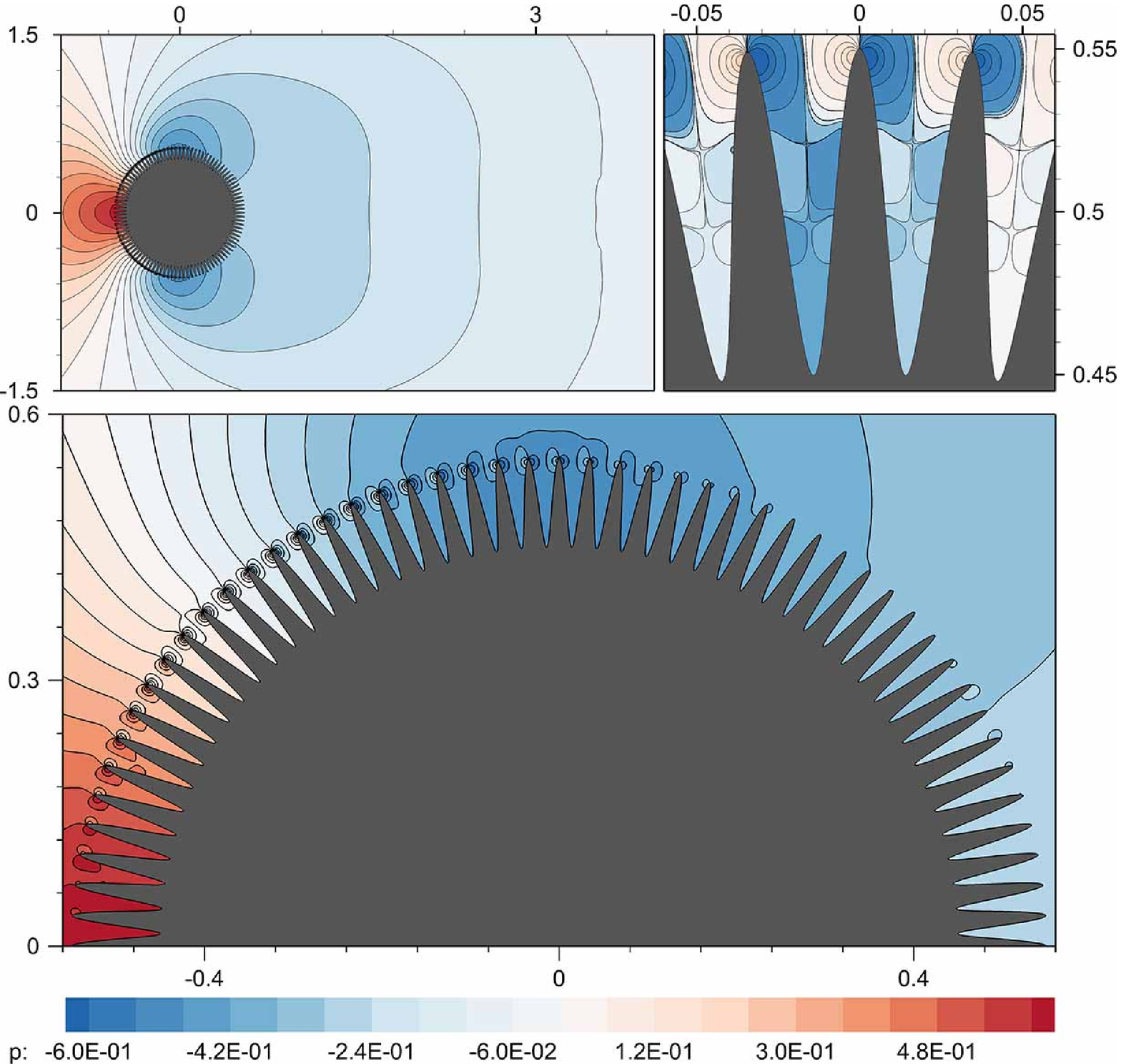}
        }
        \caption{Flow pattern around various bumpy circular cylinders with $Re=40$.
        The pressure contours around the body  for six different numbers of bumps, $m_b$ = 10, 30, 50, 70, 90, and 100.}\label{fig:bumpy pressure}
    \end{figure}

\begin{figure}[H]
        \centering
        \subfigure[$m_b=10$]
        {
            \includegraphics[width=0.4\textwidth]{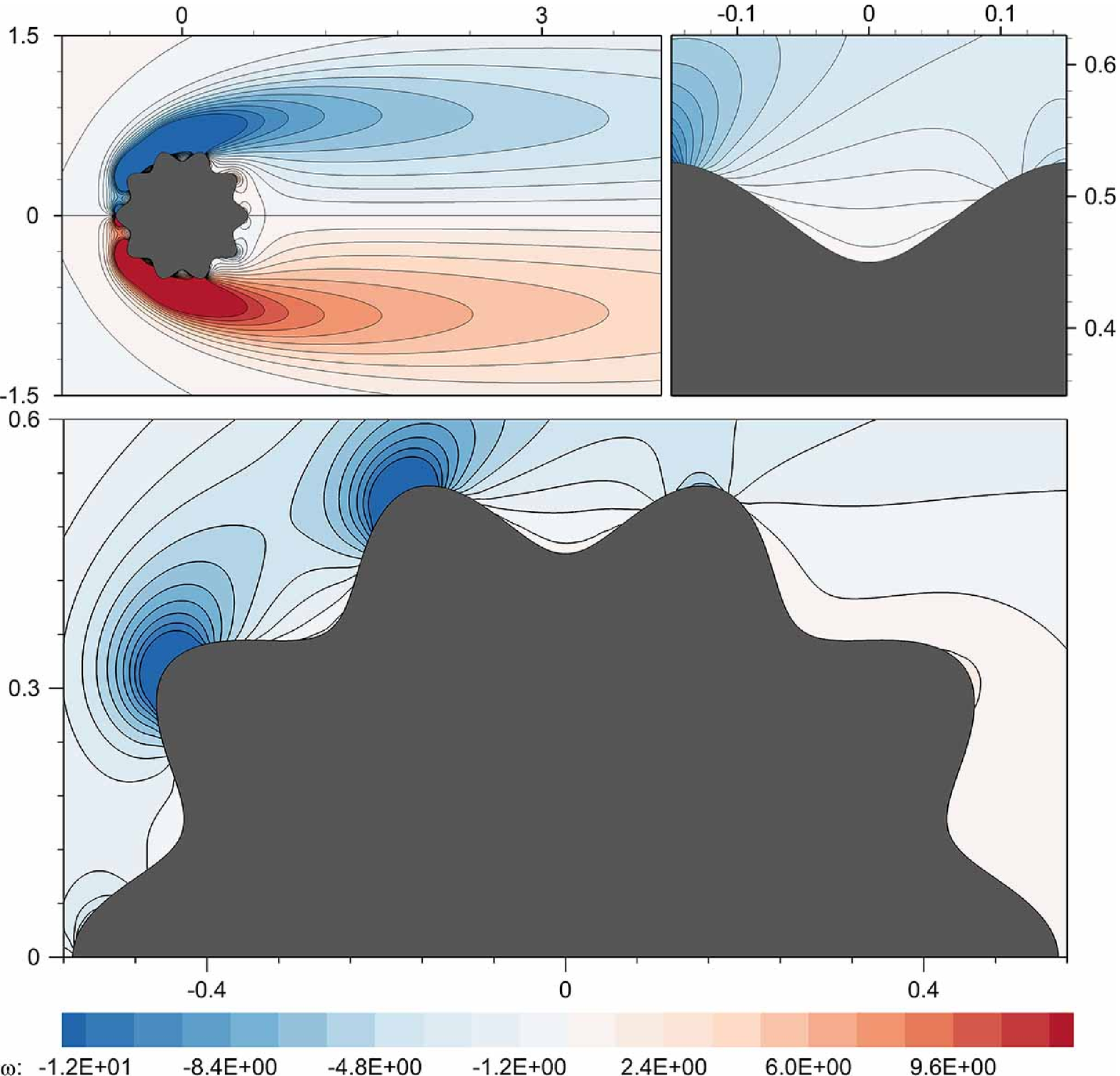}
        }\qquad
        \subfigure[$m_b=30$]
        {
            \includegraphics[width=0.4\textwidth]{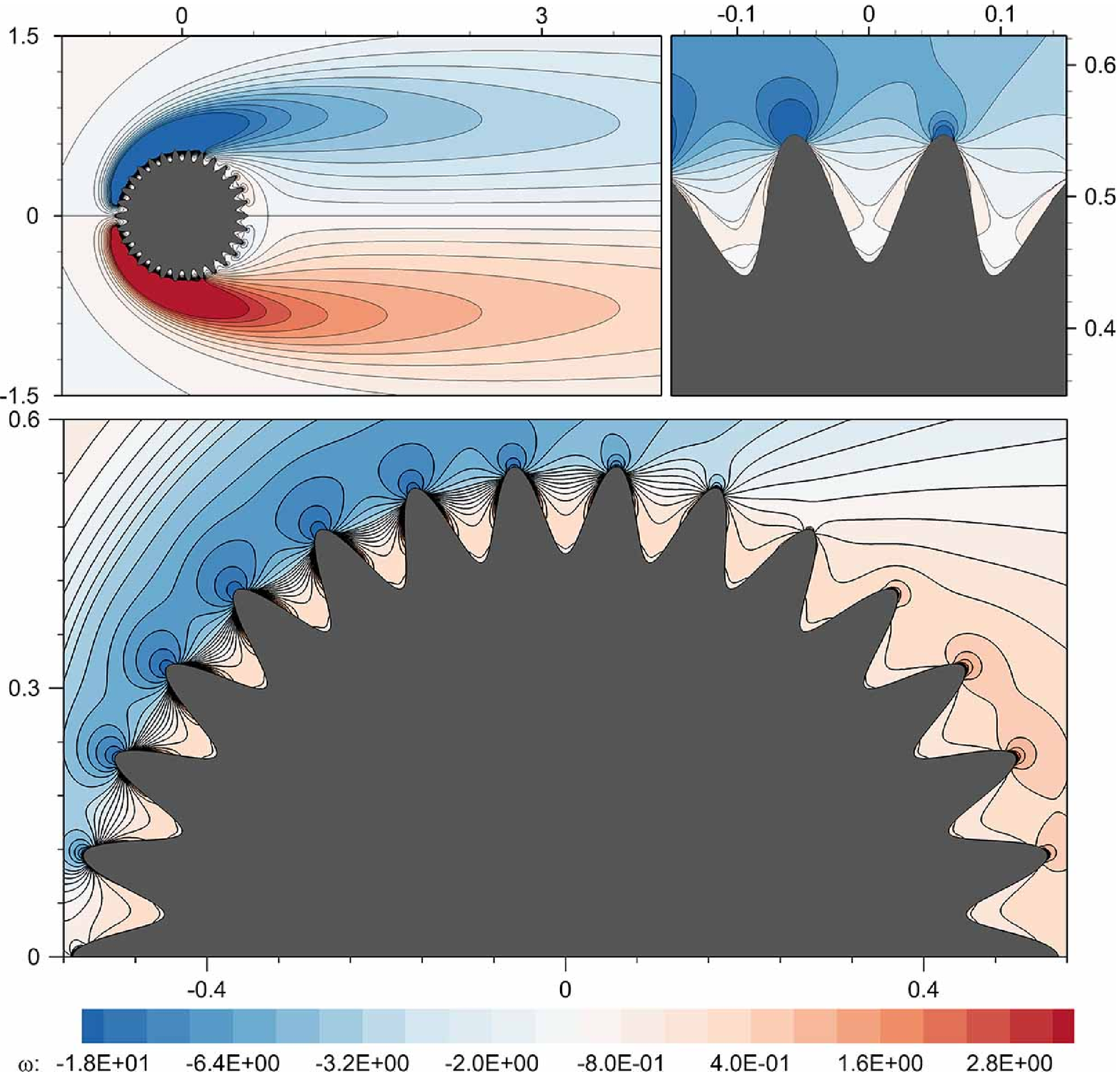}
        }\\
        \subfigure[$m_b=50$]
        {
            \includegraphics[width=0.4\textwidth]{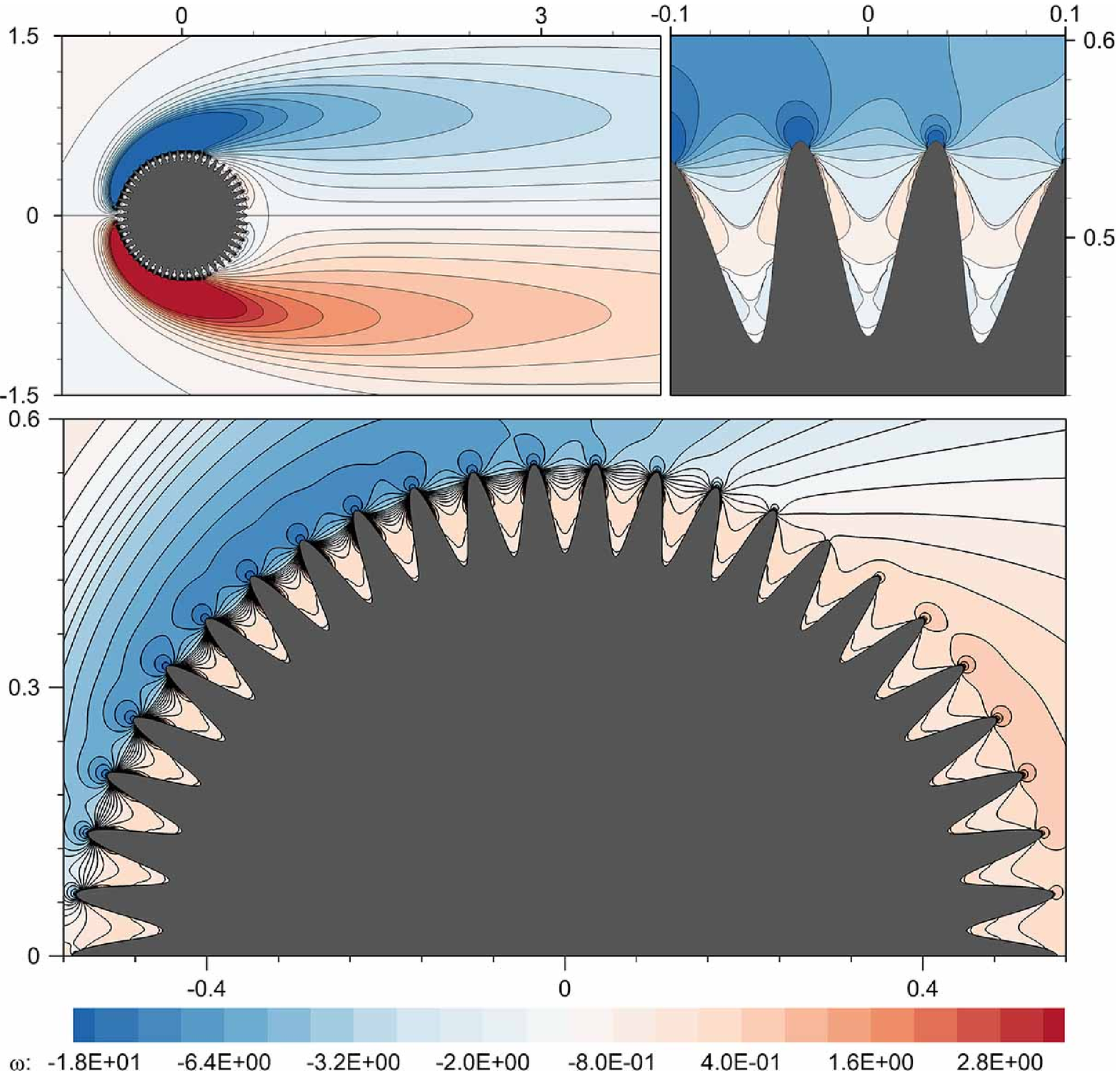}
        }\qquad
        \subfigure[$m_b=70$]
        {
            \includegraphics[width=0.4\textwidth]{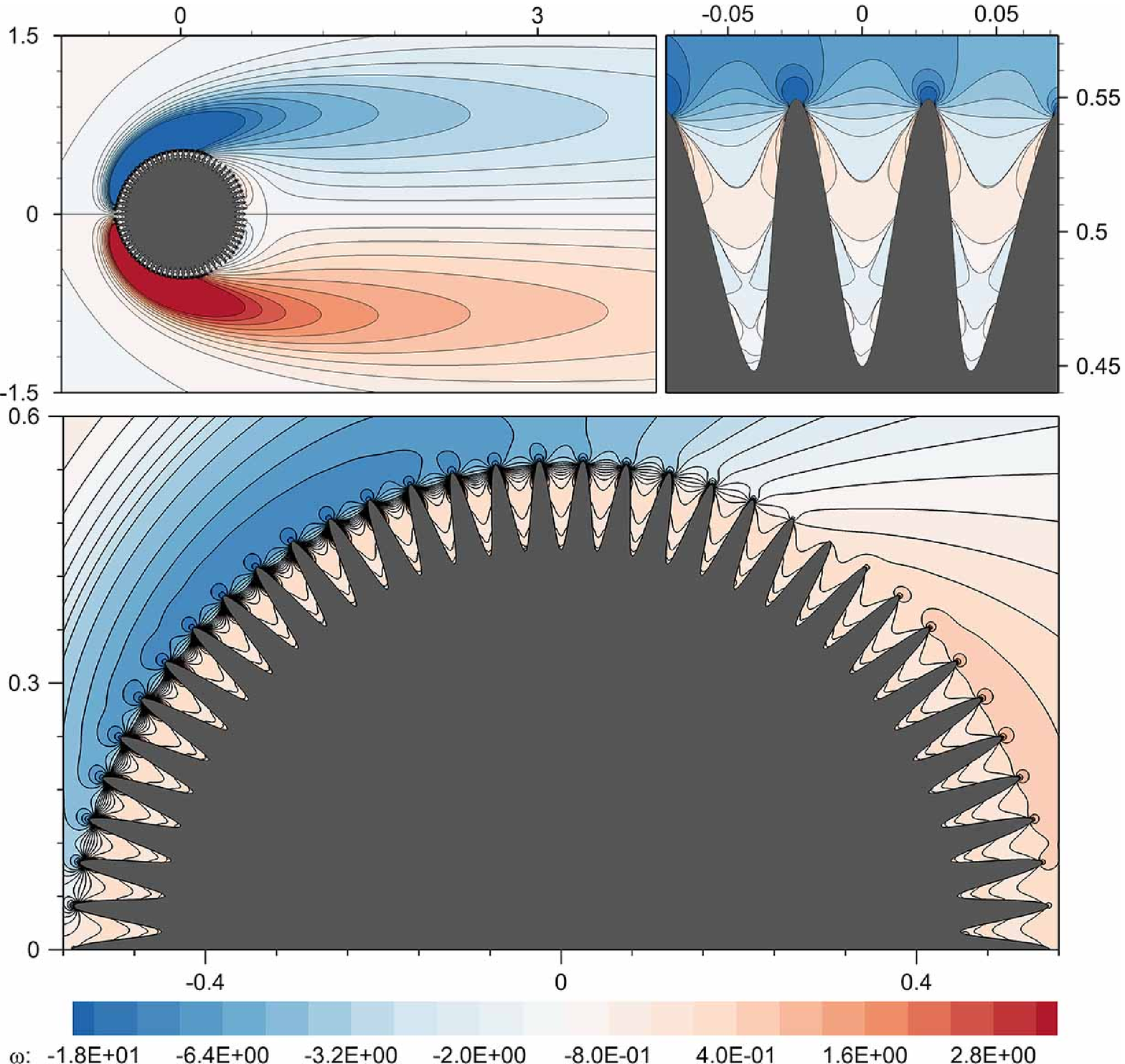}
        }\\
        \subfigure[$m_b=90$]
        {
            \includegraphics[width=0.4\textwidth]{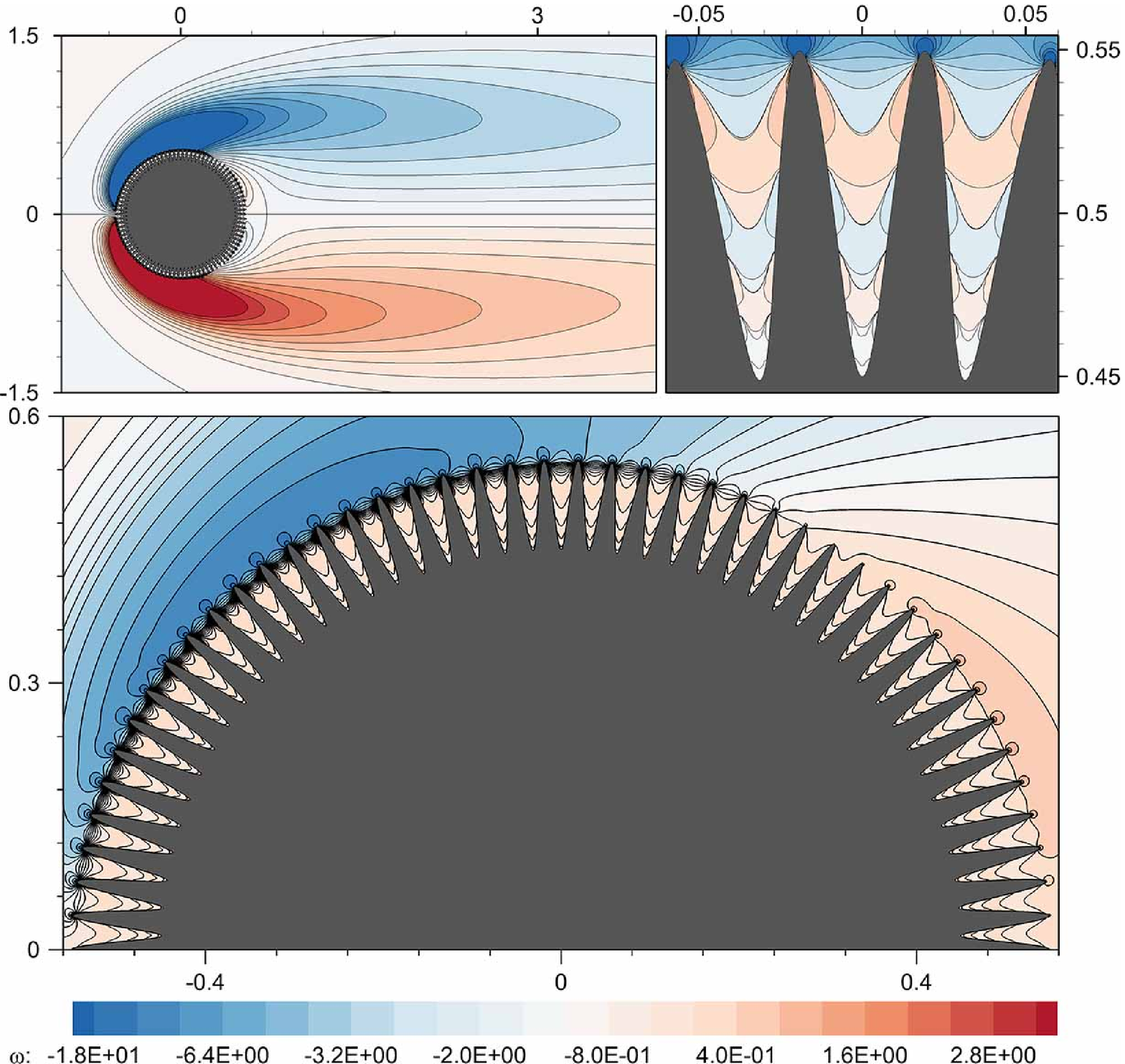}
        }\qquad
        \subfigure[$m_b=100$]
        {
            \includegraphics[width=0.4\textwidth]{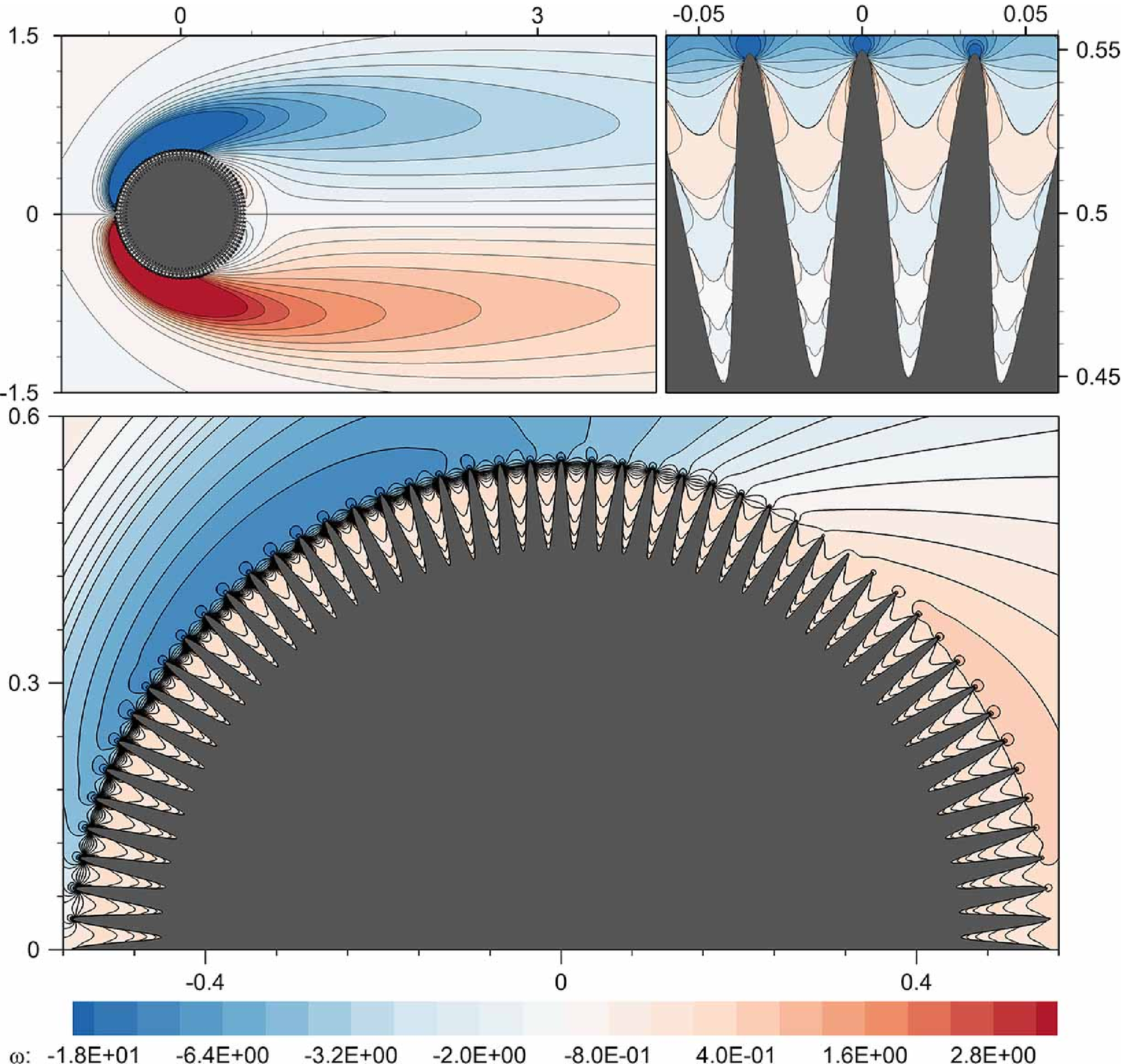}
        }
        \caption{Flow pattern around various bumpy circular cylinders with $Re=40$.
        The vorticity contours around the body  for six different numbers of bumps, $m_b$ = 10, 30, 50, 70, 90, and 100.}\label{fig:bumpy vorticity}
    \end{figure}

\begin{figure}[H]
        \centering
        \includegraphics[width=0.75\textwidth]{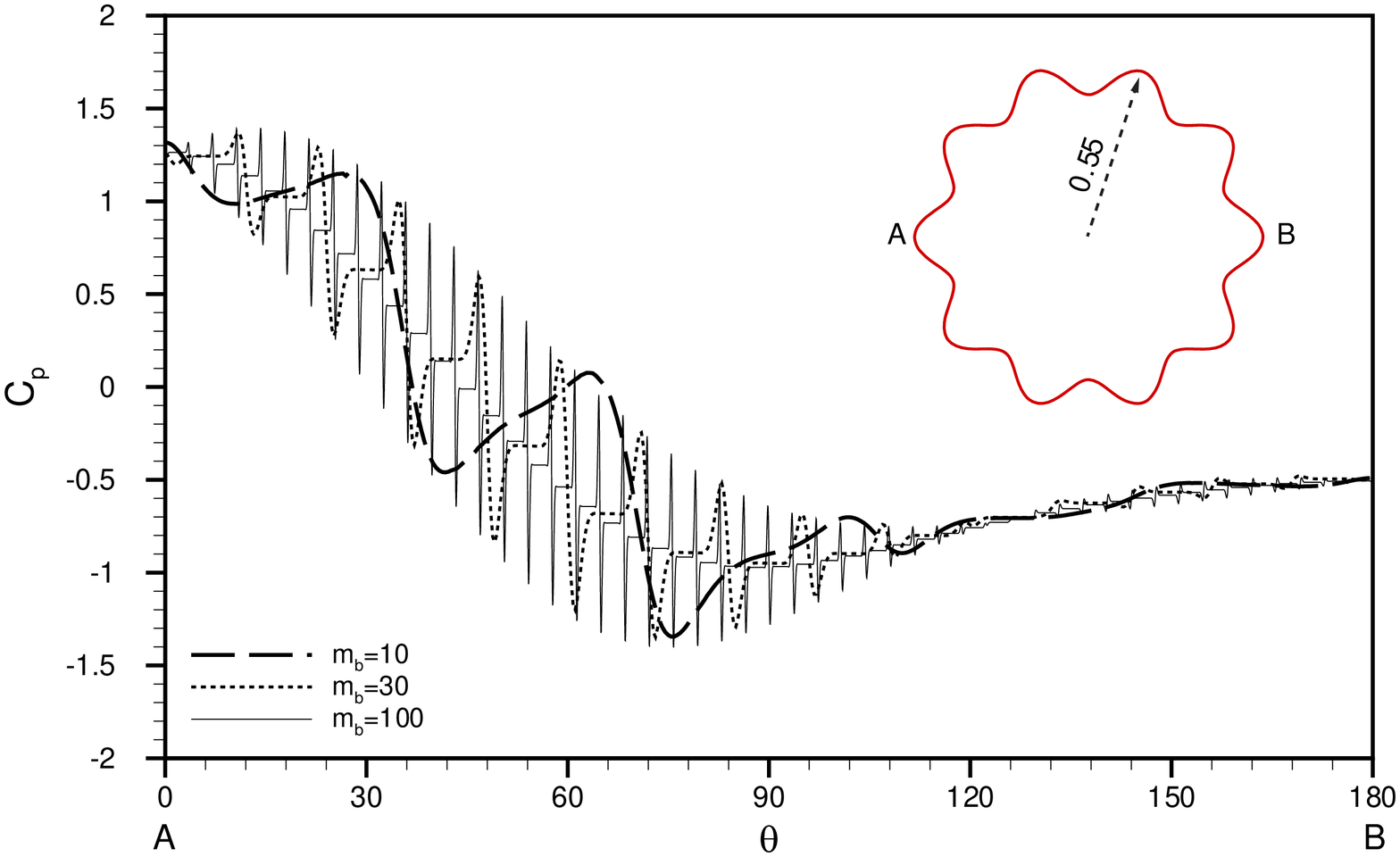}
        \caption{Pressure distributions on the bumpy surface with $Re=40$.
        A and B indicate the front and end noses.}\label{fig:bumpy Cp}
    \end{figure}

The flow pattern is illustrated in Figs.~\ref{fig:bumpy stream
functions},~\ref{fig:bumpy pressure}, and~\ref{fig:bumpy vorticity}
for the cases of $m_b$ = 10, 30, 50, 70, 90, 100. In
figure.~\ref{fig:bumpy stream functions}, the lines do not depict
equal intervals, which is intended to visualize the weak eddies
arising between bumps. These eddies may implicate the energy
transfer from flow and, as a result, the drag on the body.
In order to further convince ourselves, we show the values of
pressure($C_p$) profiles along the surfaces are given in
Fig.~\ref{fig:bumpy Cp} for three different numbers of bumps, $m_b$
=10, 30, and 100. The pressure profiles vary up and down along the
circumference, being consistent with the each bumpy shape
considered.

The bumpy circular cylinders are able to reduce the total drag in a
viscous flow via an embedded vortex inside of the bumps that imposes
a slip condition, versus a no-slip condition, where the bottom wall
in a viscous flow would normally be. Increased aspect ratio and
reduced gap height lead to better drag reduction potential at
Reynolds number,$Re=$40.

\section{Conclusion}
In this study, a new VIP method is presented by introducing the
accurate virtual interpolation point scheme as well as the virtual
local stencil approach. The present method is based on the concept
of point collocation on a virtual staggered structure together with
a fractional step method.

The virtual staggered structure consists of the virtual
interpolation points and the virtual local stencil. The use of the
virtual staggered structure arrangement, which stores all the
variables at the same physical location and employs only one set of
nodes using virtual interpolation points, reduces the numerical difficulty is caused by geometrical
complexity.

In the VIP method, the choice of an accurate interpolation scheme
satisfying the spatial approximation in the complex domain is
important because there is the virtual staggered structure for
computation of the velocities and pressure since there is no
explicit staggered structure for stability. In our proposed method,
the high order derivative approximations for constructing node-wise
difference equations are easily obtained.

Several different flow problems (decaying vortices, lid-driven
cavity, triangular cavity, flow over a circular cylinder and a bumpy
cylinder) are simulated using the virtual interpolation point method
proposed in this study. The simulation results with both the the
accurate virtual interpolation point scheme and the virtual local
stencil approach agree very well with the previous numerical and
experimental results, indicating the validity and accuracy of the
present VIP method.

\bibliographystyle{model1-num-names}

\bibliography{ref}

\end{document}